%% file: PABM_ArXiv.tex
\documentclass[11pt]{article}
 \usepackage{amsmath}
\usepackage{amssymb}
\usepackage{amsthm}

\usepackage{graphicx}
\usepackage{color}
\usepackage{multicol}
\usepackage{mathrsfs}
\usepackage{float}
\usepackage[caption = false]{subfig} 
\usepackage{dsfont}
\usepackage{times}
\usepackage{commath} 
\usepackage[linesnumbered,ruled]{algorithm2e}

\begin{document}

%%%%%%%%%%%%%%%%%%%%%%%%%%%%%%%%%%%%%%%%%%%%%%%%%%%%%%%%%%%%%%%%%%%%%%%%%%%%%%%%%%%%%%%%%%%

\input{Def1.tex}

%%%%%%%%%%%%%%%%%%%%%%%%%%%%%%%%%%%%%%%%%%%%%%%%%%%%%%%%%%%%%%%%%%%%%%%%%%%%%%%%%%%%%%%%%%%

\title{  { Estimation and Clustering in Popularity Adjusted Stochastic Block Model }}

\author{{    Majid Noroozi, Ramchandra Rimal and Marianna Pensky}   \\
         Department of Mathematics,
         University of Central Florida 
 }

\date{}

\bibliographystyle{abbrv} %{plain}
\maketitle

\begin{abstract}
The paper considers the Popularity Adjusted Block model (PABM) introduced by Sengupta and Chen (2018).
We argue that the main appeal of the PABM is the flexibility of the spectral properties of the 
graph which makes the PABM an attractive choice for modeling networks that appear in biological sciences. 
We expand the theory of PABM to the case of an arbitrary number of communities 
which possibly grows with a number of nodes in the network and is not assumed to be known. 
We produce the estimators of the probability matrix and the community structure and 
provide non-asymptotic upper bounds for the estimation and the clustering errors.
We use  the Sparse Subspace Clustering (SSC) approach to partition the network into communities,
the approach that, to the best of our knowledge,   has not been used for clustering network data.
The theory is supplemented by a simulation study. In addition, we show advantages of the PABM  
for modeling a butterfly similarity network and a human brain functional network.

\vspace{2mm} 


\vspace{2mm}{\bf AMS (2000) Subject Classification}: {Primary:  62F12, 62H30.   Secondary:  05C80  }
\end{abstract}

\section{Introduction}
\label{sec:intro}
\setcounter{equation}{0}

Statistical network analysis has become a major field of research, with applications as diverse as 
sociology, biology, genetics, ecology, information technology  to name a few. 
An overview of statistical modeling of random graphs can be found in, 
e.g.,   \cite{MAL-005} and \cite{Kolaczyk:2009:SAN:1593430}.

Consider an undirected network with $n$ nodes and no self-loops and multiple edges. 
Let $A \in \{0,1\}^{n\times n}$ be the symmetric adjacency  matrix of the network with
$A_{i,j}=1$ if there is a connection between nodes $i$ and $j$,  and $A_{i,j}=0$
otherwise. We assume that 
\be \label{eq:A_Bern_P}
 A_{i,j}\sim \mbox{Bernoulli}(P_{i,j}),  \quad 1 \leq i \leq j \leq n,
\ee
where  $A_{i,j}$  are conditionally independent 
given $P_{i,j}$ and  $A_{i,j} = A_{j,i}$, $P_{i,j} = P_{j,i}$ for $i>j$.

The block models assume that each node  in the network belongs to one of $K$ distinct blocks or communities $\calN_k$, $k=1, \cdots, K$. 
Let $c$ denote the  vector of community assignment,  with $c_i = k$ if the   node $i$ belongs to the   community $k$.
Then, the probability of connection between node $i \in \calN_k$ and node $j \in  \calN_l$  depends on the 
pair of blocks $(k,l)$ to which nodes $(i,j)$ belong. One can also consider a corresponding {\it membership} 
(or {\it clustering}) matrix $Z \in \{0,1\}^{n\times K} $     such that $Z_{i,k}=1$ iff $i \in \calN_k$, $i=1, \ldots, n$.

A classical random graph model for networks with community structure is the Stochastic Block Model (SBM)  that has been 
% introduced by Lorrain and White (1971) and subsequently 
studied by a number of authors (see, e.g., \cite{JMLR:v18:16-480}, \cite{Gao:2017:AOM:3122009.3153016}  among others). 
Under this model, all nodes belonging to a community are considered to be stochastically equivalent, 
in the sense that the probability of connection between nodes is completely defined by the communities to which they 
belong. 
Specifically,  under  the $K$-block SBM,  this probability is completely  determined by the community assignment for nodes $(i,j)$, so that
$P_{i,j} = B_{c_i,c_j}$  where $B_{k,l}$ is the probability of connection between communities $k$ and $l$.
In particular, any   nodes from the same community have the same degree distribution and the same expected degree.

Since the real-life networks usually contain a very small number of high-degree nodes while the rest 
of the nodes have very few connections (low degree), the SBM model fails to explain the structure of many networks that occur  in practice. 
The  Degree-Corrected Block Model (DCBM) addresses this deficiency by  allowing these probabilities to be 
multiplied by the node-dependent weights (see, e.g., \cite{Karrer2011StochasticBA},
\cite{chen2018} and \cite{zhao2012consistency} among others). Under the DCBM, the elements of matrix $P$ are modeled as 
$P_{i,j }= \theta_i B_{c_i,c_j}\theta_j$, where $\theta_i$, $i=1,\ldots,n$, are the degree parameters of the   nodes, and $B$ is the $(K \times K)$  
matrix of baseline interaction between communities. 
Identifiability of the parameters is usually ensured by a constraint of the form  $\sum_{i\in \mathcal{N}_k} \theta_i = 1$ for all    $k = 1,...,K$
(see, e.g., \cite{Karrer2011StochasticBA}).

A network feature that is closely associated with community structure is the popularity of nodes across communities 
defined as the number of edges between a specific node and a specific community. 
While the DCBM allows to correctly detect  the communities, 
and accurately fits the total degree   by enforcing the node-specific degree parameters, it 
enforces the node popularity to be uniformly proportional to the node degree. 
Hence, the DCBM fails to model node popularities in a flexible and realistic way. 
For this reason, recently,  \cite{RePEc:bla:jorssb:v:80:y:2018:i:2:p:365-386} 
introduced the Popularity Adjusted Stochastic Block Model (PABM) which models  the probability of a connection between nodes  
as a product of popularity parameters that depend on the communities to which the nodes belong as well as on the pair of nodes themselves.
In particular, in PABM 
\begin{equation} \label{eq:PABM-model}
P_{i,j}= V_{i,c_{j}} V_{j,c_{i}}, 
\end{equation}
where $V_{i,k}$, $1 \leq i \leq n$,  $1 \leq k \leq K, $ is the scaling parameter that identifies 
 popularity of node $i$ in class $k$,  and $0 \leq P_{i,j} \leq 1$ for any $i$ and $j$.
Specifically,  %Sengupta and Chen
 \cite{RePEc:bla:jorssb:v:80:y:2018:i:2:p:365-386} define the  popularity 
of  node  $i$ in     community $k$ as 
% $ \mu _{i,k} =   \sum_{j\in \mathcal{N}_k}  P_{i,j}$. 
 $ \mu _{i,k} =  \displaystyle \sum_{j\in \mathcal{N}_k}  P_{i,j}$. 
They  noted that the ratio of   popularities of 
the nodes $(i,j)\in \calN_k$ in the same community $k$  is equal to one for the SBM, is  independent of community $k$  
(a function of   $i$ and $j$ only) in DCBM but can vary between nodes and communities for the PABM, thus, allowing a more flexible  
modeling of connection probabilities. 
The authors showed that PABM generalizes both the SBM and the DCBM, suggested the quasi-maximum likelihood  type procedure for 
estimation and clustering and  demonstrated the improvement achieved through this new methodology.

The flexibility of PABM, however,  is not limited to modeling the popularity parameters of the nodes. 
In order to better understand the model, consider a rearranged version $P(Z,K)$ of matrix $P$ where 
its  first $n_{1}$ rows correspond to  nodes from class 1, 
the next $n_{2}$ rows correspond to  nodes from class 2 and the last $n_K$ rows correspond to  nodes from class $K$. 
% \textcolor{blue}{ 
Denote the $(k,l)$-th block of matrix $P(Z,K)$ by $P^{(k,l)} (Z,K)$. Since sub-matrix $P^{(k,l)} (Z,K) \in [0,1]^{n_k \times n_l}$ 
corresponds to pairs of nodes in communities $(k,l)$ respectively, one obtains from \eqref{eq:PABM-model} that 
$P^{(k,l)}_{i,j}= V_{i_k,l} V_{j_l,k}$ where   $i_k$ is the $i$-th  element in $\calN_k$ and $j_l$ is the $j$-th  element in $\calN_l$.
Thus, matrices $P^{(k,l)} (Z,K)$ are rank-one matrices with the unique singular vectors generating them. 
%}
Indeed, consider vectors $\Lam^{(k,l)}$ with elements $\Lam^{(k,l)}_i = V_{i_k,l}$, where $i=1, \ldots, n_k$ and $i_k \in \calN_k$.
Then,   equation  \eqref{eq:PABM-model} implies that  
\be  \label{eq:block_structure}
P^{(k,l)} (Z,K) =\Lambda^{(k,l)} \, [\Lambda^{(l,k)}]^{T}.
\ee 
% so that $P^{(k,l)} (Z,K)$ are rank-one matrices. 
Moreover,  it follows from  \eqref{eq:PABM-model}  and \eqref{eq:block_structure} 
that $P^{(k,l)} (Z,K) = [P^{(l,k)} (Z,K)]^T$ and that each pair of blocks $(k,l)$ involves a  unique combination 
of vectors $\Lambda^{(l,k)}$:
\bes
P(Z,K) =
  \begin{bmatrix}
    \Lambda^{(1,1)} (\Lambda^{(1,1)})^T & \Lambda^{(1,2)} (\Lambda^{(2,1)})^T & \cdots  & \Lambda^{(1,K)}(\Lambda^{(K,1)})^T \\
    \Lambda^{(2,1)} (\Lambda^{(1,2)})^T & \Lambda^{(2,2)} (\Lambda^{(2,2)})^T & \cdots  & \Lambda^{(2,K)} (\Lambda^{(K,2)})^T  \\
    \vdots & \vdots& \cdots& \vdots\\
    \Lambda^{(K,1)} (\Lambda^{(1,K)})^T & \Lambda^{(K,2)} (\Lambda^{(2,K)})^T & \cdots  & \Lambda^{(K,K)}(\Lambda^{(K,K)})^T \\
  \end{bmatrix}
\ees
where
\be \label{eq:Lambda}
\Lambda=
  \begin{bmatrix}
    \Lambda^{(1,1)} & \Lambda^{(1,2)} & \cdots  & \Lambda^{(1,K)}\\
    \Lambda^{(2,1)} & \Lambda^{(2,2)} & \cdots  & \Lambda^{(2,K)}  \\
    \vdots & \vdots& \cdots& \vdots\\
    \Lambda^{(K,1)} & \Lambda^{(K,2)} & \cdots  & \Lambda^{(K,K)} \\
  \end{bmatrix}
\ee

The latter implies that matrix $P(Z,K)$ is   formed by arbitrary rank one blocks 
and hence $\rank(P(Z,K)) =  \rank(P)$   can take any value between $K$ and $K^2$. 
In comparison, all other block models restrict the rank of $P$ to 
be exactly $K$. This is true not only for the  SBM and DCBM discussed above but also for their 
generalizations such as the Mixed Membership  models (MMM) (see, e.g., \cite{Airoldi:2008:MMS:1390681.1442798}  
and \cite{doi:10.1080/10618600.2016.1237362}) and the Degree Corrected Mixed Membership (DCMM) 
(see, e.g., \cite{2017arXiv170807852J}).  While the MMM and the DCMM allows more diverse structures of rank $K$ matrices 
(those matrices have to be just a product  of two rank $K$ matrices with nonnegative components while 
the PABM requires a combination of $K^2$ rank one matrices),  meaningful fitting of the MMM or DCMM relies on a variety of 
conditions (one needs to have pure nodes in the network and some identifiability conditions need to be satisfied).
In addition, while the MMM and DCMM are extremely useful for analysis of social and society-related networks 
such as publications networks, they may not be appropriate in some other applications where each node can belong to one and only one class.
The butterfly similarity network studied in this paper provides an example of such application.

In general, the  flexibility makes the PABM  an attractive choice for modeling networks that appear in biological sciences, 
especially in the situations where    memberships in multiple communities are not allowed. 
Indeed, while social networks exhibit assortative behavior due to the human tendency of forming strong associations, 
the biological networks tend to be more diverse.  
% For this reason, PABM tends to be a useful tool for modeling such networks.

However, while the PABM model is extremely valuable, the statistical inference in   
\cite{RePEc:bla:jorssb:v:80:y:2018:i:2:p:365-386} has been incomplete. In particular, the authors considered only the case 
of a small finite number of communities $K$; they  provided only asymptotic consistency results as $n \to \infty$ 
without any error bounds for finite values of  $n$; their NP-hard clustering procedure  was tailored to the case of a small  $K$.
In addition, the relaxation of this NP-hard procedure seems to be operational only in the case of $K=2$ since 
all simulations and real data examples in \cite{RePEc:bla:jorssb:v:80:y:2018:i:2:p:365-386} only tackled the case of $K=2$.

The purpose of the present paper is to  address some of those deficiencies and to advance the theory of the PABM. 
Specifically,  the main merit of our paper lies in the fact that  we recognize  that the probability matrix 
of the PABM is formed by a unique  collection of rank one matrices. This useful property has not been detected by 
 \cite{RePEc:bla:jorssb:v:80:y:2018:i:2:p:365-386} who   worked in terms of 
the Poisson likelihood and the Poisson  likelihood modularity maximizations.  This observation 
on the structure of the probability matrix leads to a variety of breakthroughs.

First, it enables us  to carry out estimation and clustering for the  PABM, without imposing any 
identifiability conditions, similarly to SBM and unlike the DCBM and mixed membership models. 
Second, our understanding of the probability matrix structure leads to  the Frobenious  norm minimization
as the basis of optimization procedure and to estimation of probability matrices by rank one approximations of the community matrices.
The latter allows us  to derive  non-asymptotic upper bounds for the estimation error, even in the case when 
the number of communities is unknown and is possibly growing with $n$. 
In addition, we  use the accuracy of approximation of the adjacency matrix 
for various number of communities, to identify the number of communities in the network.
Moreover, we   formulate  detectability conditions that guarantee that communities are identifiable, 
i.e.,  for the true probability matrix, the solution of the optimization problem is given by the true community assignment. 
Under those conditions, we provide a non-asymptotic upper bound on the proportion of the misclassified nodes
when the clustering is based on the solution of the optimization problem above.

Furthermore, we note that, under the detectability condition, the columns of the probability matrix that correspond to any of the  communities
lie in a  $K$-dimensional subspace, which is different from subspaces corresponding to all other communities. 
The   latter conclusion results in the  introduction of the Sparse Subspace Clustering (SSC) approach for partitioning the network into communities.
While the SSC is widely used in computer vision,  to the best of our knowledge, it has never been used for clustering network data.
The advantage of the SSC procedure  is that it is known to work very well in practice and  has several well studied versions 
(see Section~\ref{sec:SSC_review} for the discussion of the SSC algorithms). Moreover, unlike the Extreme Point algorithm 
which  \cite{RePEc:bla:jorssb:v:80:y:2018:i:2:p:365-386} managed to implement only in the case of $K=2$,
the SSC works well for an arbitrary number of communities. Our simulation study, as well as the real data examples,  
handle various number of communities  between 2 and 6. In particular, we demonstrate the advantages of the PABM for 
modeling  networks that appear in biological sciences.

We show that, under  the detectability condition, the SSC delivers the correct community assignment at population level.
We discuss the state of the art results for  the accuracy of the SSC approach and point out why 
they cannot be applied directly in  the case of the independent Bernoulli errors.
Investigation of the precision of the SSC for such errors is the matter of future work.

The rest of the paper is organized as follows. Section~\ref{sec:est_clust}
considers estimation and clustering in PABM as a solution of a penalized optimization procedure,
and investigates its accuracy. Specifically, 
Section~\ref{sec:notation} introduces notations  used throughout the paper.
Section~\ref{sec:opt_proc} formulates estimation and clustering as solutions of an optimization procedure.
Section~\ref{sec:est_errors} derives an upper bound for the estimation error in the case when the number of communities in the PABM is unknown.
Section~\ref{sec:sparse_PABM} delivers an upper bound for estimation errors in the case when all probabilities of connections are uniformly small, and 
 also  discusses advantages of the PABM for  modeling sparsity when this assumption is not true.
Section~\ref{sec:cl_detectability} provides detectability conditions at the population level.
Section~~\ref{sec:cl_err_ideal}  offers sufficient conditions for the proportion of misclassified 
nodes to be bounded above by a pre-specified quantity $\rho_n$ with a high probability. 
Since the optimization problem in Section~\ref{sec:opt_proc}  is NP-hard, Section~\ref{sec:SSC}
presents a computationally tractable way of finding communities by the Subspace Clustering.
In particular, Section~\ref{sec:SSC_review} reviews the Sparse Subspace Clustering (SSC) methodologies  and 
elaborates on what kind of SSC procedure we employ in this paper.  Section~\ref{sec:SSC_correctness}
shows that the SSC delivers correct community assignment  at the population level while  Section~\ref{sec:SSC-errors}
investigates this question in the case when the SSC is applied to the adjacency matrix.
Section~\ref{sec:simul_real} deliberates about 
practical implementation of clustering and provides a simulation study and real data examples.  
Finally, Section~\ref{sec:AppendA} presents the proofs of all statements in the paper.

%%%%%%%%%%%%%%%%%%%%%%%%%%%%%%%%%%%%%%%%%%%%%%%%%%%%%%%%%%%%%%%%%%%%%%%%%%%%%%%%%%%%%%%%%%%%%%%%%%%%%%%%%%%%%%%%%%%%%%%%%%%%%%
%%%%%%%%%%%%%%%%%%%%%%%%%%%%%%%%%%%%%%%%%%%%%%%%%%%%%%%%%%%%%%%%%%%%%%%%%%%%%%%%%%%%%%%%%%%%%%%%%%%%%%%%%%%%%%%%%%%%%%%%%%%%%%

\section{Estimation and  clustering}
\label{sec:est_clust}
\setcounter{equation}{0}

%%%%%%%%%%%%%%%%%%%%%%%%%%%%%%%%%%%%%%%%%%%%%%%%%%%%%%%%%%%%%%%%%%%%%%%%%%%%%%%%%%%%%%%%%%%%%%%%%%%%%%%%%%%%%%%%%%%%%%%%%%

\subsection{Notation}
\label{sec:notation}

For any two positive sequences $\{ a_n\}$ and $\{ b_n\}$, $a_n \lesssim b_n$ and $a_n \asymp b_n$ mean  that 
there exists a constant $C>0$ independent of $n$ such that, respectively,  $a_n \leq C b_n$  and 
$C^{-1} a_n \leq b_n \leq C a_n$ for any $n$. For any set $\Om$, denote cardinality of $\Om$ by $|\Om|$.
For any numbers $a$ and $b$, $a \wedge b = \min (a,b)$.
% For any $x$, $[x]$ is the largest integer no larger than $x$. 
% 
For any vector $ t \in \RR^p$, denote  its $\ell_2$, $\ell_1$, $\ell_0$ and $\ell_\infty$ norms by, 
respectively,  $\|  t\|$, $\|  t\|_1$,  $\|  t\|_0$ and $\|  t\|_\infty$.
% Denote by $\| t_1 -  t_2\|_H$ the Hamming distance between vectors $t_1$ and $t_2$.
Denote by $1_m$  the $m$-dimensional column vector with all components equal to one.  
For any matrix $A$,  denote its spectral and Frobenius norms by, respectively,  $\|  A \|_{op}$ and $\|  A \|_F$.
Let $\vect(A)$ be the vector obtained from matrix $A$ by sequentially stacking its columns.

Denote by $\calM_{n,K}$ a collection of   clustering  matrices $Z \in \{0,1\}^{n\times K} $  
 such that $Z_{i,k}=1$ iff $i \in \calN_k$, $i=1, \ldots, n$, and 
$Z^T Z = \diag (n_1, \ldots, n_K)$ where  $n_k = |\calN_k|$ is the size of community $k$, where $k=1, \ldots, K$.
Denote $n_{\min} = \di \min_k n_k$. Denote by  $\scrPZK \in \{0,1\}^{n \times n}$  
the   permutation matrix corresponding to $Z \in \calM_{n,K}$ that rearranges  
any matrix   $B \in R^{n,n}$, so that its  first $n_{1}$ rows correspond to  nodes from class 1, 
the next $n_{2}$ rows correspond to  nodes from class 2 and the last $n_K$ rows correspond to  nodes from class $K$. 
Recall that $\scrPZK$ is an orthogonal matrix with $\scrPZK^{-1} = \scrPZK^T$. 
For any $\scrPZK$ and any matrix $B \in \RR^{n \times n}$ denote the permuted matrix and its blocks by, respectively,
$B(Z,K)$ and $B^{(k,l)} (Z,K)$, where  $B^{(k,l)}(Z,K) \in \RR^{n_k \times n_l}$,  $k,l=1, \ldots, K$,  and
\be  \label{eq:permute}
B(Z,K) = \scrPZK^T B \scrPZK,  \quad   \quad B = \scrPZK B(Z,K) \scrPZK^T.
\ee  
Also, throughout the paper, we use the star symbol to identify the true quantities. In particular, we 
denote the true matrix of connection probabilities by $P_*$, the true number of classes by $K_*$
and the true clustering matrix that partitions $n$ nodes into $K_*$ communities by $Z_*$.

%%%%%%%%%%%%%%%%%%%%%%%%%%%%%%%%%%%%%%%%%%%%%%%%%%%%%%%%%%%%%%%%%%%%%%%%%%%%%%%%%%%%%%%%%%%%%%%%%%%%%%%%%%%%%%%%%%%%%%%%%%

\subsection{Optimization procedure for estimation and clustering}
\label{sec:opt_proc}

\begin{figure}[t]
\[\hspace{8mm}\includegraphics[height = 4.2cm]{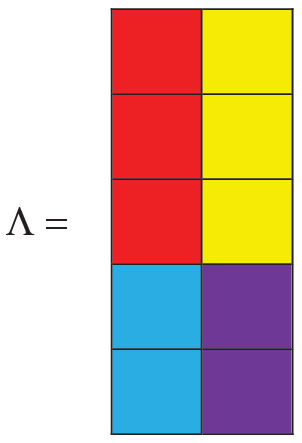} \hspace{15mm}  \includegraphics[height = 4.2cm]{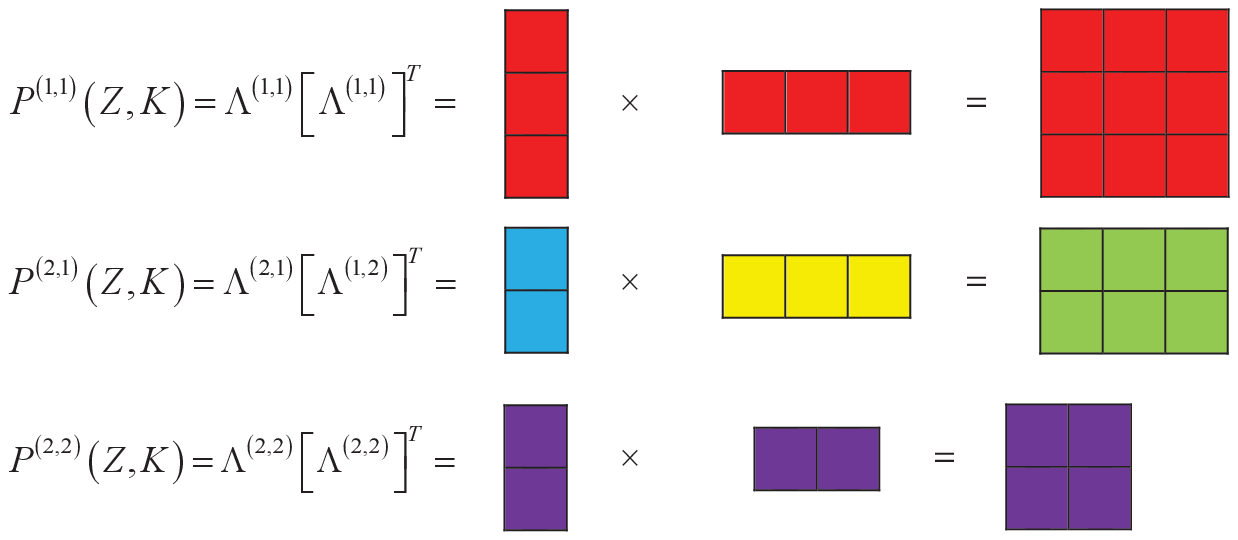} \]
\[\hspace{-1mm}\includegraphics[height = 4.2cm]{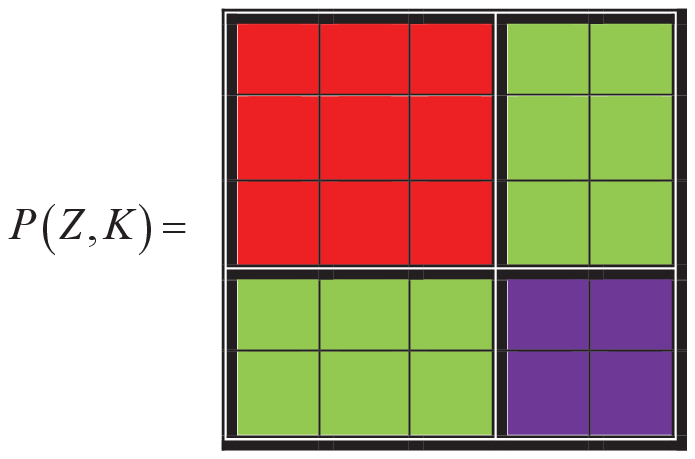}   \hspace{33mm} \includegraphics[height = 4.2cm]{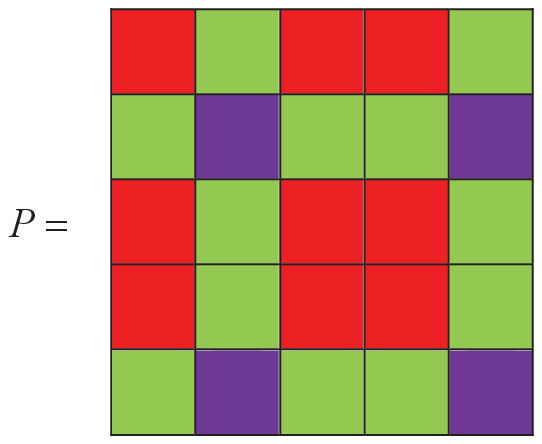}  \hspace{21mm} \]
\caption{Matrices $\Lam$, $P (Z,K)$ and $P$ in the case of $n=5$ and $K=2$.     Matrix  $\Lambda$  (top left):  $\Lambda^{(1,1)}$  (red),  
$\Lambda^{(2,1)}$  (blue), $\Lambda^{(1,2)}$  (yellow), $ \Lambda^{(2,2)}$ (violet). 
Assembling re-organized probability matrix  $P(Z,K)$ (top right):  $P^{(1,1)}(Z,K)$  (red),
 $P^{(2,1)}(Z,K)$ (green), $P^{(2,2)}(Z,K)$ (violet). 
Re-organized probability matrix  $P(Z,K)$  (bottom left): $P^{(1,1)}(Z,K)$  (red),
 $P^{(2,1)}(Z,K)$ and $P^{(1,2)}(Z,K)$  (green), $P^{(2,2)}(Z,K)$ (violet). 
Probability matrix  $P$ (bottom right):    nodes  1,3,4 are in community 1;  
nodes 2 and 5 are in community 2.
}
\label{fig:PZK}
\end{figure}

In   this section we consider estimation of the true probability matrix $P_*$.
Consider block $P_*^{(k,l)}(Z_*,K_*)$ of the rearranged version $P_*(Z_*,K_*)$ of $P_*$. 
Let $\Lam \equiv \Lam (Z_*,K_*) \in [0,1]^{n \times K_*}$ be a  block matrix  with each column $l$ partitioned into $K_*$
blocks $\Lam^{(k,l)} \equiv  \Lam^{(k,l)}(Z_*,K_*)  \in [0,1]^{n_k}$. Then, due to \eqref{eq:block_structure}, 
$P_*^{(k,l)} (Z_*,K_*)$ are rank-one matrices
such that $P_*^{(k,l)} (Z_*,K_*) = [P_*^{(l,k)} (Z_*,K_*)]^T$ and that each pair of blocks $(k,l)$ involves a  unique combination 
of vectors $\Lambda^{(k,l)}$. The structures of matrices  $P_*(Z_*,K_*)$, $\Lam$ and $P_*$ are illustrated in Figure \ref{fig:PZK}.

Observe that although matrices $P_*^{(k,l)} (Z_*,K_*)$ in \eqref{eq:block_structure} are well defined, 
vectors $\Lambda^{(k,l)}$ and $\Lambda^{(l,k)}$ can be determined only up to a 
multiplicative constant. In particular, under  the constraint 
\be  \label{eq:SC_constraint}
1_{n_k}^T \Lambda^{(k,l)}  = 1_{n_l}^T \Lambda^{(l,k)},
\ee  
\cite{RePEc:bla:jorssb:v:80:y:2018:i:2:p:365-386}  obtained explicit expressions for vectors $\Lambda^{(k,l)}$
and $\Lambda^{(l,k)}$ in   \eqref{eq:block_structure}.
\ignore{ 
\be \label{eq:Lambda_recov_SC}
\Lambda^{(k,l)}(Z,K) = \frac{P_*^{(k,l)} (Z_*,K_*) 1_{n_l}}{\sqrt{1_{n_k}^T P_*^{(k,l)} (Z_*,K_*) 1_{n_l}}};\quad 
\Lambda^{(l,k)}(Z,K) = \frac{(P_*^{(k,l)} (Z_*,K_*))^T 1_{n_k}}{\sqrt{1_{n_k}^T P_*^{(k,l)} (Z_*,K_*) 1_{n_l}}}
\ee
}
In reality,   $K_*$ and matrices  $Z_*$ and $P_*$ are unknown and need to be recovered.
If $K_*$ were known, in order to estimate $Z_*$ and $P_*$, one could permute the rows and the columns of the adjacency 
matrix $A$ using permutation matrix $\mathscr{P}_{Z,K_*}$ obtaining matrix 
$A(Z,K_*) = \mathscr{P}_{Z,K_*}^T A \mathscr{P}_{Z,K_*}$  
and then, following assumption  \eqref{eq:block_structure},  
minimize some divergence measure between blocks of $A(Z,K_*)$ and 
the products $\Lambda^{(k,l)} \, [\Lambda^{(l,k)}]^{T}$. 
One of such measures is the Bregman divergence between $A(Z,K_*)$
and $\Lambda^{(k,l)} \,  [\Lambda^{(l,k)}]^{T}$.

The Bregman divergence between vectors $x$ and $y$  associated with a continuously-differentiable, 
strictly convex function  $F$ is defined as 
\bes
D_{F}(x,y)=F(x)-F(y)-\langle \nabla F(y),x-y\rangle  
\ees
where $\nabla F(y)$ is the gradient of $F$ with respect to $y$. 
The Bregman divergence between any matrices $X$ and $Y$ of the same dimension can be defined as 
the Bregman divergence between their vectorized versions:  $D_{F}(X,Y) = D_{F}(\vect(X),\vect(Y))$.
It is well known that  $D_{F}(X,Y) \geq 0$ for any $X$ and $Y$ and  $D_{F}(X,Y) = 0$ iff $X=Y$.
In particular, the Poisson log-likelihood maximization used in \cite{RePEc:bla:jorssb:v:80:y:2018:i:2:p:365-386}  
corresponds to minimizing the Bregman divergence with 
% $F(x) = {\displaystyle  \sum_{i} (x_i \ln x_i - x_i)}$. 
\bes
F(x) =   \sum_{i} (x_i \ln x_i - x_i). 
\ees
Under the assumption \eqref{eq:block_structure} and  the constraint \eqref{eq:SC_constraint} 
of \cite{RePEc:bla:jorssb:v:80:y:2018:i:2:p:365-386},
the latter leads to maximization over  $\Lambda^{(k,l)}$ and $Z \in \calM_{n,K_*}$   of the following quantity
\be \label{eq:poiss}
l(\Lambda| A) = - D_{F}(A,\Lambda) = \sum_{k,l=1}^{K_*} \sum_{i=1}^{n_k} \sum_{j=1}^{n_l} 
\lkv  A^{(k,l)}_{i,j}  \ln \lkr \Lambda^{(k,l)}_{i}   \,  \Lambda^{(l,k)}_j \rkr  -  
\lkr \Lambda^{(k,l)}_{i}   \,  \Lambda^{(l,k)}_j \rkr \rkv.
\ee 
where  $A^{(k,l)}$ stands for $A^{(k,l)}  (Z,K_*)$, the $(k,l)$-th block of matrix $A (Z,K_*)$.  
It is easy to see that the expression \eqref{eq:poiss} 
coincides with the Poisson log-likelihood up to a term which depends on matrix $A$ only, and  
is independent of $P,Z$ and $K_*$.  Maximization of  \eqref{eq:poiss}  over $\Lambda$,
under condition \eqref{eq:SC_constraint}, for given $Z$ and $K_*$,
leads to the estimators of $\Lambda$ obtained in \cite{RePEc:bla:jorssb:v:80:y:2018:i:2:p:365-386}
\be \label{eq:Lambda_est}
\widehat{\Lambda}^{(k,l)}  = \frac{A^{(k,l)} (Z,K_*) 1_{n_l}}{\sqrt{1_{n_k}^T A^{(k,l)} (Z,K_*) 1_{n_l}}};\quad 
\widehat{\Lambda}^{(l,k)}  = \frac{(A^{(k,l)} (Z,K_*))^T 1_{n_k}}{\sqrt{1_{n_k}^T A^{(k,l)} (Z,K_*) 1_{n_l}}}.
\ee
% which correspond  to replacing matrix $P_*$  by $A (Z,K_*)$ in the expressions \eqref{eq:Lambda_recov_SC}.
Afterwards, \ignore{Sengupta and Chen} \cite{RePEc:bla:jorssb:v:80:y:2018:i:2:p:365-386} plug   the estimators   \eqref{eq:Lambda_est} 
into \eqref{eq:poiss}, thus,  obtaining the likelihood modularity function which they 
further maximize in order to obtain community assignments.

In the present paper, we use the Bregman divergence associated with the Euclidean distance $(F(x) = \|x\|^2)$ which, for a given $K$, leads to 
the following optimization problem
\bes % \label{eq:naive_opt}
% \begin{array}{ll}
  (\widehat{\Lambda}, \widehat{Z}) \in   \underset{\Lambda, Z}{\text{argmin}}  
\left\{\displaystyle   \sum_{k,l = 1}^K \norm {A^{(k,l)}(Z,K)  - \Lambda^{(k,l)} [\Lambda^{(l,k)}]^T }_{F}^2  \right\}   
\quad  \text{s.t.}\quad A(Z,K) = \scrPZK^T A \scrPZK 
% \end{array}
 \ees 
Note that recovery of the components  $\Lambda^{(k,l)}$ and $\Lambda^{(l,k)}$ of the products above relies on 
an identifiability condition of the type \eqref{eq:SC_constraint}.
Since these conditions can be imposed in a variety of ways, we denote 
$\Theta^{(k,l)} =   \Lambda^{(k,l)} [\Lambda^{(l,k)}]^T$ and recover the uniquely defined rank one  matrix $\Theta^{(k,l)}$.
In addition, since   the number of clusters  $K$ is unknown, we  impose a penalty   on   $K$
in order to safeguard against choosing too many clusters. 
Hence, we need to solve the following  optimization problem 
\be  \label{eq:opt_main}  
\begin{array}{ll} 
 (\hat{\Theta}, \hat{Z},\hat{K}) \in & \underset{\Theta, Z,K}{\text{argmin}} 
\left\{\displaystyle \sum_{k,l = 1}^K \norm {A^{(k,l)}(Z,K)  - \Theta^{(k,l)} }_{F}^2 + \Pen (n,K) \right\} \\
 & \text{s.t.}\quad A(Z,K) = \scrPZK^T A \scrPZK,\quad \text{rank} (\Theta^{(k,l)}) =1 ;\quad k,l = 1,2, \cdots , K.
 \end{array}
 \ee
Here, $\hat{\Theta}$ is the block matrix with blocks $ \hat\Theta^{(k,l)}$, $k,l=1, \ldots, \hat{K}$ and $\Pen (n,K)$
will be defined later.

Observe that, if $\hat{Z}$ and $\hat{K}$   were known, the best solution of problem \eqref{eq:opt_main} would be given by 
the rank one approximations   $\hat{\Theta}^{(k,l)}$ of matrices $A^{(k,l)}(\hat{Z},\hat{K})$ 
\be \label{eq:Theta_est}
\hat{\Theta}^{(k,l)} (\hat{Z},\hat{K}) =  \Pi_{\hat{u}, \hat{v}} \left(A^{(k,l)}(\hat{Z},\hat{K})\right) 
= \hat{\sigma}_{1}^{(k,l)}\hat {u}^{(k,l)}(\hat{Z},\hat{K})(\hat{v}^{(k,l)} (\hat{Z},\hat{K}))^T,
\ee
where $\hat{\sigma}_{1}^{(k,l)}$ are the largest singular values of matrices $A^{(k,l)}(\hat{Z},\hat{K}))$;  
$\hat {u}^{(k,l)}(\hat{Z},\hat{K})$, $\hat{v}^{(k,l)} (\hat{Z},\hat{K})$ are the corresponding 
singular vectors, and $\Pi_{\hat{u}, \hat{v}} \left(A^{(k,l)}(\hat{Z},\hat{K})\right)$ is the rank one projection of matrix 
$A^{(k,l)}(\hat{Z},\hat{K})$ (see Lemma~\ref{lem:lowrank_approx} in   Section~\ref{sec:AppendA}  for the exact expression).
%Rao and Rao
Due to the Perron-Frobenius theorem ( \cite{rao_rao_1998}, {\bf P.15.1.14}), $\hat{\sigma}_{1}^{(k,l)}>0$
and elements of vectors $\hat {u}^{(k,l)}(\hat{Z},\hat{K})$  and  $\hat{v}^{(k,l)} (\hat{Z},\hat{K})$ are
non-negative. 
Plugging \eqref{eq:Theta_est} into \eqref{eq:opt_main}, we rewrite optimization problem 
\eqref{eq:opt_main} as 
\be  \label{eq:opt_ZK}
\begin{array}{ll}
(\hat{Z},\hat{K}) \in & \underset{ Z,K}{\text{argmin}}  \left\{ \displaystyle  \sum_{k,l = 1}^K \norm {A^{(k,l)}(Z,K) -  
\Pi_{\hat{u}, \hat{v}} \left(A^{(k,l)}(Z,K)\right)}_{F}^2  + \Pen (n,K)\right\}\\
 & \text{s.t.}\quad A(Z,K) = \scrPZK^T A \scrPZK 
\end{array}
 \ee
In order to obtain $(\hat{Z},\hat{K})$, one needs to solve optimization problem \eqref{eq:opt_ZK}   for every $K$, 
obtaining
\be \label{eq:opt_ZK3}
\begin{array}{ll}
\hat{Z}_K \in & \underset{Z \in \calM_{n,K}}{\text{argmin}}  \left\{ \displaystyle \sum_{k,l = 1}^K \norm {A^{(k,l)}(Z,K) -  
\Pi_{\hat{u}, \hat{v}} \left(A^{(k,l)}(Z,K)\right)}_{F}^2 \right\} 
\end{array}
\ee 
and then find $\hat{K}$ as 
\be \label{eq:opt_for_K}
\hat{K}  \in   \underset{K}{\text{argmin}}  \left\{ \sum_{k,l = 1}^K \norm {A^{(k,l)}(\hat{Z}_K,K) -  
\Pi_{\hat{u}, \hat{v}} \left(A^{(k,l)}(\hat{Z}_K,K)\right)}_{F}^2  + \Pen (n,K)\right\}.
\ee 
Note that if the true number of clusters $K_*$ were known, the penalty in \eqref{eq:opt_main}  and \eqref{eq:opt_ZK}
would be unnecessary.

%%%%%%%%%%%%%%%%%%%%%%%%%%%%%%%%%%%%%%%%%%%%%%%%%%%%%%%%%%%%%%%%%%%%%%%%%%%%%%%%%%%%%%%%%%%%%%%%%%%%%%%%%%%%%%%55
%%%%%%%%%%%%%%%%%%%%%%%%%%%%%%%%%%%%%%%%%%%%%%%%%%%%%%%%%%%%%%%%%%%%%%%%%%%%%%%%%%%%%%%%%%%%%%%%%%%%%%%%%%%%%%%55

\subsection{The  penalty and the   estimation errors }
\label{sec:est_errors}

In this section we evaluate  the estimation and the clustering errors. We choose the penalty which, with high probability, 
exceeds the random errors. In particular, we denote
\be \label{eq:pen1}
\Pen (n,K) = H_1 n K + H_2 K^2 \ln n + H_3  n \ln K,
\ee 
where $H_1, H_2$ and $H_3$ are positive absolute constants that can be evaluated. 
Then, the following statement holds.

\begin{thm}  \label{th:oracle}
Let $(\hat{\Theta}, \hat{Z},\hat{K})$ be a solution of optimization problem  \eqref{eq:opt_main}.
Construct the estimator $\hat{P}$  of $P_*$ of the form 
 \be \label{eq:P_total_est}      
\hat{P}  =\mathscr{P}_{\hat{Z},\hat{K}} \hat\Theta(\hat{Z},\hat{K})\mathscr{P}_{\hat{Z},\hat{K}}^T 
\ee
where $\mathscr{P}_{\hat{Z},\hat{K}}$ is the permutation matrix corresponding to  $(\hat{Z},\hat{K})$. 
Then, for any $t >0$ and some absolute positive constants $H$  and  $\tilde{C}$, one has 
\be \label{eq:oracle}
\PP \lfi  n^{-2}\, \norm{\hat{P}  -P_{*}}_F^2  \leq   n^{-2}\,  H\, \Pen (n,K_{*})    + n^{-2}\,    \tilde{C}\,t  \rfi
\geq  1 - 3 e^{-t},
\ee 
\be \label{eq:oracle Expectation}
n^{-2}\,  \EE\norm{\hat{P}  - P_{*}}_F^2  \leq   
n^{-2}\, H\,  \Pen (n,K_{*})  +  n^{-2}\, \tilde{C}.  
\ee 
\end{thm}

%%%%%%%%%%%%%%%%%%%%%%%%%%%%%%%%%%%%%%%%%%%%%%%%%%%%%%%%%%%%%%%%%%%%%%%%%%%%%%%%%%%%%%%%%%%%%%%%%%%%%%%%%%%%%%%55

The exact values of $H$ and $\tilde{C}$ can be found in the proof of Theorem~\ref{th:oracle}.
Observe that estimation is  always consistent as long as $K/n \to 0$.  
Note   also that  the estimation errors in 
\eqref{eq:oracle} and \eqref{eq:oracle Expectation} are proportional to the right hand side of \eqref{eq:pen1}.
The first term  in \eqref{eq:pen1} corresponds to the error of estimating $n K$ unknown entries 
of matrix $\Lambda$, the second term is  associated with estimation of rank $K^2$ matrix while the last term is due
to the clustering of $n$ nodes into $K$ communities. If $K$ grows with $n$, i.e., $K = K(n) \to \infty$ as $n \to \infty$, then 
the first term in \eqref{eq:pen1} dominates the other two terms. However, in the case of a fixed $K$, 
the first and the third terms grow at the same rate as $n \to \infty$. The second term is always of a smaller order provided 
$K(n)/n \to 0$.

%%%%%%%%%%%%%%%%%%%%%%%%%%%%%%%%%%%%%%%%%%%%%%%%%%%%%%%%%%%%%%%%%%%%%%%%%%%%%%%%%%%%%%%%%%%%%%%%%%%%%%%%%%%%%%%55

\subsection{The  sparse PABM }
\label{sec:sparse_PABM}

The real life networks are usually sparse in a sense that a large number of nodes have small degrees.
One of the advantages of the PABM is that it allows flexible modeling of sparsity.
Traditionally, in most statistical models, sparsity of a vector means that a large proportion of its components
is  equal to zero. One of the shortcomings of both the SBM and the DCBM is that they do not allow to 
impose the condition that some of the connection probabilities are equal to zero. 
Naturally, for the SBM, it is not realistic to assume that all nodes in a pair of communities have no connections.
Neither can one set any of the node-specific weight to zero, since this   will 
force the respective node to be totally disconnected from the network.
For this reason, unlike in other numerous statistical settings,   sparsity in block models is defined as a 
low maximum probability of connections between the nodes:   
\be \label{eq:uniform_sparse}
\max_{i,j} P_{i,j} \leq \tau_n,
\ee
where  $\tau_n$ is small when $n$ is large   (see, e.g., \cite{KlTsVe2017} and \cite{lei2015}).

There are several shortcomings of  this definition of sparsity. 
First, even in the  context of the simplest model, the SBM, in order to take a full advantage of assumption \eqref{eq:uniform_sparse},
one needs to carry out the estimation under the restriction that all entries of the matrix $\hat{P}$ are bounded above 
by $\tau_n$  (see  \cite{KlTsVe2017}), which is an unknown quantity.

On the other hand, in the context of the PABM, one can take advantage of sparsity in a much more natural way.  
Indeed, unlike the SBM and the DCBM,  the PABM setting allows some connection probabilities to be zero while keeping 
average connection probabilities between classes above certain level and the network connected. 
This is certainly true since  setting $\Lambda^{(k,l)}_i = 0$   in the   PABM
simply means that that node $i$ in class $k$ is not active  (``popular") in class $l$. 
The latter  does not prevent   node $i$  from having high probability of connection with nodes in another class.

Therefore, the PABM, similarly  to other sparse statistical  settings, allows structural sparsity
where small parameters are set to zero rather than considered to be infinitesimally small.
Setting some of the connection probabilities to zero, rather than bounding all 
of the connection probabilities  by  a very small  number, as in \eqref{eq:uniform_sparse}, 
not only leads to   better understanding of network topology but also allows  more precise 
estimation of the probability matrix $P_*$. Furthermore, this approach enables one to 
handle the unknown number of communities that is possibly growing with $n$.
While we do not consider the structurally sparse PABM in this paper,
we investigate the structurally sparse PABM  in depth in our subsequent publication \cite{noroozi2019sparse}.

%%%%%%%%%%%%%%%%%%%%%%%%%%%%%%%%%%%%%%%%%%%%%%%%%%%%%%%%%%%%%%%%%%%%%%%%%%%%%%%%%%%%%%%

Below, we briefly consider the case of the uniformly sparse PABM satisfying condition  \eqref{eq:uniform_sparse}.
In this case, the main error term $nK$ in \eqref{eq:oracle} and \eqref{eq:oracle Expectation} is replaced  by $\tau_n n K$, 
which can significantly reduce the error if $K = K(n) \to \infty$. However, the drawback of this approach is that one needs 
to know either the sparsity level $\tau_n$ or the number of communities $K$.
The reason for this is that the penalty term, which offsets the random error, should contain a component  
$C \tau_n n K$  where $C$ is an absolute constant. The latter quantity 
may not be monotone  since  $K = K(n)$  is growing with $n$, while $\tau_n$ is decreasing with $n$.

For this reason, we derive the estimation error under a more common scenario that the 
number of communities is known: $K=K_*$. 
In this case, penalty is unnecessary and one can just solve optimization problem  \eqref{eq:opt_ZK3}
for the known number of communities.

 \begin{thm}  \label{th:oracle_sparse}
Let $K = K_*$ be known,   $\hat{Z}$ be a solution of optimization problem \eqref{eq:opt_ZK3}
with $K = K_*$. Let $\hat\Theta= \hat\Theta(\hat{Z})$ be the matrix with blocks $\hat{\Theta}^{(k,l)}$ given by
\eqref{eq:Theta_est}. 
Construct the estimator $\hat{P}$  of $P_*$ of the form \eqref{eq:P_total_est} 
where $\mathscr{P}_{\hat{Z}}$ is the permutation matrix corresponding to  $\hat{Z}$. 
Assume that $n_{\min} = \di \min_k (n_k)$ is large enough, so that  
\be \label{eq:nmin_cond}
\log(2 n_{\min}) \leq (2\, n_{\min})^{2/13}, \quad 
\tau_n \geq C_\tau  \log(2\, n_{\min})/ n_{\min} 
 \ee
 for some absolute constant $C_{\tau}  >0$ and $\tau_n$ in \eqref{eq:uniform_sparse}.
Then, for any $t >0$ and some  absolute positive constants  $\tilde{H}$, $H_1, H_2$ and $H_3$, one has 
\be \label{eq:oracle_sparse}
\PP \lfi  n^{-2}\, \|\hat{P}  -P_{*}\|_F^2    \leq    H_1 \tau_n n^{-1}\, K + H_2 n^{-2}\, K^2   + H_3  n^{-1}\, \ln K   
+ n^{-2}\,    \tilde{H}\,t  \rfi \geq  1 - 3 e^{-t}, 
\ee 
\be  \label{eq:oracle Expectation_sparse}
n^{-2}\,  \EE \|\hat{P}  - P_{*}\|_F^2   \leq   
 H_1 \tau_n n^{-1}\, K + H_2 n^{-2}\, K^2   + H_3  n^{-1}\, \ln K  +  n^{-2}\, \tilde{H}.   
\ee 
\end{thm}

\noindent
The advantage of Theorem~\ref{th:oracle_sparse} is that it replaces the main error term $O(n K)$
in Theorem~\ref{th:oracle} by the smaller quantity $O(\tau_n n K)$ and this is done without any knowledge 
of $\tau_n$. If $\tau_n \to 0$ as $n \to \infty$, the latter may be significantly smaller than the former. 
This reduction, however, comes at a price. First, application of Theorem~\ref{th:oracle_sparse} requires 
the knowledge of the number of communities $K$. Second, while results in Theorem~\ref{th:oracle} are non-asymptotic 
and are valid for any combination of $n$ and $K$, Theorem~\ref{th:oracle_sparse} requires not only $n$ but also $n_{\min}$ 
to be large via conditions \eqref{eq:nmin_cond}.

%%%%%%%%%%%%%%%%%%%%%%%%%%%%%%%%%%%%%%%%%%%%%%%%%%%%%%%%%%%%%%%%%%%%%%%%%%%%%%%%%%%%%%%%%%%%%%%%%%%%%%%%%%%%%%%%%%%%%%%%%%
%%%%%%%%%%%%%%%%%%%%%%%%%%%%%%%%%%%%%%%%%%%%%%%%%%%%%%%%%%%%%%%%%%%%%%%%%%%%%%%%%%%%%%%%%%%%%%%%%%%%%%%%%%%%%%%%%%%%%%%%%%

\subsection{Detectability of clusters  }  
\label{sec:cl_detectability}

In order to evaluate the clustering error, we assume that the true number of communities $K = K_*$ is known.
Let $Z_* \in \calM_{n, K_*}$ be the true clustering matrix. 
% 
% The first question that needs to be solved  is under what conditions the communities  are identifiable. 
% Indeed, if the data follow the SBM, then  trying to find communities under the PABM assumption is guaranteed to fail since 
%  the probability matrix can be partitioned into rank one blocks that do not follow the true communities. 
In order for clustering to be successful, one needs a  detectability condition  that guarantees that communities are identifiable.
\\

{\bf   A1. \ }  There exists $k$, $1 \leq k \leq  K_*$, such that vectors ${\Lambda}^{(k,1)}, \ldots,  {\Lambda}^{(k,K_*)}$ 
are linearly independent and have all positive components.
\\

Assumption {\bf A1} is an alternative formulation of  the Detectability Assumption 4.4 of %Sengupta and  Chen 
\cite{RePEc:bla:jorssb:v:80:y:2018:i:2:p:365-386} which  states that, for any two nodes $j_1$ and $j_2$ 
that belong to different communities, the set $ \lfi  P_{i,j_1}/P_{i,j_2}  \rfi_{i=1}^n$ assumes at least 
$(K+1)$ distinct values.  Similarly, to \cite{RePEc:bla:jorssb:v:80:y:2018:i:2:p:365-386}, Assumption {\bf A1} 
guarantees that, for the true $K = K_*$, expression \eqref{eq:opt_ZK3} is minimized  at $Z=Z_*$.

Note that the assumption, that  all elements  of matrix $\Lambda$ are positive,   is necessary. 
Indeed, consider a PABM with $K=2$ and matrix $\Lambda$ such that ${\Lambda}^{(1,1)} =  {\Lambda}^{(1,2)}=u$  and 
${\Lambda}^{(2,1)}=v$ and ${\Lambda}^{(1,2)}=w$, where $v$ and $w$ are linearly independent.
If $u_i =0$ and $w_j=0$ for some $i$ and $j$, then matrix $P$ has two proportional  columns, $i$ and $j$, of the form $(cu, 0)^T$,
and   nodes $i$ and $j$ can be placed in any of the two communities.
In order to avoid the condition that all elements of matrix $\Lam$ are positive, one can use an alternative assumption.
\\

{\bf Assumption A1*. \ }  For any $k=1,\ldots, K_*$, vectors ${\Lambda}^{(k,1)}, \ldots,  {\Lambda}^{(k,K_*)}$ 
are linearly independent.
\\

% Then the following statement (similar  to Lemma 4.2 of \cite{RePEc:bla:jorssb:v:80:y:2018:i:2:p:365-386}) is true.

\begin{lem} \label{lem:detect}
Let Assumption {\bf A1}  or Assumption {\bf A1*}  holds.
Let $Z_* \in \calM_{n, K_*}$ be  the true clustering matrix and $Z \in \calM_{n, K_*}$
 be an arbitrary clustering matrix. Then, 
\be \label{eq:detect}
\sum_{k,l = 1}^K \|P_*^{(k,l)}(Z_*) -  \Pi_{(1)} (P_*^{(k,l)}(Z_*)\|_{F}^2 
\leq  \sum_{k,l = 1}^K \|P_*^{(k,l)}(Z) -  \Pi_{(1)} (P_*^{(k,l)}(Z)\|_{F}^2
\ee 
where, for any matrix $B$, $\Pi_{(1)} (B)$ is its rank one approximation. 
Moreover, equality in \eqref{eq:detect} occurs if and only if matrices  $Z$ and $Z_*$ coincide up to a permutation of columns.
\end{lem}

Lemma~\ref{lem:detect} implies that if $K = K_*$ is known, then optimization problem \eqref{eq:opt_main} leads to the 
true clustering assignment at a population level. The next section explores the clustering errors in the case 
when optimization procedure \eqref{eq:opt_main} is applied to the adjacency matrix.

%%%%%%%%%%%%%%%%%%%%%%%%%%%%%%%%%%%%%%%%%%%%%%%%%%%%%%%%%%%%%%%%%%%%%%%%%%%%%%%%%%%%%%%%%%%%%%%%%%%%%%%%%%%%%%%%%%%%%%%%%%
%%%%%%%%%%%%%%%%%%%%%%%%%%%%%%%%%%%%%%%%%%%%%%%%%%%%%%%%%%%%%%%%%%%%%%%%%%%%%%%%%%%%%%%%%%%%%%%%%%%%%%%%%%%%%%%%%%%%%%%%%%

\subsection{The clustering errors  }  
\label{sec:cl_err_ideal}

Note that if $Z_{*}$ is the true clustering matrix  and $Z$ is any other clustering matrix, then the proportion 
of misclassified nodes can be evaluated as 
 \be\label{eq:misclustered}
 \Err(Z, Z_{*}) =  (2n)^{-1}\, \underset{\mathscr{P}_K \in \mathcal{P}_K} {\min}  \norm{Z \mathscr{P}_K - Z_{*}} _1
= (2n)^{-1}\,  \underset{\mathscr{P}_K \in \mathcal{P}_K} {\min}  \norm{Z \mathscr{P}_K - Z_{*}} _F^2  
 \ee
where $\mathcal{P}_K$ is the set of permutation matrices  $\mathscr{P}_K: \{ 1,2,\cdots,K\} \longrightarrow \{ 1,2,\cdots,K\}$.
 Let 
\be\label{eq:def;Lambda}
 \Upsilon(Z_{*},\rho) = \lfi Z \in \mathcal{M}_{n,K } :  (2n)^{-1}\,  \underset{\mathscr{P}_K \in \mathcal{P}_K} 
{\min}  \norm{Z \mathscr{P}_K - Z_{*}} _1 \geq \rho \rfi 
 \ee
be the set of clustering matrices with the proportion of misclassified  nodes being at least  $\rho_n$, $0 < \rho_n < 1$.

The success of clustering in \eqref{eq:opt_ZK3} relies upon the fact that matrix $P_*$ is a collection of $K^2$ rank one blocks,
so that the operator and the Frobenius norms of each block are the same.  On the other hand, if clustering were incorrect,
the ranks of the blocks would increase which would lead to the discrepancy between their operator and Frobenius norms. 
In particular, the following statement is true.

\begin{thm}  \label{th:clust}
Let  $K = K_*$ be the true number of clusters and $Z_* \in \calM_{n, K_*}$ be the true clustering matrix.
Let Assumption {\bf A1}  or  {\bf A1*}  holds.
Let $\hat{Z}  \equiv \hat{Z}_K$ be a solution of the optimization problem \eqref{eq:opt_ZK3}.  
If for some  $\alpha_n  \in (0,1/2)$  and $\rho_n  \in (0,1)$, one has 
 \be \label{eq:cond_check}
 \norm{P_{*}}_{F}^2 -  (1 + \alpha_n)\,  \underset{Z \in \Upsilon(Z_{*},\rho_n)} {\max} 
\sum_{k,l = 1}^K \norm{  P_{*}^{(k,l)} (Z) }_{op}^2  
 \geq \frac{H}{\alpha_n}\,   (n K + K^2 \ln  n), 
\ee
where $H$ is an absolute positive constant independent of $K$,$n$, $\rho_n$ and $\alpha_n$,
then, with probability at least $1 - 2 e^{-n}$, the proportion of the nodes, misclassified  by $\hat{Z}$,
is at most $\rho_n$.
\end{thm}

In order to see what condition \eqref{eq:cond_check} means, we consider a simple example 
of the SBM   with $K=2$ and $P_{i,j} = b$ when nodes $i$ and $j$ belong to the same 
community and $P_{i,j} = r$  when they belong to different ones. 
Then,  condition \eqref{eq:cond_check} reduces to the following inequality.

\begin{lem} \label{lem:SBM_clust_er}
Consider the SBM   with $K=2$, $P = Z B Z^T$ where $B_{1,1} = B_{2,2} =b$, $B_{1,2} =r$
and $Z_* \in \calM_{n,2}$ with equal size communities.  Then,  for $0 < \rho \equiv  \rho_n  <\min(1,r^2/b^2)$
and $\alpha\equiv \alpha_n$, one has 
\be \label{eq:cond_SBM}
 \|P_{*}\|_{F}^2 -  (1 + \alpha)\,  \underset{Z \in \Upsilon(Z_{*},\rho)} {\max} 
\sum_{k,l = 1}^2 \|P_{*}^{(k,l)} (Z)\|_{op}^2  
 \geq 
 \tilC\, \lkv n^2 \,\rho^2\, (b^2 - r^2)^2\,  b^{-2} -  4\, \alpha \, n^2\, b^2 \rkv, 
\ee
where $\tilC$ is an absolute constant.
\end{lem}

%\noindent
% Setting $\alpha_n =   \rho_n^2 (1 - r^2/b^2)^2/8$, we obtain an upper bound on $\rho_n$.
%\\

\begin{cor} \label{cor:SBM_clust_er}
For the SBM in Lemma~\ref{lem:SBM_clust_er}, one has 
\be \label{eq:rho_bound}
\rho_n^2 \lesssim n^{-1/2}\, b^{-1} \, (b^2 - r^2)^{-2}. 
\ee
\end{cor}

The example above shows that condition \eqref{eq:cond_check} is less sensitive than  conditions
that are based on the difference between mean vectors of probabilities of connections between the communities in the case of the SBM (or scaled 
mean vectors of the communities in the case of the DCBM). Indeed,  it follows from \cite{Gao:2017:AOM:3122009.3153016}
that for  the SBM in Lemma~\ref{lem:SBM_clust_er} one can attain the misclassification rate   
\be \label{eq:Gao}
\rho_n \asymp \exp \lkr - \frac{n (b-r)^2}{2 K b} \rkr
\ee    
which is much smaller than $\rho_n$ in \eqref{eq:rho_bound}.
Nevertheless, achieving the misclassification rate   \eqref{eq:Gao}
depends upon not only the knowledge that the data is generated by the SBM, but also that this SBM is strongly assortative 
and balanced, and, in addition, requires  handling the Bernoulli likelihood. On the contrary, Theorem~\ref{th:clust}
is designed to work in the case where the communities are not characterized by their means and are not necessarily  assortative or balanced.
In addition, our procedure is based on minimizing the Frobenius norm which is much more computationally efficient 
but is less sensitive than the Bernoulli likelihood maximization.

On the other hand, we believe that the assessment of Theorem~\ref{th:clust} 
is valuable since it allows one  to upper bound   the misclassification rate  
rather than just stating that it tends to zero when the number of nodes in the network grows,
as it is routinely done in the papers that draw clustering assignments on the 
basis of modularity maximizations (see, e.g.,   Bickel and Chen \cite{Bickel21068}, 
Zhao {\it et al.} \cite{zhao2012consistency} and  \cite{RePEc:bla:jorssb:v:80:y:2018:i:2:p:365-386}).

%%%%%%%%%%%%%%%%%%%%%%%%%%%%%%%%%%%%%%%%%%%%%%%%%%%%%%%%%%%%%%%%%%%%%%%%%%%%%%%%%%%%%%%%%%%%%%%%%%%%%%%%%%%%%%%55
%%%%%%%%%%%%%%%%%%%%%%%%%%%%%%%%%%%%%%%%%%%%%%%%%%%%%%%%%%%%%%%%%%%%%%%%%%%%%%%%%%%%%%%%%%%%%%%%%%%%%%%%%%%%%%%55

\section{Sparse subspace clustering }
\label{sec:SSC}
\setcounter{equation}{0}

In Section~\ref{sec:est_clust}, we obtained an estimator $\hat{Z}$ of the true clustering matrix $Z_*$ 
as a solution of optimization problem \eqref{eq:opt_ZK}. 
Minimization in \eqref{eq:opt_ZK} is somewhat similar to modularity maximization
in \cite{Bickel21068},  \cite{RePEc:bla:jorssb:v:80:y:2018:i:2:p:365-386}  or \cite{zhao2012consistency}, in the sense that 
modularity maximization as well as minimization in \eqref{eq:opt_ZK} 
are NP-hard, and, hence,  require some relaxation in order to obtain  an implementable clustering solution.

In the case of the SBM and the DCBM,  possible relaxations  include semidefinite programming 
(see, e.g., \cite{DBLP:journals/corr/AminiL14} and references therein),
variational methods (\cite{celisse2012}) and  spectral clustering and its versions 
% such as  the approximate $K$-means, the spherical $K$-median or regularized spectral clustering 
(see, e.g.,    \cite{joseph2016}, \cite{lei2015}  and   \cite{rohe2011spectral} among others).
Since in the case of PABM, columns of matrix $P_*$ that correspond to nodes in the same class are neither identical, nor proportional,
application of spectral clustering (and its versions such as spherical spectral clustering) 
to matrix $P_*$ directly does not deliver the partition of the nodes.

However, it is easy to see that the columns of   matrix $P_*$ that correspond 
to nodes in the same class form a matrix with $K$ rank-one blocks,
hence, those columns lie in the subspace of the dimension  at most $K$. 
Therefore, matrix $P_*$ is constructed of $K$ clusters of 
columns (rows) that lie in the union of $K$   subspaces, each of the dimension $K$.
Under Assumption {\bf A1*}, those subspaces are independent in the sense that
the dimension of their union (the rank of $P_*$) is equal to the sum  $K^2$ of the dimensions $K$ of individual subspaces, 
and they can be recovered.
For this reason, the subspace clustering presents a technique for 
obtaining a fast and reliable solution of optimization problem~\eqref{eq:opt_ZK} (or \eqref{eq:opt_ZK3}.

%%%%%%%%%%%%%%%%%%%%%%%%%%%%%%%%%%%%%%%%%%%%%%%%%%%%%%%%%%%%%%%%%%%%%%%%%%%%%%%%%%%%%%%%%%%%%%%%%%%%%%%%%%%%%%%%%%%%%%%%%%
%%%%%%%%%%%%%%%%%%%%%%%%%%%%%%%%%%%%%%%%%%%%%%%%%%%%%%%%%%%%%%%%%%%%%%%%%%%%%%%%%%%%%%%%%%%%%%%%%%%%%%%%%%%%%%%%%%%%%%%%%%

\subsection{Review of the  subspace clustering}
\label{sec:SSC_review}

Subspace clustering has been widely used in computer vision and, for this reason, it is a very
well studied and developed technique in comparison with the    Extreme Points algorithm 
used in \cite{RePEc:bla:jorssb:v:80:y:2018:i:2:p:365-386}.   
Subspace clustering is designed for separation of points that lie in the union of subspaces. 
Let $\{ X_{j} \in \R^{D} \}_{j=1}^{n} $ be a given set of points drawn from an unknown union of 
$K \geqslant 1$ linear or affine subspaces $\{\calS_{i} \}_{i=1}^{K} $ of unknown dimensions 
$d_{i}= \text{dim}(\calS_{i})$, $0<d_{i} <D$, $i=1,...,K$. In the case of linear subspaces, the subspaces can be described as 
$$\calS_{i}=\{ \bm{x} \in \R^{D} : \bm{x}=  \bm{U}_{i}\bm{y} \} , \hspace{5mm} i=1,...,K $$
where $\bm{U}_{i} \in \R^{D \times d_{i}}$ is a basis for subspace $\calS_{i}$ 
and $\bm{y} \in \R^{d_{i}}$ is a low-dimensional representation for point $\bm{x}$. 
The goal of subspace clustering is to find the number of subspaces $K$, their dimensions 
$\{ d_{i} \}_{i=1}^{K}$, the subspace bases $\{ \bm{U}_{i} \}_{i=1}^{K}$, and the segmentation 
of the points according to the subspaces.

Several methods have been developed to implement subspace clustering such as  algebraic methods 
(\cite{inproceedings}, \cite{Ma:2008:ESA:1405158.1405160}, \cite{vidal2005generalized}), 
iterative methods (\cite{P.AgarwalandN.Mustafa2004}, \cite{Bradley:2000:KC:596077.596262},  \cite{tseng2000nearest}), 
and spectral clustering based methods (\cite{Vidal:2009aa},  \cite{Elhamifar:2013:SSC:2554063.2554078}, 
\cite{Favaro:2011:CFS:2191740.2191857},  \cite{Liu:2013:RRS:2412386.2412936}, 
\cite{Liu2010RobustSS},  \cite{soltanolkotabi2014},  
 \cite{vidal2011subspace}). 
In this paper, we shall use the latter group of techniques.

Spectral clustering algorithms rely on  construction of  an affinity matrix 
whose entries are based on some distance measures between the points. 
In particular, in the case of the SBM, adjacency matrix itself serves as the affinity matrix, while for the DCBM, 
the affinity matrix is obtained by normalizing rows/columns of $A$. 
In the case of the  subspace clustering problem, one cannot use the typical distance-based affinity
because two points could be very close to each other, but lie in different subspaces, while they 
could be far from each other, but lie in the same subspace. One of the solutions is to construct 
the affinity matrix using self-representation of the points with the expectation that a point is more likely to 
be presented as a linear combination of points in its own subspace rather than from a different one.  
A number of approaches such as Low Rank Representation  (see, e.g., \cite{Liu:2013:RRS:2412386.2412936} and 
\cite{Liu2010RobustSS}) and Sparse Subspace Clustering  (see, e.g., 
 \cite{Vidal:2009aa} and \cite{Elhamifar:2013:SSC:2554063.2554078})
have been proposed in the past decade   for the solution of this problem.

In this paper, we use Sparse Subspace Clustering (SSC) since it allows one to take advantage of the knowledge that, for a given $K$,
columns of matrix $P_*$  lie in the union of $K$ distinct subspaces, each of the dimension at most $K$. 
If matrix $P_*$ were known, the weight matrix $W$ would be based on writing every data point 
as a sparse linear combination of all other points by minimizing the number of nonzero coefficients 
\begin{equation}  \label{mn:opt_prob1}
\mathop\text{min}_{W_{j}}  \|W_{j}\|_{0}  \hspace{5mm}   \mbox{s.t}  \hspace{5mm}  (P_*)_{j}=\sum_{k \ne j} W_{kj} (P_*)_{k}
\end{equation}
where, for any matrix $B$,  $B_{j}$ is its   $j$-th column. The affinity matrix  
of the SSC is the symmetrised version of the weight matrix $W$. 
If the subspaces are linearly independent, then the solution to the optimization problem  \eqref{mn:opt_prob1}  
is such that $W_{k,j}\neq 0$ only if points $k$ and $j$ are in the same subspace. 
In the case of data contaminated by noise, the SSC algorithm does not attempt to write  
data as an exact linear combination of other points. 
Instead, SSC  is based on the solution of the following optimization problem
\begin{equation} \label{mn:NP_Hard_1}  
\widehat{W}_j \in \underset{W_{j}}{\text{argmin}}  
\lfi   \|W_{j}\|_{0} + \gamma \|{A_{j}-AW_{j}}\|_{2}^{2}  \quad \mbox{s.t.} \quad   W_{jj}=0 \rfi, \quad  j=1,...,n, 
\end{equation}
where $\gamma > 0$ is a tuning parameter. Problem \eqref{mn:NP_Hard_1}  can be rewritten in an equivalent form as 
\begin{equation} \label{mn:NP_Hard_2}
\widehat{W}_j \in \underset{W_{j}}{\text{argmin}}  
\lfi   \|{A_{j}-AW_{j}}\|_{2}^{2} \quad \mbox{s.t.} \quad \|W_{j}\|_{0} \leq {L}, \quad W_{jj}=0 \rfi, \quad j=1,...,n,
\end{equation}
where $L$ is  the maximum number of nonzero elements in each column of $W$; in our case $L = K$.
We solve  \eqref{mn:NP_Hard_2} using the Orthogonal Matching Pursuit (OMP) algorithm  (\cite{Mallat:1993:MPT:2198030.2203996}, 
\cite{weisberg2005applied}) implemented  in  SPAMS Matlab toolbox (see \cite{mairal2014spams}).  
Given $\widehat{W}$,  the affinity (similarity) matrix is defined as 
\be  \label{eq:similarity}
S  = |\widehat{W}| + |\widehat{W}^{T}|
\ee 
where, for any matrix $B$, matrix $|B|$ has absolute values of elements of $B$ as its entries.

The similarity matrix allows to construct the {\it similarity graph} 
$G = (V;E)$ where $(i,j) \in E$ if and only if $S_{i,j} >0$ (see, e.g., \cite{pmlr-v51-wang16b}).

%%%%%%%%%%%%%%%%%%%%%%%%%%%%%%%%%%%%%%%%%%%%%%%%%%%%%%%%%%%%%%%%%%%%%%%%%%%%%%%%%%%%%%%%%%%%%%%%%%%%%%%%%%%%%%%%%%%%%%%%%%%%%%
%%%%%%%%%%%%%%%%%%%%%%%%%%%%%%%%%%%%%%%%%%%%%%%%%%%%%%%%%%%%%%%%%%%%%%%%%%%%%%%%%%%%%%%%%%%%%%%%%%%%%%%%%%%%%%%%%%%%%%%%%%%%%%

\begin{rem}\label{rem:comput_complex}
{\bf (Computational complexity of the SSC.)\ \   }
{\rm 
Implementation of the SSC consists of two parts. The first part, evaluation of the  matrix $\widehat{W}$.
The second part is spectral clustering of the similarity matrix \eqref{eq:similarity}. %$|\widehat{W}| + |\widehat{W}^{T}|$. 
While the first component of the problem is more computationally expensive (since it requires solution of $n$ sparse regression problems),
it is also the portion that can be easily carried out via parallel computing. Indeed, evaluation of each of the vectors $W_j$, 
$j = 1,\ldots, n$, is completely independent from evaluation of all the others. Hence, if one has $m$ CPUs available, evaluation of matrix $W$ 
can be accomplished $m$ times faster. As a result, with the adequate facilities available, 
the computational limits of the technique is similar to the ones in the more traditional block models such as SBM and DCBM.
}
\end{rem}

%%%%%%%%%%%%%%%%%%%%%%%%%%%%%%%%%%%%%%%%%%%%%%%%%%%%%%%%%%%%%%%%%%%%%%%%%%%%%%%%%%%%%%%%%%%%%%%%%%%%%%%%%%%%%%%%%%%%%%%%%%
%%%%%%%%%%%%%%%%%%%%%%%%%%%%%%%%%%%%%%%%%%%%%%%%%%%%%%%%%%%%%%%%%%%%%%%%%%%%%%%%%%%%%%%%%%%%%%%%%%%%%%%%%%%%%%%%%%%%%%%%%%

\subsection{Correctness of the SSC at population level}
\label{sec:SSC_correctness}

%%%%%%%%%%%%%%%%%%%%%%%%%%%%%%%%%%%%%%%%%%%%%%
%
\begin{algorithm}[t] 
\caption{\ Consistent noiseless SSC (\cite{pmlr-v51-wang16b}) }
% \caption{\ Clustering consistent noiseless SSC \cite{pmlr-v51-wang16b} }
\label{alg: SSC_clust}
\begin{flushleft} 
{\bf Input:} The noiseless data matrix $P_*$ \\
{\bf Output:} Vector of community assignments $\hat{c}$, clustering matrix $\hat{Z}$ and 
recovered subspaces  $\hat{\calS}_{k}$, $k=1, \ldots, K$  \\
{\bf Steps:}\\
% {\bf 1.\ Initialization:}  Normalize each column of $P$ so that it has unit two norm {\bf NO NEED!!}  \\
{\bf 1.\ Constructing the similarity graph:} Solve the optimization problem \eqref{mn:opt_prob1}  
and construct the similarity matrix $S$ defined in \eqref{eq:similarity}. %  = |\widehat{W}| + |\widehat{W}^{T}|$. 
Construct  the  similarity graph $G = (V;E)$ where $(i,j) \in E$ if and only if $S_{i,j} >0$.\\
{\bf 2.\ Subspace recovery:}  For each connected   component
$G_r = (V_r; E_r)$ of $G$, compute $\hat{\calS}_{(r)}  = \Range(P_{V_r})$ using any convenient linear algebraic
method. Let $\{\hat{\calS}^{(k)}\}_{k=1}^K$  be the $K$ unique subspaces
in $\hat{\calS}_{(r)}$.\\
{\bf 3. \ Final clustering:} For each connected component $V_r$ with $\hat{\calS}_{(r)} = \hat{\calS}^{(k)}$, set
$\hat{c}_i = k$, $\hat{Z}_{i,l} = \II(l=k)$ for $i \in V_r$ and $l=1, \ldots, K$. \\ 
\end{flushleft} 
\end{algorithm}
%
%%%%%%%%%%%%%%%%%%%%%%%%%%%%%%%%%%%%%%%%%%%%

In order to apply the SSC for clustering in the PABM, we need to show that the SSC can detect   
communities correctly, at least at the population level.
Let $P_*$ be the true probability matrix. Since $P_*$ contains no errors, one can obtain the 
 coefficients matrix $W$ as a solution of the optimization problem \eqref{mn:opt_prob1}. 
It turns out that Assumption {\bf A1*} in Section~\ref{sec:cl_detectability} guarantees 
the correct community assignment   (up to permutations of class labels).

Recovery of the clustering matrix $Z_*$ relies on the fact that each column of the matrix $P_*$
is represented as a linear combination of points in its own subspace rather than from a different one. 
This is formalized as the {\it Self-Expressiveness Property} ({\bf SEP}) of the similarity graph $S$:
$S_{i,j} >0$ implies that nodes $i$ and $j$ belong to the same cluster (see, e.g.   
\cite{Elhamifar:2013:SSC:2554063.2554078}). Note that the reverse is not necessarily true:
the fact that nodes $i$ and $j$ are in the same cluster does not necessarily imply that $S_{i,j} >0$.
Under Assumption A1*, Theorem 1 of  \cite{Elhamifar:2013:SSC:2554063.2554078} 
ensures that the similarity matrix $S$ obtained as a solution of 
optimization problem \eqref{mn:opt_prob1} satisfies the SEP.

Nevertheless,  the SEP alone does not lead to the perfect clustering
because the obtained similarity graph $G$ could be
poorly connected (see, e.g.,  \cite{ICPR_2011_Nasihatkon}). It appears however
that a simple post-clustering procedure (Algorithm 1) suggested in \cite{pmlr-v51-wang16b}, 
guarantees the correct recovery. % We present the algorithm below.

\begin{thm}  \label{th:correctness}
Assume that the correct number of communities $K=K_*$ is known. Then, under 
Assumption~A1*, Algorithm 1 recovers communities correctly up to a permutation, that is, there exists a permutation $\pi$ on 
$\{1, \ldots, K\}$ such that $\hat{Z}_{i,\pi(k)} = (Z_*)_{i,k}$ for every $i=1, \ldots, n$ and $k=1, \ldots, K$.
\end{thm}

%%%%%%%%%%%%%%%%%%%%%%%%%%%%%%%%%%%%%%%%%%%%%%%%%%%%%%%%%%%%%%%%%%%%%%%%%%%%%%%%%%%%%%%%%%%%%%%%%%%%%%%%%%%%%%%55
%%%%%%%%%%%%%%%%%%%%%%%%%%%%%%%%%%%%%%%%%%%%%%%%%%%%%%%%%%%%%%%%%%%%%%%%%%%%%%%%%%%%%%%%%%%%%%%%%%%%%%%%%%%%%%%55

\subsection{Accuracy of the data-based  SSC}
\label{sec:SSC-errors}

While there are many papers that evaluate clustering errors in the case of the k-means algorithm 
and spectral clustering, as well as  their relaxations, there are very few results available for the accuracy  
of the SSC, and those results are quite recent. 
As it is evident from Section~\ref{sec:SSC_correctness}, the successful clustering relies
on the fact that the SEP condition is satisfied with the high probability and that the similarity 
graph is sufficiently connected.

The main effort of the scientific community was devoted to establishing the SEP condition.
Initially, this effort was directed to its justification when the true matrix is measured without errors
\cite{Vidal:2009aa} or with  outliers \cite{soltanolkotabi2012}. The latter paper
assumes that the columns of the data matrix are generated at random using the bases of the 
respective sub-space (semi-random model). Furthermore,  the bases may themselves be generated  uniformly, at random,
from the  unit sphere (random model), and the outliers are also uniformly distributed on the unit sphere.
In the subsequent paper, \cite{soltanolkotabi2014} handle the case where the data matrix contains 
small Gaussian errors. Few later papers refine the results of the last two publications. 
Specifically, \cite{JMLR:v17:13-354} extend results of \cite{soltanolkotabi2014} to the case 
of the deterministic model with random noise, i.e., the model where true matrix is  fixed 
in advance and is not generated at random. The noise vectors are assumed to be i.i.d. spherically symmetric 
with the lengths bounded above by a small quantity. Results of the latter paper are also featured in a 
very recent monograph of \cite{Shi_2019} which refines results of \cite{JMLR:v17:13-354} 
in the case of missing observations or Gaussian noise.

The graph connectivity in the SSC has been much less studied. 
Indeed, it is mentioned in  \cite{ICPR_2011_Nasihatkon} that the similarity graph may 
satisfy the SEP condition but be poorly connected.
Fortunately, this issue can be addressed by   post-processing procedures suggested in  
 \cite{pmlr-v51-wang16b}. One of the procedures is presented in Algorithm 1.

Since our optimization problem guarantees that coefficient vector has $K$ non-zero components, 
and, hence, leads to a sufficiently well connected similarity graph, it is the SEP condition
that presents the hardest  challenge. Indeed, as it follows from the review above,
there are two types of derivations of the SEP conditions in the existing literature.
Specifically, in the case where $X \in \RR^{n \times N}$ is a true matrix with columns drawn from $K$ 
different subspaces    $\calS_{i},  i=1,...,K,$ and $Y = X + \Xi$\ is its noisy version, 
the papers differ on whether they treat elements of matrix $\Xi$ as deterministic or random.
In both cases, the procedures start with scaling columns of matrix $Y$ to the unit length.

The case of the random errors handles either  the normally distributed errors (\cite{Shi_2019}, \cite{soltanolkotabi2014}),
or, more broadly, i.i.d errors having a spherically symmetric distribution (\cite{JMLR:v17:13-354}). Moreover, this assumption
constitutes the cornerstone of the proofs since the arguments there are based 
on the fact that the errors are invariant under an orthogonal transformation.
In addition, for both random and deterministic errors (see, e.g., \cite{pmlr-v51-wang16b}), 
it is assumed  that for any columns $X_i$ and $\Xi_i$ of matrices
$X$ and $\Xi$, respectively, one has 
\be \label{eq:p-xi}
\|\Xi_i\| \leq \delta_n \| X_i \|, \quad i=1, \ldots, n, 
\ee 
where $\delta_n \to 0$ as $n \to \infty$ with high probability.
The latter implies that $\| X_i \| = 1 + o(1)$ as $n \to \infty$.

Note that neither of the above assumptions are true in the case of Bernoulli errors. 
It is easy to see that Bernoulli errors are not i.i.d and that
the columns of matrix $\Xi$ are not spherically symmetric.  
Moreover, assumption \eqref{eq:p-xi} is not true in the case of the   Bernoulli data.
To understand this, consider a vector $p \in [0,1]^n$, a vector $a$ of 
independent Bernoulli variables $a_i \sim \mbox{Bernoulli}(p_i)$, $i=1,\ldots,n$,
and   $\xi = a-p$. Then,  for $p_i \leq 1/2$ one has $p_i (1 - p_i) \geq p_i^2$ and
\bes
\EE \|\xi\|^2 = \sum_{i=1}^n (p_i - p_i^2) \geq   \|p\|^2,
\ees
so that    inequality \eqref{eq:p-xi} does not hold.

To make matters worse, scaling of the columns of the matrix $A$ in \eqref{eq:A_Bern_P} to unit length is  itself
problematic.  Indeed, since components of a Bernoulli vector are either zeros or ones, one has 
$\|a\|^2  = \sum_{i=1}^n a_i = a^T 1$, so $\|a\|^2$ is a very poor approximation of $\|p\|^2$.
As a matter of fact,  scaling column $A_i$ of matrix $A$ to unit length amounts to dividing this column 
by $\sqrt{d(i)}$ where $d(i)$ is the degree of the node $i$.

In conclusion, the existing error analysis for the SSC cannot be used in the case of 
Bernoulli data and one needs to establish a brand new theory. Development of such theory is a matter of 
future investigation.

%%%%%%%%%%%%%%%%%%%%%%%%%%%%%%%%%%%%%%%%%%%%%%%%%%%%%%%%%%%%%%%%%%%%%%%%%%%%%%%%%%%%%%%%%%%%%%%%%%%%%%%%%%%%%%%55
% This is the shortened version of the section
% The complete version of this section is in Clust_Nov21_2018.tex

\section{Simulations and   real data examples}
\label{sec:simul_real}
\setcounter{equation}{0}

%In this section we demonstrate the performance of PABM model by simulation and its performance on the real data examples.
%%%%%%%%%%%%%%%%%%%%%%%%%%%%%%%%%%%%%%%%%%%%%%%%%%%%%%%%%%%%%%%%%%%%%%%%%%

\subsection{Simulations on synthetic networks}
\label{sec:simulations}

\ignore{
\begin{figure}[h!]

\centerline{ \includegraphics[width = 18.3cm]{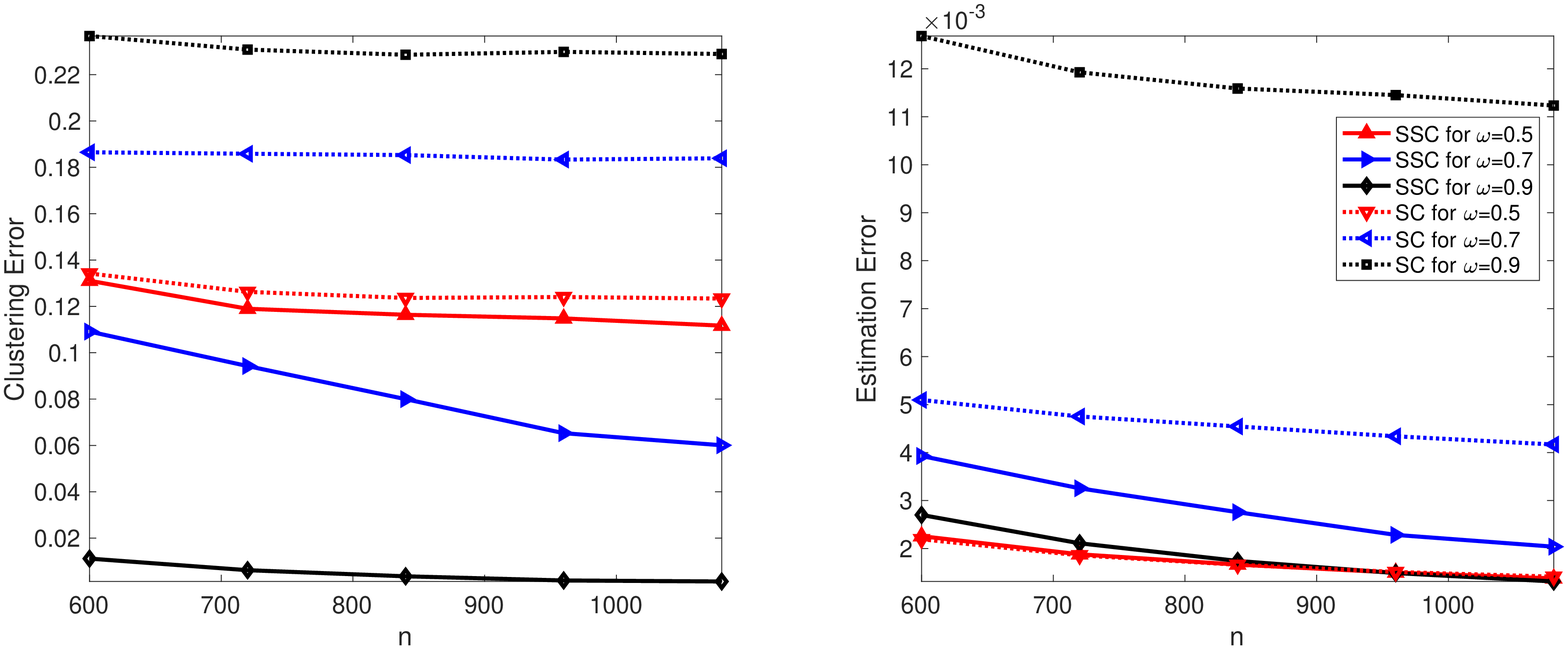}  }
\centerline{ \includegraphics[width = 18.3cm]{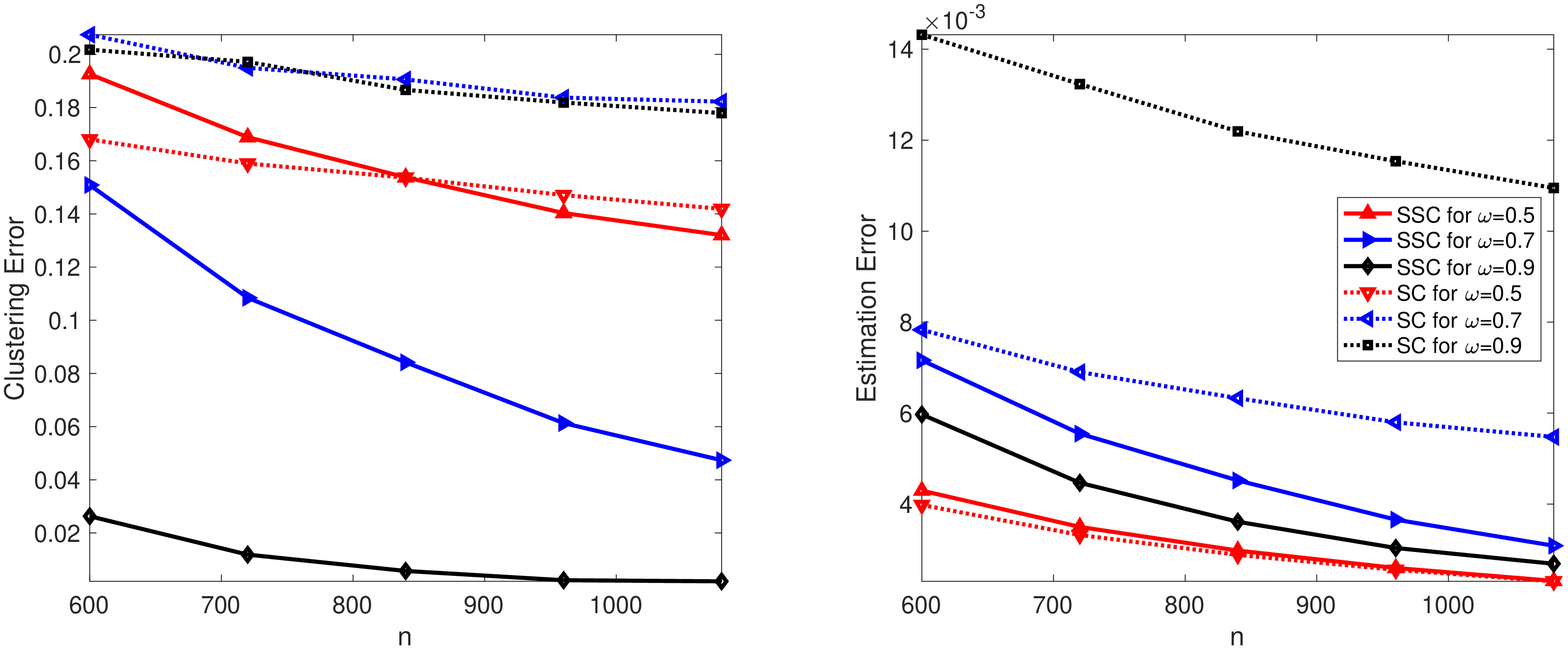}  }
\caption{ The clustering errors $\Err (\hat{Z},Z)$ defined in \eqref{eq:misclustered}  (left panels)
and the estimation errors $n^{-2}\, \|{\hat{P}-P}\|_{F}^{2}$ (right panels) for $K=4$ (top)  and 
$K=8$ (bottom) clusters. The errors are evaluated over 50 simulation runs. 
The number of nodes ranges from $n=600$ to $n=1080$ with the increments of 120. 
SSC results are represented by the solid lines; SC results are represented by the dotted lines:
$\om=0.5$  (red), $\om=0.7$ (blue) and $\om=0.9$ (black).
}
\label{mn:fig1}
\end{figure}
} % ignore

\ignore{
\begin{figure}[h!]

\centerline{ \includegraphics[width = 18.3cm]{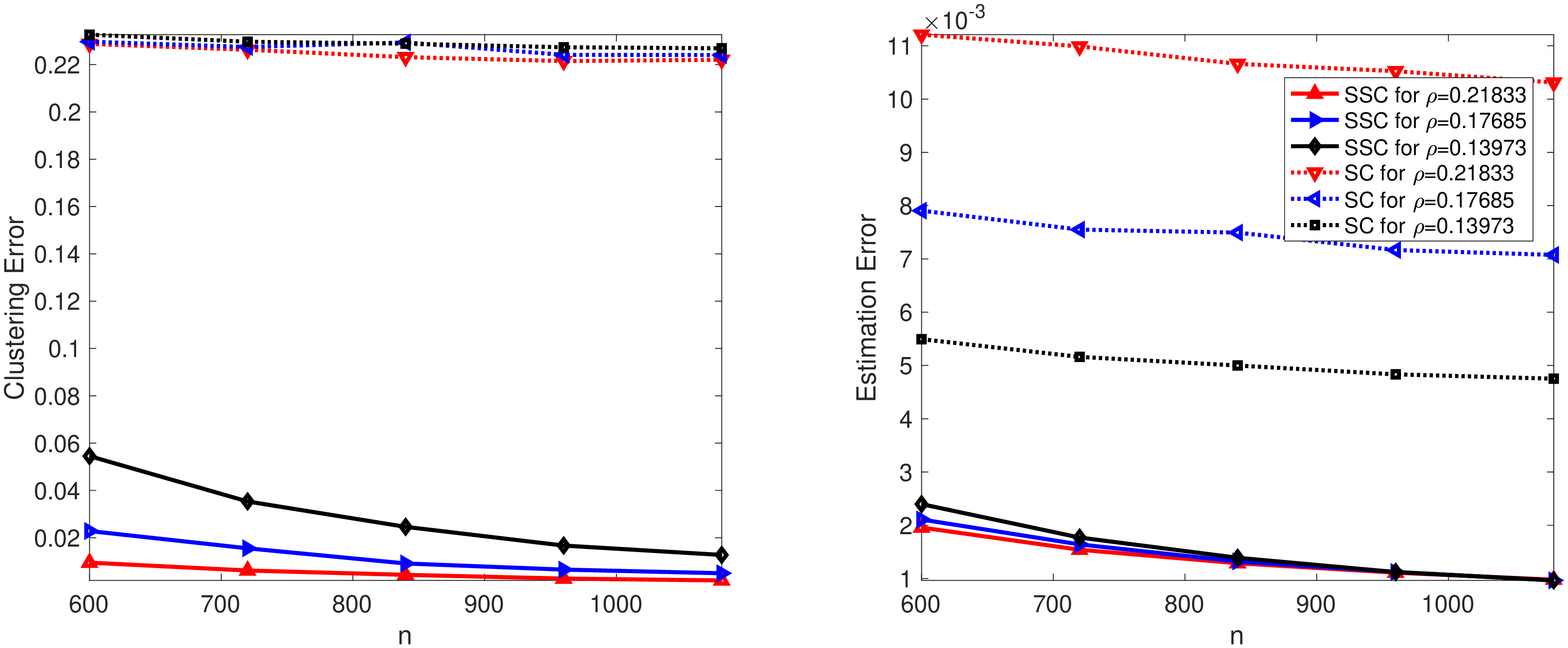}  }
\centerline{ \includegraphics[width = 18.3cm]{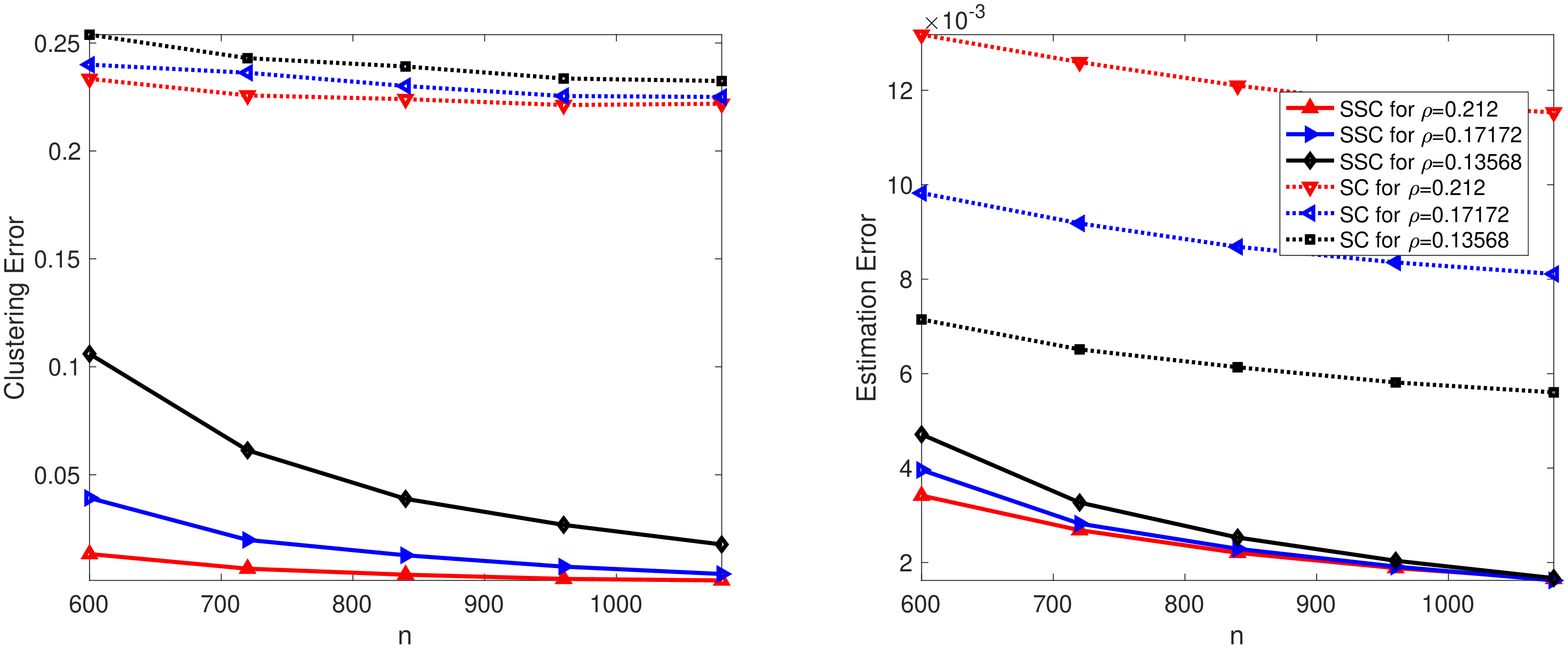}  }
\caption{ The clustering errors $\Err (\hat{Z},Z)$ defined in \eqref{eq:misclustered}  (left panels)
and the estimation errors $n^{-2}\, \|{\hat{P}-P}\|_{F}^{2}$ (right panels) for $K=3$ (top)  and 
$K=5$ (bottom) clusters. The errors are evaluated over 50 simulation runs. 
The number of nodes ranges from $n=600$ to $n=1080$ with the increments of 120. 
SSC results are represented by the solid lines; SC results are represented by the dotted lines:
for three different values of $\rho$ and fixed $\om=0.9$.
}
\label{mn:fig2}
\end{figure}
} % ignore

In this section we evaluate the performance of our method using synthetic networks. 
We assume that the number of communities (clusters) $K$ is known and for simplicity consider 
a perfectly balanced model with $n/K$ nodes in each cluster.
We generate each network  from a random graph model with a  symmetric probability  matrix $P$
given by the PABM model with a clustering matrix $Z$ and a block matrix $\Lambda$.

Sengupta and Chen (2018), in their simulations, considered networks with $K=2$ communities of equal sizes and matrices $\Lam$
in \eqref{eq:PABM-model} with elements 
$\Lambda_{i,r}=\alpha_i \sqrt{\frac{h}{1+h}}$ when node $i$ lies in class $r$, and $\Lambda_{i,r}=\beta_{i}\sqrt{\frac{1}{1+h}}$ 
otherwise, where $h$ is the homophily factor. The factors $\alpha_i$ and $\beta_i$ were set to 0.8 for half of the nodes 
in each class and to 0.2 for another half at random, and $h$ ranges between 1.5 and 4.0. 
Note that, although the data generated by the procedure above follows PABM, the probability matrix has constant blocks, 
for which the spectral clustering is known to deliver accurate results. In particular, the setting above  leads to  the SBM with four blocks.
However, the spectral clustering incurs some difficulties as the probabilities of connections in every community become more diverse. 
In this paper, we make sure to  generate networks that follow PABM with diverse probabilities of connections.

To generate a more diverse synthetic network, we start by producing a block matrix $\Lambda$ in \eqref{eq:Lambda} 
with random entries  on the interval $(0,a)$, $0<a<1$. 
We multiply the non-diagonal blocks of $\Lambda$ by $\omega$, $0<\omega <1$,  to ensure that most nodes in the same community have 
larger probability of interactions. Then matrix $P(Z,K)$ with blocks $P^{(k,l)}_{Z,K} =\Lambda^{(k,l)}  (\Lambda^{(l,k)})^{T}$,
$k,l = 1, \ldots, K$, mostly has larger entries in the diagonal blocks than in the non-diagonal blocks. 
The parameter $\om$ is the heterogeneity parameter. Indeed, if $\om=0$, 
the matrix $P_*$ is strictly block-diagonal, while in the case of $\om =1$, there is no difference between diagonal and non-diagonal blocks.
Next, we generate a random clustering matrix $Z \in \calM_{n,K}$ corresponding to the case of equal community sizes and the 
permutation matrix $\mathscr{P}(Z,K)$ corresponding to the clustering matrix $Z$. Subsequently, we scramble rows and columns 
of $P(Z, K)$ to create the probability matrix  $P=\scrPZK P(Z,K) \scrPZK^T$.  
Finally we generate the lower half of the adjacency matrix $A$ as independent Bernoulli variables  
$A_{ij} \sim \text{Ber}(P_{ij})$, $i=1, \ldots, n, j=1, \ldots, i-1$, and set  $A_{ij} = A_{ji}$ when $j >i$.
In practice, the diagonal $\text{diag}(A)$ of matrix $A$ is unavailable, so we estimate $\diag (P)$ 
without its knowledge. 

 \begin{figure}[h!]

\centerline{ \includegraphics[width = 16cm]{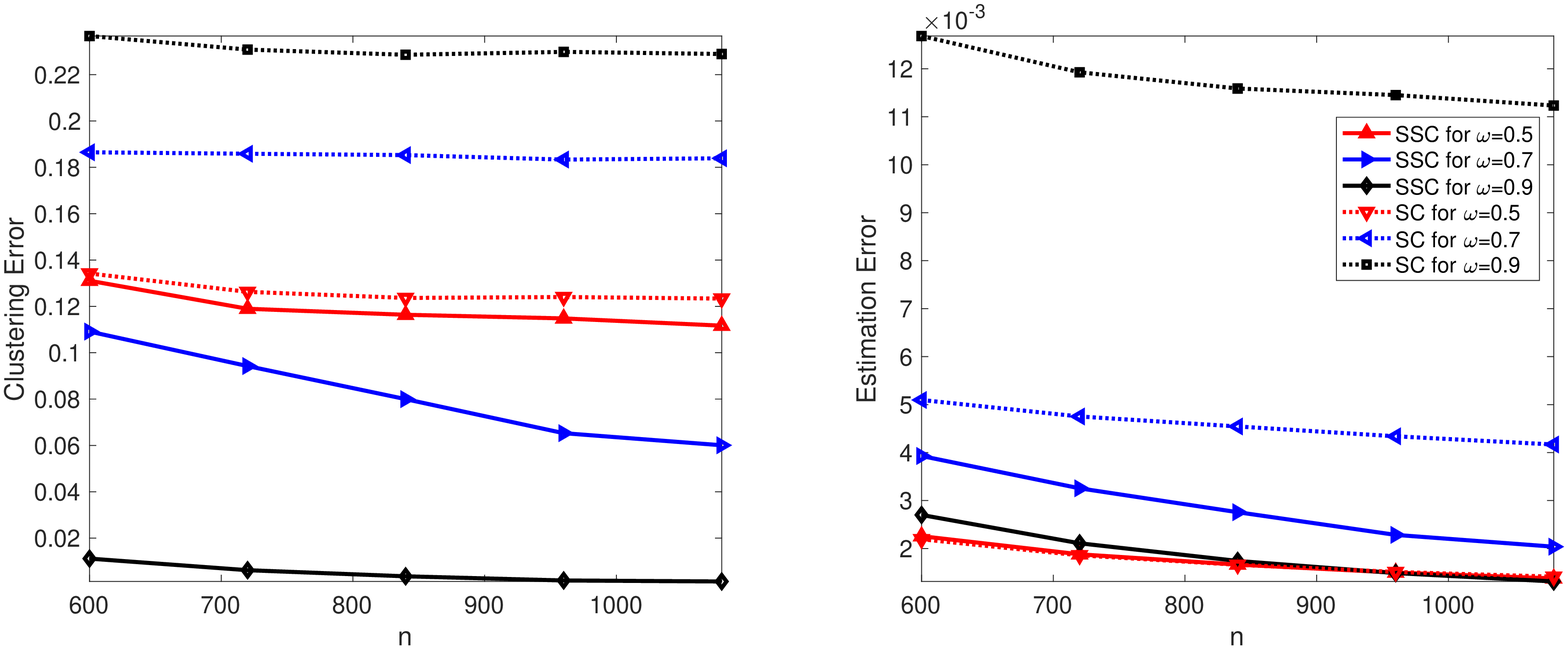}  }
\centerline{ \includegraphics[width = 16cm]{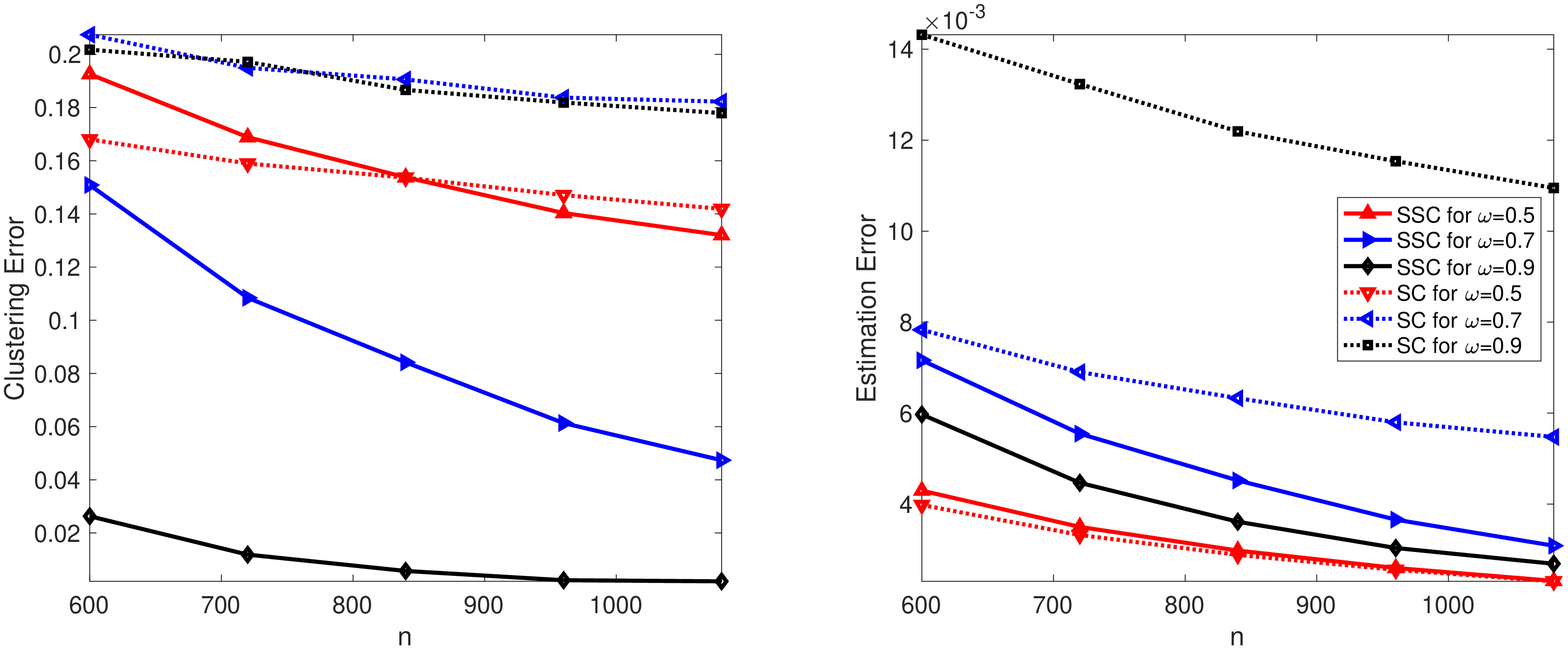}  }
\caption{ The clustering errors $\Err (\hat{Z},Z)$ defined in \eqref{eq:misclustered}  (left panels)
and the estimation errors $n^{-2}\, \|{\hat{P}-P}\|_{F}^{2}$ (right panels) for $K=4$ (top)  and 
$K=8$ (bottom) clusters. The errors are evaluated over 50 simulation runs. 
The number of nodes ranges from $n=600$ to $n=1080$ with the increments of 120. 
SSC results are represented by the solid lines; SC results are represented by the dotted lines:
$\om=0.5$  (red), $\om=0.7$ (blue) and $\om=0.9$ (black).
}
\label{mn:fig1}
\end{figure}

\begin{figure}[h!]

\centerline{ \includegraphics[width = 16cm]{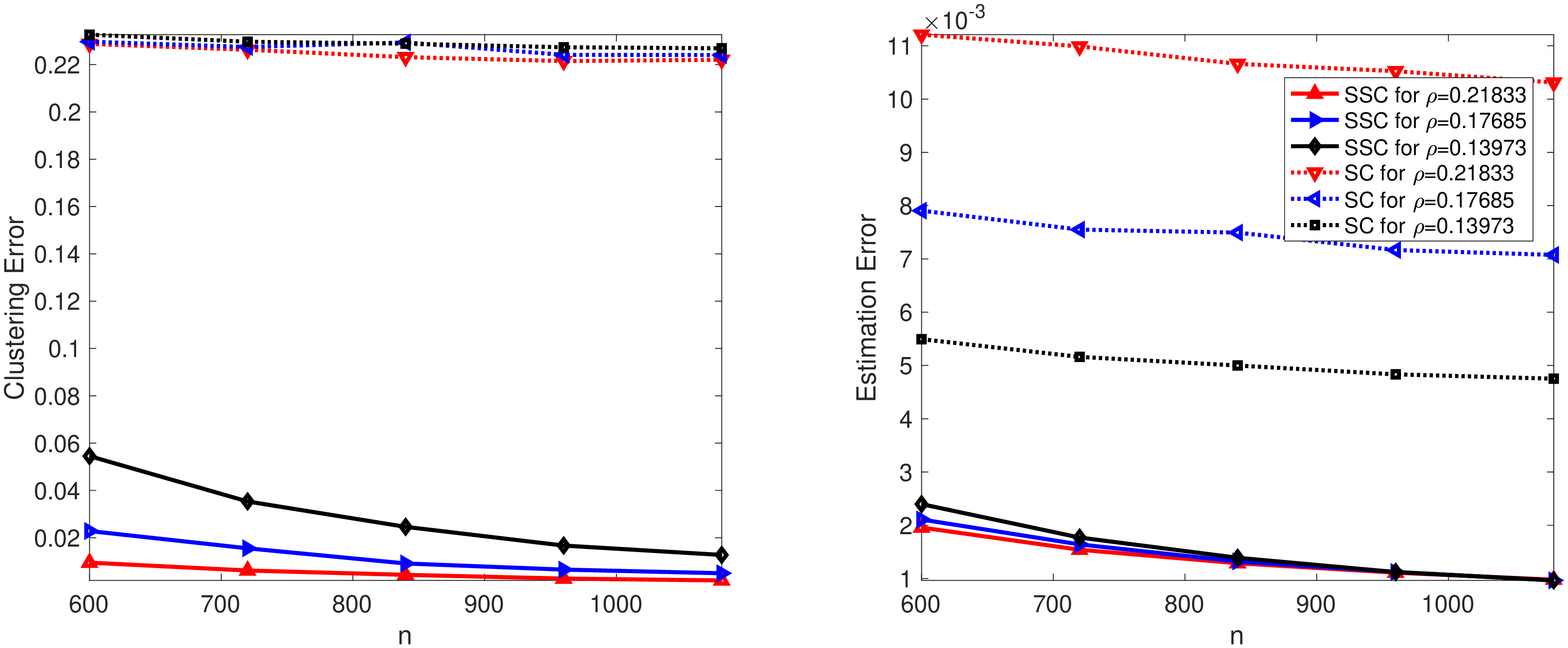}  }
\centerline{ \includegraphics[width = 16cm]{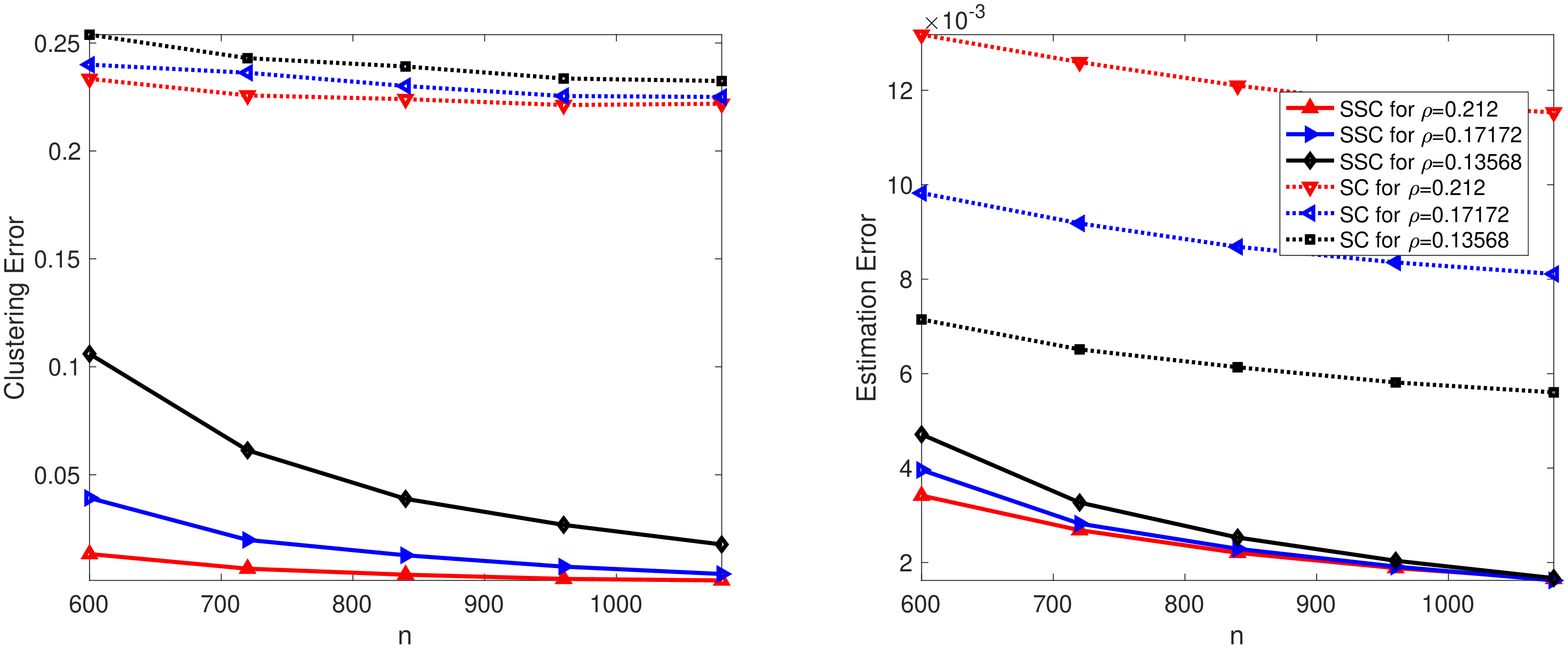}  }
\caption{ The clustering errors $\Err (\hat{Z},Z)$ defined in \eqref{eq:misclustered}  (left panels)
and the estimation errors $n^{-2}\, \|{\hat{P}-P}\|_{F}^{2}$ (right panels) for $K=3$ (top)  and 
$K=5$ (bottom) clusters. The errors are evaluated over 50 simulation runs. 
The number of nodes ranges from $n=600$ to $n=1080$ with the increments of 120. 
SSC results are represented by the solid lines; SC results are represented by the dotted lines:
for three different values of $\rho$ and fixed $\om=0.9$.
}
\label{mn:fig2}
\end{figure}

\cite{RePEc:bla:jorssb:v:80:y:2018:i:2:p:365-386} used  the   Extreme Points (EP) algorithm, introduced in  \cite{le2016},
 as a clustering procedure.  For $K = 2$, the EP algorithm computes the two leading eigenvectors of the adjacency
matrix $A$, and finds the candidate assignments associated with the extreme points of the
projection of the cube $[-1, 1]^n$ onto the space spanned by the two leading eigenvectors
of $A$. The technique is becoming problematic when $K$ grows and the probabilities of connections    are getting 
more diverse, hence, \cite{RePEc:bla:jorssb:v:80:y:2018:i:2:p:365-386} have only studied 
 performances of estimation and clustering  in  the case of $K=2$ and the choices of probability matrix 
$P$ described above.
As we have mentioned before, these are the settings  for which  the spectral clustering  procedure 
allows to identify the communities. Considering that we are interested in studying $K>2$ and the more diverse probabilities of connections, 
we use the spectral clustering directly (SC thereafter) and compare its precision with the sparse subspace clustering (SSC)
procedure.

Since the diagonal elements of matrix $A$ are unavailable,  we initially set $A_{ii}=0$, $i=1,...,n$. 
We solve optimization problem (\ref{mn:NP_Hard_2})  using the Orthogonal Matching Pursuit (OMP) algorithm.
% implemented  in  SPAMS MatLab toolbox (see \cite{mairal2014spams}).  
% 
After matrix $\widehat{W}$ of weights is  evaluated, we obtain the  clustering matrix
$\hat{Z}$ by applying spectral clustering to $|\widehat{W}| + |\widehat{W}^{T}|$, as it was described in Section~\ref{sec:SSC_review}.
Given $\hat{Z}$, we generate matrix $A(\hat{Z}) = \scrP_{\hat{Z}}^T A \scrP_{\hat{Z}}$ with blocks $A^{(k,l)}(\hat{Z})$,
$k,l=1,\ldots, K$, and obtain $\hat{\Theta}^{(k,l)}(\hat{Z}, \hat{K})$ by using the rank one approximation for each of the blocks. 
Finally, we  estimate  matrix $P$ by $\hat{P} =  \hat{P}(\hat{Z}, \hat{K})$ using formula~\eqref{eq:P_total_est} with $\hat{K}=K$.

\begin{figure}[h!]
\centerline{ \includegraphics[width=4.5in ]{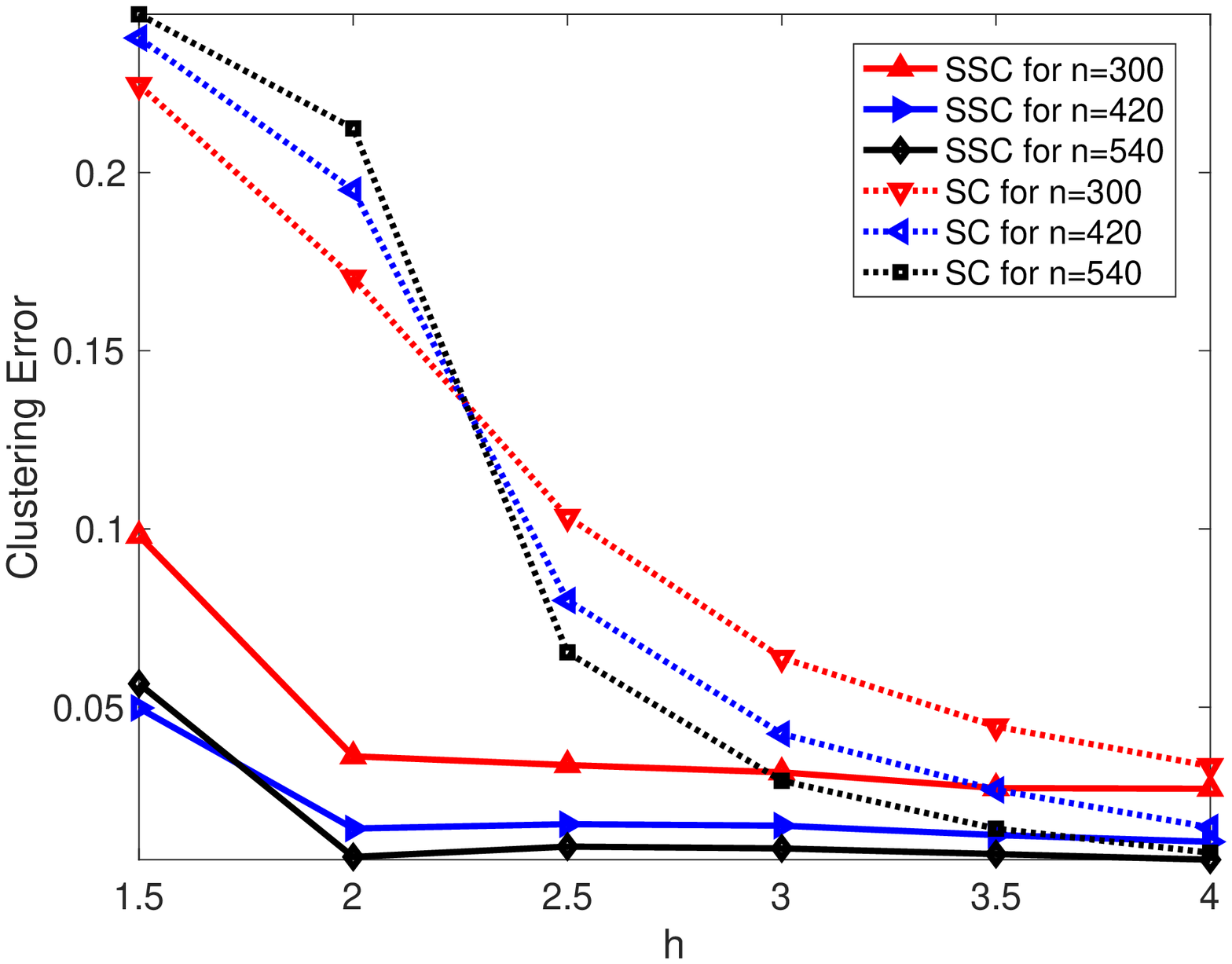}}%{San_Chen_Plot} }%{Sengupta_and_Chen} %{Sen_and_Chen} 
\caption{Clustering errors of SC and SSC for $K=2$ clusters and $n = 300, 420$ and 540 nodes in the simulations setting of 
% Sengupta and Chen
 \cite{RePEc:bla:jorssb:v:80:y:2018:i:2:p:365-386}.  The homophily factor $h$  ranges from 1.5 to 4 with increments of 0.5
}
\label{mn:fig4}
\end{figure}

%%%%%%%%%%%%%%%  Simulations results %%%%%%%%%%%%%%%%%%%%%%%%%%%%%%%%%%%%%%

We compared the accuracy of SSC and SC methods in terms of 
the average estimation errors $n^{-2}\, \|{\hat{P}-P}\|_{F}^{2}$
and the average clustering errors $\Err (\hat{Z},Z)$ defined in \eqref{eq:misclustered}.
Figures~\ref{mn:fig1}~and~~\ref{mn:fig2}  show  the  results of these comparisons for 
% $K=4$ and $K=8$, and figure~\ref{mn:fig2} shows the  results  for $K=3$ and $K=5$, and 
the number of nodes ranging from $n=600$  to $n =1080$ with the increments of 120.
The left panels  display the clustering errors $\Err (\hat{Z},Z)$ 
while the right ones exhibit the estimation errors $n^{-2}\, \|{\hat{P}-P}\|_{F}^{2}$,  as   functions of the number of nodes.
All errors are averaged over 50 simulation runs. Figure~\ref{mn:fig1} explores the effects of 
heterogeneity on the precision of estimation and clustering by carrying out simulations for  $K=4$ and $K=8$, and  
for three different values of the parameter $\omega$: $\omega =0.5$, 0.7, and 0.9.  
Block matrix $\Lambda$ in \eqref{eq:Lambda} has random entries   the interval $(0,1)$ in this case. 
Figure~\ref{mn:fig1}  confirms  that the   SSC  is becoming more and more  accurate in comparison with SC as $\omega$ grows.
The latter is due to the fact that the SSC is more suitable for handling heterogeneous connections probabilities.

Figure~\ref{mn:fig2} examines the impact of sparsity on clustering precision. 
For this round of simulations we used $K=3$ and $K=5$, fixed  $\om=0.9$ and generated entries of the block 
matrix $\Lambda$ in \eqref{eq:Lambda} on the intervals (0,1), (0,0.9) and (0,0.8).
The results are sorted by the average connection probability $\rho$
which, in the above simulation setting, takes values 0.21833, 0.17685 and 0.13973
for $K=3$, and 0.21200, 0.17172 and 0.13568 for $K=5$. Figure~\ref{mn:fig2} shows that 
the clustering errors decrease  as $\rho$ and $n$ increase, with the effect 
of growth of $n$ on the accuaracy of clustering  being much more significant in the case of the SSC.

Figure \ref{mn:fig4} presents the results of comparison of the clustering errors of SSC and SC in the simulations settings of 
%  Sengupta and Chen
 \cite{RePEc:bla:jorssb:v:80:y:2018:i:2:p:365-386}. It is easy to see that, while for larger values of the homophily factor $h$ 
both methods perform almost equally well, the accuracy of SC deteriorates as $h$ is getting smaller, due to the fact that the 
differences between probabilities of connections within and between clusters become less significant. The latter shows that 
the SSC approach is beneficial for clustering in PABM model. Indeed, it delivers more accurate  results than the SC 
when probabilities of connections are more diverse.  On the other hand, SSC is still applicable when the PABM reduces to the SBM,
although SC is more accurate in the case of the SBM since it does not require an additional step of evaluating the affinity matrix.

%%%%%%%%%%%%%%%%%%%%%%%%%%%%%%%%%%%%%%%%%%%%%%%%%%%%%%%%%%%%%%%%%%%%%%%%%%%%%%%%%%%%%%%%%

\begin{rem}\label{rem:SC_versus_SSC}
{\bf Spectral Clustering Versus Sparse Subspace Clustering. }{\rm
It is worth noting  that when the  matrix of probabilities $P_*$ is close to being block diagonal,
the spectral clustering can be still used for recovering community assignments, even if $P_*$ does not follow the SBM. 
The latter is due to the fact that, in this situation,  the graph can be well approximated by a union of distinct connected components,
and, therefore, SC allows  to identify the true clusters. Moreover, in such situation, SC has an advantage of not requiring 
an additional step of self-representation, which is computationally costly and produces additional errors. 
On the other hand, as we shall see from examples below, when   probabilities of connections become more heterogeneous,
SSC turns to be more precise than SC.
In addition, since  PABM has more unknown parameters than SBM, its correct fitting requires
sufficient number of nodes per class (see, e.g., \cite{soltanolkotabi2014});
otherwise, its accuracy declines. 
}
\end{rem}

%%%%%%%%%%%%%%%%%%%%%%%%%  Unknown K  %%%%%%%%%%%%%%%%%%%%%%%%%%%%%%%%%%%%%%%%%%%%%%%%%%%%

\begin{rem}\label{rem:Unknown_K}
{\bf Unknown number of clusters. }{\rm
In our previous simulations we treated  the true number of clusters as a known quantity. 
However, we can actually use $\hat{P}$ to obtain an   estimator $\hat{K}$  of $K$ 
by solving, for every suitable $K$, the   optimization problem  \eqref{eq:opt_for_K}, 
which can be equivalently rewritten as 
\be  \label{mn:Unknown_K}
\hat{K}=\mathop\text{argmin}_{K} \{ \|{\hat{P}-A}\|_{F}^{2} + \Pen (n,K)\}.
\ee
The penalty  $\Pen (n,K)$ defined in \eqref{eq:penalty} is, however, motivated by 
the objective of setting it above the noise level with a very high probability. 
In our simulations, we also study the selection of an unknown $K$  using somewhat smaller penalty 
\be \label{eq:smaller_pen}
\Pen(n,K)=\rho(A) n K \sqrt{\ln n\, (\ln K)^{3}}
\ee
where $\rho(A)$ is the density of matrix $A$,   the proportion of  nonzero entries of $A$. 

In order to assess the accuracy of  $\hat{K}$ as an estimator of $K$,  we evaluated 
$\hat{K}$ as a solution of optimization problem \eqref{mn:Unknown_K} with the penalty \eqref{eq:smaller_pen}
in each of the previous simulations settings over 50 simulation runs.
Table~\ref{mn:table1} presents the   relative frequencies of the estimators $\hat{K}$ of $K_*$ for $K_*$ 
ranging from 3 to 6,   $n=420$ and $n =840$ and $\om = 0.5$, 0.7 and 0.9.
Table~\ref{mn:table1} confirms that for majority of settings, $\hat{K} = K_*$, the true number of clusters, 
with high probability.  Moreover, the estimator   $\hat{K}$ of $K$ is more reliable for higher values of $\om$
and larger number of nodes per cluster.
}
\end{rem}

%%%%%%%%%%%%%%%%%%%%%%%%%%%%%%%%%%%%%%%%%%%%%%%%%%%%%%%%%%%%%%%%%%%%%%%%%%%%%%%%%%%%%%%%%%%%%%%%%%%%%%%%%%%%%%%%%%%%%%%%%%%%%%
%%%%%%%%%%%%%%%%%%%%%%%%%%%%%%%%%%%%%%%%%%%%%%%%%%%%%%%%%%%%%%%%%%%%%%%%%%%%%%%%%%%%%%%%%%%%%%%%%%%%%%%%%%%%%%%%%%%%%%%%%%%%%%

\begin{table} [H]
%\begin{center}
\centering
\caption{The relative frequencies of the estimators $\hat{K}$ of $K_*$ for $K_*$ ranging from 3 to 6,   $n=420$ and $n =840$
and $\om = 0.5$, 0.7 and 0.9.}
%\scalebox{0.9}{
%\resizebox{0.5\textwidth}{
%\resizebox{\columnwidth}{
%\resizebox{0.7}{0.48}{
\begin{tabular}{|c| c|c |c| c| c|c|c| }
\multicolumn{8}{ l }{  }\\
\hline \hline 

& & \multicolumn{3}{ |c |}{n=420} & \multicolumn{3}{c|}{n=840}\\
\hline
$K_*$ & $\hat{K}$ & $\omega=0.5$   & $\omega=0.7$ & $\omega=0.9$   & $\omega=0.5$   & $\omega=0.7$ & $\omega=0.9$   \\
\cline{1-8} 
    & 2   & 0   & 0   & 0   & 0   & 0   &   0    \\  \cline{2-8}
    & 3   & \textbf{0.76}  & \textbf{0.80}   &  \textbf{0.90}  & \textbf{0.52}   &  \textbf{0.60}  & \textbf{0.80}     \\   \cline{2-8}     
 \textbf{3}  & 4   &  0.24  &  0.16  & 0.10   & 0.36   & 0.26   &   0.16    \\   \cline{2-8}
    & 5   & 0   & 0.04   &  0  &  0.12  & 0.14   &  0.02     \\     \cline{2-8}       
    & 6   & 0   & 0   &  0  & 0   & 0   &  0.02     \\      
\cline{2-8}    
\hline
\hline
\multicolumn{8}{ l }{  }\\
\hline \hline 

& & \multicolumn{3}{ |c |}{n=420} & \multicolumn{3}{c|}{n=840}\\
\hline
$K_*$ & $\hat{K}$ & $\omega=0.5$   & $\omega=0.7$ & $\omega=0.9$   & $\omega=0.5$    & $\omega=0.7$ & $\omega=0.9$   \\\cline{1-8}
 
    & 2   & 0   & 0   & 0   &  0  &  0  &  0     \\  \cline{2-8}
    & 3   & 0.06   & 0.14    & 0   & 0.02   &  0.02  & 0      \\   \cline{2-8}     
 \textbf{4}  & 4   & \textbf{0.64}   & \textbf{0.66}   &  \textbf{0.96}  & \textbf{0.56}   & \textbf{0.64}   & \textbf{0.76}      \\   \cline{2-8}
    & 5   & 0.28   & 0.16    &  0.04  &  0.30  & 0.26  &  0.22     \\     \cline{2-8}       
    & 6   & 0.02   &  0.04  & 0   & 0.12    & 0.08   & 0.02      \\      
\cline{2-8}      
\hline
\hline
\multicolumn{8}{ l }{  }\\
\hline \hline
& & \multicolumn{3}{ |c |}{n=420} & \multicolumn{3}{c|}{n=840}\\
\hline
$K_*$ & $\hat{K}$ & $\omega=0.5$   & $\omega=0.7$ & $\omega=0.9$   & $\omega=0.5$    & $\omega=0.7$ & $\omega=0.9$   \\ \cline{1-8}
 
    & 2   &  0  & 0.02   & 0   &  0  &  0  &  0     \\  \cline{2-8}
    & 3   &  0.02  &  0  & 0.02   & 0  & 0   & 0      \\   \cline{2-8}     
 \textbf{5}  & 4   & 0.14   & 0.16   & 0.04   & 0.04   &  0.04  &  0     \\   \cline{2-8}
    & 5   & \textbf{0.64}   &  \textbf{0.66}  &  \textbf{0.82}  &  \textbf{0.78}  & \textbf{0.68}  &  \textbf{0.90}     \\     \cline{2-8}       
    & 6   & 0.20   & 0.16   &  0.12  & 0.18   &  0.28  &  0.10     \\      
\cline{2-8}      
\hline
\hline
\multicolumn{8}{ l }{  }\\
\hline \hline 
& & \multicolumn{3}{ |c |}{n=420} & \multicolumn{3}{c|}{n=840}\\
\hline
$K_*$ & $\hat{K}$ & $\omega=0.5$   & $\omega=0.7$ & $\omega=0.9$   & $\omega=0.5$    & $\omega=0.7$ & $\omega=0.9$   \\\cline{1-8}
 
    & 2   &  0  &  0.04  & 0   &  0  & 0   &  0     \\  \cline{2-8}
    & 3   & 0.06   &  0.18  &  0.02  &  0  & 0   & 0      \\   \cline{2-8}     
 \textbf{6}  & 4   &  0.18  & 0.22   & 0.02   &  0  & 0   &  0     \\   \cline{2-8}
    & 5   &  0.28  & 0.22   &  0.08  &  0.12  & 0.16  & 0.10      \\     \cline{2-8}       
    & 6   &  \textbf{0.48}  & \textbf{0.34}   &  \textbf{0.88}  & \textbf{0.88}   &  \textbf{0.84}  &  \textbf{0.90}     \\      \cline{2-8}      
\hline
\hline
\end{tabular} 
%\end{center}
\label{mn:table1}
\end{table}

%%%%%%%%%%%%%%%%%%%%%%%%%%%%%%%%%%%%%%%%%%%%%%%%%%%%%%%%%%%%%%%%%%%%%%%%%%%%%%%%%%%%%%%%%%%%%%%%%%%%%%%%%%%%%%%%%%%%%%%%%%%%%%
%%%%%%%%%%%%%%%%%%%%%%%%%%%%%%%%%%%%%%%%%%%%%%%%%%%%%%%%%%%%%%%%%%%%%%%%%%%%%%%%%%%%%%%%%%%%%%%%%%%%%%%%%%%%%%%%%%%%%%%%%%%%%%

\ignore{
%%%%%%%%%%%%%%%%%%%%%%%%%  Unknown K  %%%%%%%%%%%%%%%%%%%%%%%%%%%%%%%%%%%%%%%%%%%%%%%%%%%%

\begin{rem}\label{rem:Unknown_K}
{\bf Unknown number of clusters. }{\rm
In our previous simulations we treated  the true number of clusters as a known quantity. 
However, we can actually use $\hat{P}$ to obtain an   estimator $\hat{K}$  of $K$ 
by solving, for every suitable $K$, the   optimization problem  \eqref{eq:opt_for_K}, 
which can be equivalently rewritten as 
\be  \label{mn:Unknown_K}
\hat{K}=\mathop\text{argmin}_{K} \{ \|{\hat{P}-A}\|_{F}^{2} + \Pen (n,K)\}.
\ee
The penalty  $\Pen (n,K)$ defined in \eqref{eq:penalty} is, however, motivated by 
the objective of setting it above the noise level with a very high probability. 
In our simulations, we also study the selection of an unknown $K$  using somewhat smaller penalty 
\be \label{eq:smaller_pen}
\Pen(n,K)=\rho(A) n K \sqrt{\ln n\, (\ln K)^{3}}
\ee
where $\rho(A)$ is the density of matrix $A$,   the proportion of  nonzero entries of $A$. 

In order to assess the accuracy of  $\hat{K}$ as an estimator of $K$,  we evaluated 
$\hat{K}$ as a solution of optimization problem \eqref{mn:Unknown_K} with the penalty \eqref{eq:smaller_pen}
in each of the previous simulations settings over 50 simulation runs.
Table~\ref{mn:table1} presents the   relative frequencies of the estimators $\hat{K}$ of $K_*$ for $K_*$ 
ranging from 3 to 6,   $n=420$ and $n =840$ and $\om = 0.5$, 0.7 and 0.9.
Table~\ref{mn:table1} confirms that for majority of settings, $\hat{K} = K_*$, the true number of clusters, 
with high probability.  Moreover, the estimator   $\hat{K}$ of $K$ is more reliable for higher values of $\om$
and larger number of nodes per cluster.
}
\end{rem}

}

%%%%%%%%%%%%%%%%%%%%%%%%%%%%%%%%%%%%%%%%%%%%%%%%%%%%%%%%%%%%%%%%%%%%%%%%%%%%%%%%%%%%%%%%%%%%%%%%%%%%%%%%%%%%%%%%%%%%%%%%%%%%%%
%%%%%%%%%%%%%%%%%%%%%%%%%%%%%%%%%%%%%%%%%%%%%%%%%%%%%%%%%%%%%%%%%%%%%%%%%%%%%%%%%%%%%%%%%%%%%%%%%%%%%%%%%%%%%%%%%%%%%%%%%%%%%

\subsection{Real data examples}
\label{sec:real_data} 

\begin{figure}[h!]
%  \begin{center}
\[ \includegraphics[height=6.0cm] {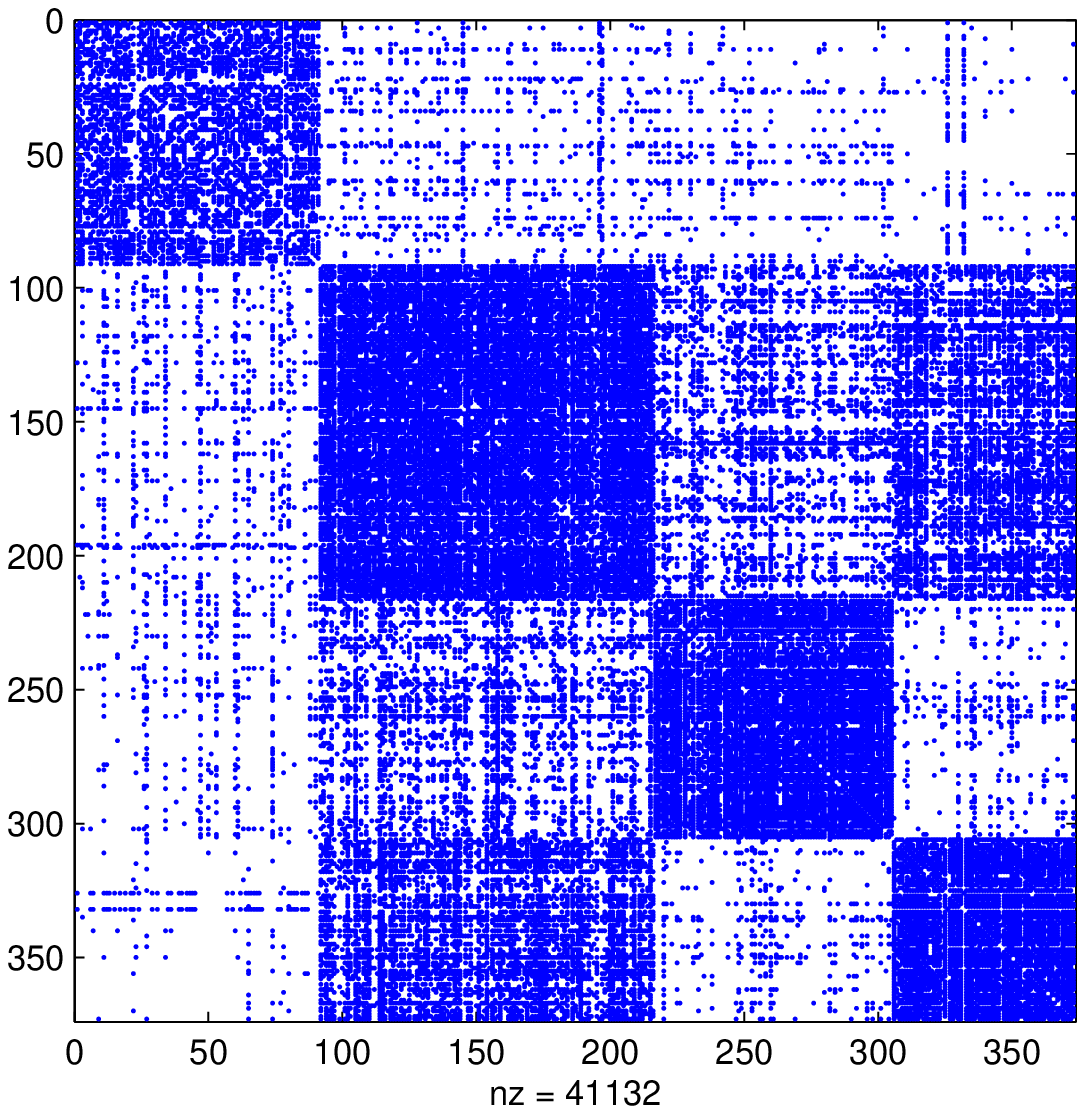} \hspace{4mm} 
\includegraphics[height=6.0cm]{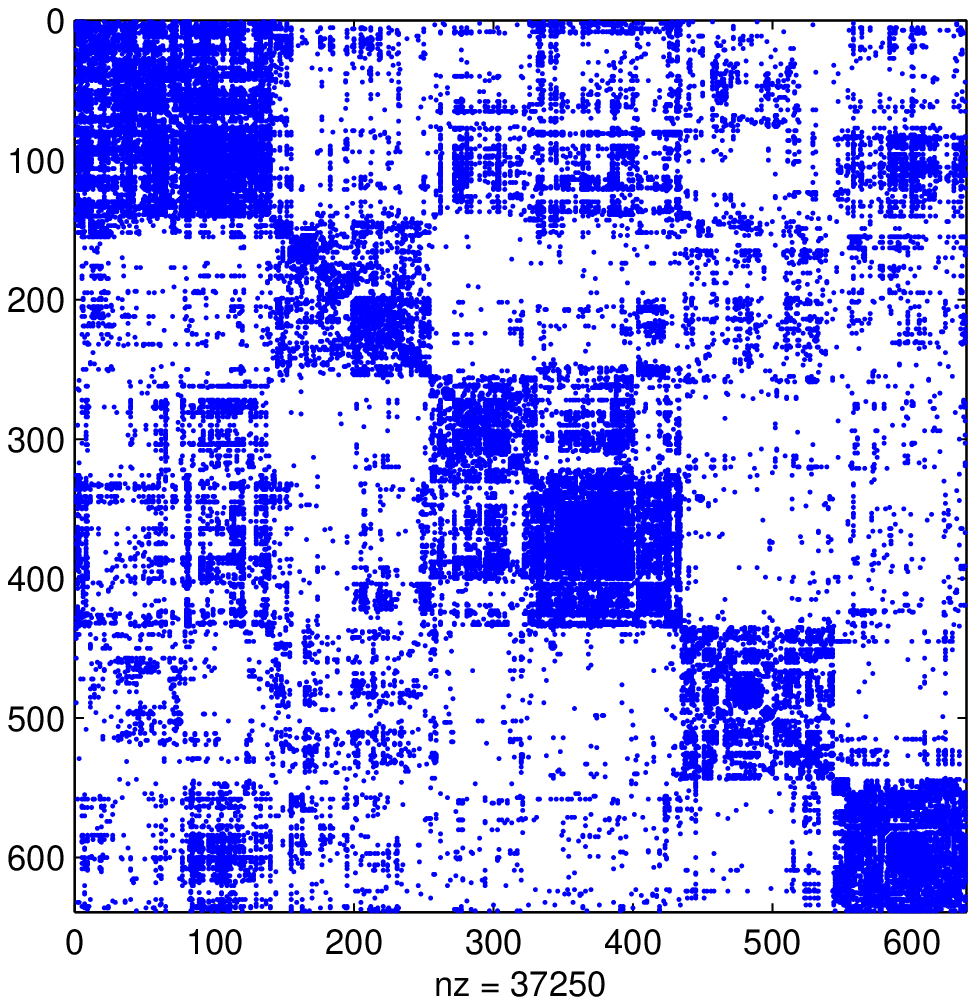} 
\] 
%    \subfigure{\label{fig:edge-a}\includegraphics[scale=0.29]{Butterfly_K_2}}
%    \subfigure{\label{fig:edge-b}\includegraphics[scale=0.29]{Butterfly_K_4}}   \\
%    \end{center}
  \caption{Adjacency matrices of the butterfly similarity network with 41132 nonzero 
entries and 4 clusters (left) and the brain network  with 37250 nonzero entries and 6 clusters (right)}
  \label{mn:fig_Adjacency matrices}
\end{figure}

In   this section, we report the performances of SSC and SC in studying real life networks. 
The social networks usually exhibit strong assortative behavior, the phenomenon which  
 is possibly due to the tendency of humans to form strong associations.
Perhaps, for this reason,   the political blogs network, the British Twitter network, 
and the Digital Bibliography \& Library Project network which have been analyzed by Sengupta and Chen (2018) have nearly block-diagonal 
adjacency matrices, so SC exhibits good performance in clustering of those networks (see Remark~\ref{rem:SC_versus_SSC}).

However, PABM provides a more accurate description of more diverse networks,
in particular, the networks that appear in biological sciences. Below, we  consider  
a butterfly similarity network extracted from the Leeds Butterfly  dataset described in % Wang {\it et al.} 
\cite{wang2018network}. 
Leeds Butterfly dataset  contains fine-grained  images of 832 butterfly species that belong to 10 different classes, 
with each class containing between 55 and 100 images. In this network, the nodes represent butterfly species 
and edges represent visual similarities between them. Visual similarities are evaluated on the basis of 
butterfly images and range from 0 to 1. 
We study a network  by extracting the four largest     classes as a simple graph  with  373 nodes and 20566 edges.
We draw an edge between the nodes if  the visual similarity between those nodes is greater than zero.
\ignore{
We carried out clustering of the nodes using the SSC, the SC, and the weighted $k$-median algorithm, one of the popular clustering methods for the  DCBM   
used in  \cite{lei2015} and \cite{gao2018community}. We compared the clustering assignments of those methods with the true 
class specifications of the species using the  adjusted Rand index  that measures the agreement between two clustering assignments. 
The value  of the  adjusted Rand index between the true 
class specifications and the clustering assignments obtained by  the SSC is 0.73; the weighted $k$-median algorithm is 0.67; and the SC is 0.61. 
%The SSC and the weighted $k$-median algorithm provide 89\% and 86\% accuracy, respectively, while SC is correct only in 64\% of cases. 
In addition, we applied formula (\ref{mn:Unknown_K}) with $K$ ranging from 2 to 6 and obtained the true number of clusters. % with 100\% accuracy.  
}

Classification of species on the basis of their visual similarities is a very important task.
In many applications, the  goal is to classify species automatically on the basis of their images,   captured by  a remote camera.
This type of monitoring is essential for surveying bio-diversity and tracking abundance and habitats 
of species that may be affected by climate change and human activities. 
While the related species may look similar, classification of species does not allow mixed memberships:
each of the actual species belong to one and only one class.

Figure \ref{mn:fig_Adjacency matrices} (left) shows the adjacency matrix of the graph (after clustering),
which suggests that the PABM is a reasonable model to fit to the network.  The latter is due to the fact that, 
 since the phenotype of the species in the same class can vary, the SBM may not 
provide an adequate summary for the class similarities.  Replacing the SBM by the DCBM does 
not solve the problem either, since it is unlikely that few butterflies are ``more similar'' to the others than the rest.
On the other hand, the PABM allows some of the butterflies in one class to be ``more similar'' 
to species of another specific class than the others, thus, justifying application of the PABM.

We carried out clustering of the nodes using the SSC, the SC, and the weighted $k$-median algorithm, one of the popular clustering methods for the  DCBM   
used in  \cite{lei2015} and \cite{gao2018community}. We compared the clustering assignments of those methods with the true 
class specifications of the species using the  adjusted Rand index  that measures the agreement between two clustering assignments. 
The value  of the  adjusted Rand index between the true 
class specifications and the clustering assignments obtained by  the SSC is 0.73; the weighted $k$-median algorithm is 0.67; and the SC is 0.61. 
%The SSC and the weighted $k$-median algorithm provide 89\% and 86\% accuracy, respectively, while SC is correct only in 64\% of cases. 
In addition, we applied formula (\ref{mn:Unknown_K}) with $K$ ranging from 2 to 6 and obtained the true number of clusters. % with 100\% accuracy.  

As the second real network, we analyze  a human brain functional network, measured using the resting-state functional MRI (fMRI). 
In particular, we use the co-activation matrix of  the brain connectivity dataset, described in % Crossley {\it et al.} 
\cite{crossley2013cognitive}. 
In this dataset, the brain is partitioned into  638 distinct regions and a weighted graph is used to characterize the network topology.
In our analysis, we set all nonzero weights to one, obtaining the network with  18625 undirected edges.  Since, for this network, 
the true clustering  as well as the true number of clusters are  unknown, we first applied formula (\ref{mn:Unknown_K}) 
with $K$ ranging from 2 to 10 to find the number of clusters obtaining  $\hat{K} = 6$. This agrees with the assessment in 
\cite{crossley2013cognitive}  where the authors partitioned the network into 6 groups (if one considers the ``rich-club''
communities as separate clusters).
Subsequently, we applied the SSC for partitioning the network into blocks and derived the estimator $\hat{P}$ of $P_*$.  
Figure \ref{mn:fig_Adjacency matrices} (right) shows the adjacency matrix of the graph after clustering.
The true probability matrix $P_*$ is unknown, we can only report that $ n^{-2}\, \|{\hat{P}-A}\|_{F}^{2} = 0.05$,
which indicates high agreement between the two matrices. 
We also carried out clustering of the nodes using the weighted $k$-median algorithm and the SC, that correspond, respectively, to modeling 
$P$ via the DCBM and the SBM,  and calculated the adjusted Rand index between the clustering assignments obtained by the three clustering methods. 
The adjusted Rand index between the clustering assignments obtained by  the SSC and the weighted $k$-median algorithm is 0.47; 
the SSC and the SC is 0.64; and the weighted $k$-median algorithm and the SC is 0.51. 

%%%%%%%%%%%%%%%%%%%%%%%%%%%%%%%%%%%%%%%%%%%%%%%%%%%%%%%%%%%%%%%%%%%%%%%%%%%%%%%%%%%%%%%%%%%%%%%%%%%%%%%%%%%%%%%55
%%%%%%%%%%%%%%%%%%%%%%%%%%%%%%%%%%%%%%%%%%%%%%%%%%%%%%%%%%%%%%%%%%%%%%%%%%%%%%%%%%%%%%%%%%%%%%%%%%%%%%%%%%%%%%%55

 \section*{acknowledgements}
All three authors of the paper were  partially supported by National Science Foundation
(NSF)  grant DMS-1712977.

%%%%%%%%%%%%%%%%%%%%%%%%%%%%%%%%%%%%%%%%%%%%%%%%%%%%%%%%%%%%%%%%%%%%%%%%%%%%%%%%%%%%%%%%%%%%%%%%%%%%%%%%%%%%%%%55
%%%%%%%%%%%%%%%%%%%%%%%%%%%%%%%%%%%%%%%%%%%%%%%%%%%%%%%%%%%%%%%%%%%%%%%%%%%%%%%%%%%%%%%%%%%%%%%%%%%%%%%%%%%%%%%55

 \bibliographystyle{amsplain}
\bibliography{PABM.bib} 
 
%%%%%%%%%%%%%%%%%%%%%%%%%%%%%%%%%%%%%%%%%%%%%%%%%%%%%%%%%%%%%%%%%%%%%%%%%%%%%%%%%%%%%%%%%%%%%%%%%%%%%%%%%%%%%%%55
%%%%%%%%%%%%%%%%%%%%%%%%%%%%%%%%%%%%%%%%%%%%%%%%%%%%%%%%%%%%%%%%%%%%%%%%%%%%%%%%%%%%%%%%%%%%%%%%%%%%%%%%%%%%%%%55

%\subsection{Adding Citations and a References List}
%Please use a \verb|.bib| file to store your references. When using Overleaf to prepare your manuscript, you can upload a \verb|.bib| file or import your Mendeley, CiteULike or Zotero library directly as a \verb|.bib| file\footnote{see \url{https://www.overleaf.com/blog/184}}. You can then cite entries from it, like this: \cite{lees2010theoretical}. Just remember to specify a bibliography style, as well as the filename of the \verb|.bib|.

%You can find a video tutorial here to learn more about BibTeX: \url{https://www.overleaf.com/help/97-how-to-include-a-bibliography-using-bibtex}.

%This template provides two options for the citation and reference list style: 
%\begin{description}
%\item[Numerical style] Use %\verb|\documentclass[...,num-refs]{wiley-article}|
%\item[Author-year style] Use %\verb|\documentclass[...,alpha-refs]{wiley-article}|
%\end{description}
% Contains proofs

%%%%%%%%%%%%%%%%%%%%%%%%%%%%%%%%%%%%%%%%%%%%%%%%%%%%%%%%%%%%%%%%%%%%%%%%%%%%%%%%%%%%%%%%%%%%%%%%%%%%%%%%%%%%%%%55
%%%%%%%%%%%%%%%%%%%%%%%%%%%%%%%%%%%%%%%%%%%%%%%%%%%%%%%%%%%%%%%%%%%%%%%%%%%%%%%%%%%%%%%%%%%%%%%%%%%%%%%%%%%%%%%55

\section{Proofs}
\label{sec:AppendA}
\setcounter{equation}{0}

\subsection {Proof  of Theorem~\ref{th:oracle}}
{\bf Overview. } The proof follows standard oracle inequality strategy.
We bound the error  $\|\hat{P} - P_{*}\|_F^2$  by the random error term 
$\Delta(\hat{Z},\hat{K}) = 2 \Tr [(A-P_{*})^T (\hat{P}  -P_{*})]$ 
plus the difference $\Pen(n,K_{*}) - \Pen(n,\hat{K})$ between the values of the penalty function at $K_*$ and $\hat{K}$.
Subsequently, we show that the random error term is bounded above by the sum of the  $\Pen(n,\hat{K})$
and a small multiple of  $\|\hat{P} - P_{*}\|_F^2$   with high probability. 
The latter leads to the conclusion that $\|\hat{P} - P_{*}\|_F^2$ is smaller than a multiple of $\Pen(n,K_{*})$   with high probability. 
The details of the proof is as follows.
\\
\\
{\bf Proof. }  Consider functions 
\beqn \label{eq:F1}
 F_1(n,K) & = &  C_1n K + C_2 K^2 \ln (ne) +C_3( \ln n + n \ln K) \\
\label{eq:F2}
F_2(n,K) & = &  2  \ln n + 2 n \ln K,
\eeqn
where $C_1, C_2$ and $C_3$  are  absolute constants.  
Define the penalty of the form
\be \label{eq:penalty}
\Pen (n,K) =  \lkr 2 +  16\, \beta_1^{-1} \rkr F_1(n,K)   +  \beta_2^{-1}  F_2(n,K),
\ee
where positive parameters $\beta_1$ and $\beta_2$ are such that $\beta_1 + \beta_2 <1$.
It is easy to see that by rearranging and combining the terms, the penalty in \eqref{eq:penalty} can be written in the form
\eqref{eq:pen1}, so we shall carry out the proof for the penalty given by \eqref{eq:penalty}.

Denote $\Xi = A- P_*$ and recall that,  given matrix $P_*$, entries  
$\Xi_{i,j} = A_{i,j}-(P_{*})_{ij}$ of $\Xi$  are the independent Bernoulli errors for $1 \leq i \leq j \leq n$ and
$ A_{i,j}= A_{j,i}$. Then, following notation \eqref{eq:permute}, for any $Z$ and $K$
\bes 
\Xi(Z,K) = \mathscr{P}_{Z,K}^T\Xi \mathscr{P}_{Z,K} \quad \mbox{and} \quad
 P_{*} (Z,K) = \mathscr{P}_{Z,K}^T P_{*} \mathscr{P}_{Z,K}.
\ees
Then it follows from \eqref{eq:opt_main} that 
\bes
\norm{\mathscr{P}_{\hat{Z},\hat{K}}^T A\mathscr{P}_{\hat{Z},\hat{K}}-\hat\Theta(\hat{Z},\hat{K}) }_F^2 + \Pen(n,\hat{K}) 
\leq  \norm{\mathscr{P}_{Z_{*},K_{*}}^T A \mathscr{P}_{Z_{*},K_{*}}  - \mathscr{P}_{Z_{*},K_{*}}^T P_{*}  \mathscr{P}_{Z_{*},K_{*}}}_F^2 + \Pen(n,K_{*}) 
\ees
Using the fact that permutation matrices are orthogonal, we can rewrite the previous inequality as 
\be \label{eq:main_ineq0}
\norm{ A -\mathscr{P}_{\hat{Z},\hat{K}}\hat\Theta(\hat{Z},\hat{K}) \mathscr{P}_{\hat{Z},\hat{K}}^T }_F^2 
 + \Pen(n,\hat{K})\leq \norm{A - P_{*}}_F^2  + \Pen(n,K_{*}).
\ee 
Hence,  \eqref{eq:main_ineq0} and  \eqref{eq:P_total_est}  yield 
\be  \label{eq:main_ineq1}
\norm{ A - \hat{P}}_F^2 \leq \norm{A - P_{*}}_F^2  + \Pen(n,K_{*}) - \Pen(n,\hat{K}) 
\ee
Subtracting and adding $P_*$ in the norm of the left-hand side of  \eqref{eq:main_ineq1},
we rewrite \eqref{eq:main_ineq1} as 
\be \label{eq:tot_err}
  \norm{\hat{P}  -P_{*}}_F^2 \leq  \Delta(\hat{Z},\hat{K}) + \Pen(n,K_{*}) - \Pen(n,\hat{K}),  
\ee
where
\be \label{eq:DelZK}
  \Delta(\hat{Z},\hat{K}) = 2 \Tr\left[\Xi^T (\hat{P}  -P_{*})\right].
\ee

Again, using orthogonality of the permutation matrices, we can rewrite 
\bes
\Delta(\hat{Z},\hat{K}) = 2 \langle\Xi(\hat{Z},\hat{K}),(\hat\Theta(\hat{Z},\hat{K}) - P_{*} (\hat{Z},\hat{K}))\rangle,
\ees
where $\langle A,B\rangle = \Tr(A^TB)$.
Then, in the block form, $\Delta(\hat{Z},\hat{K})$ appears as 
\be \label{eq:Del_blockform}
\Delta(\hat{Z},\hat{K}) =  \displaystyle \sum_{k,l=1}^{\hat{K}} \Delta^{(k,l)}(\hat{Z},\hat{K})
\ee
where 
\bes 
\Delta^{(k,l)}(\hat{Z},\hat{K}) = 2 \left\langle\Xi^{(k,l)}(\hat{Z},\hat{K}), \Pi_{\hat{u}, \hat{v}}\left( A^{(k,l)}(\hat{Z},\hat{K})\right) - 
P_{*}^{(k,l)} (\hat{Z},\hat{K}) \right\rangle 
\ees 
and $\Pi_{\hat{u}, \hat{v}}$ is defined in   \eqref{eq:Pi_uv} of Lemma~\ref{lem:lowrank_approx}.

Let $\tilde {u}=\tilde {u}^{(k,l)}(\hat{Z},\hat{K}), \tilde {v}={\tilde{v}^{(k,l)}}(\hat{Z},\hat{K})$ be the singular vectors 
of $P_{*}^{(k,l)}(\hat{Z},\hat{K})$ corresponding to the largest singular value  of 
$P_{*}^{(k,l)}(\hat{Z},\hat{K})$. % which are denoted by $\tilde{\sigma}_{\max}^{(k,l)}$.
Then, according to Lemma~\ref{lem:lowrank_approx}
% for  any $K$, any matrix $Z \in \calM_{n,K}$ and any matrix $P$ 
\be \label{not_r1approx_kl_block of P}
 \Pi_{\tilde{u}, \tilde{v}} \left(P_{*}^{(k,l)}(\hat{Z},\hat{K})\right) = 
\tilde {u}^{(k,l)}(\hat{Z},\hat{K})(\tilde{u}^{(k,l)}(\hat{Z},\hat{K}))^T P_{*}^{(k,l)}(\hat{Z},\hat{K}) 
\tilde{v}^{(k,l)}(\hat{Z},\hat{K}) (\tilde{v}^{(k,l)} (\hat{Z},\hat{K}))^T
\ee
Recall that 
\bes
 \Pi_{\hat{u}, \hat{v}} (A^{(k,l)}(\hat{Z},\hat{K})) = \Pi_{\hat{u},\hat{v}}  
\left[P_{*}^{(k,l)} (\hat{Z},\hat{K}) + \Xi^{(k,l)}(\hat{Z},\hat{K})\right], 
\ees
Then, $\Delta^{(k,l)}(\hat{Z},\hat{K})$ can be partitioned into the sums of three components
\be \label{eq:Delkl_sum}
\Delta^{(k,l)}(\hat{Z},\hat{K}) = \Delta_1^{(k,l)}(\hat{Z},\hat{K}) + \Delta_2^{(k,l)}(\hat{Z},\hat{K}) + 
\Delta_3^{(k,l)}(\hat{Z},\hat{K}), \quad 
k,l = 1,2,\cdots, K,
\ee
where 
\beqn 
\label{kl_blockdelta1_main_est_error}
\Delta_1^{(k,l)}(\hat{Z},\hat{K}) & = &  2 \langle \Xi^{(k,l)}(\hat{Z},\hat{K}),\Pi_{\hat{u}, \hat{v}} 
( \Xi^{(k,l)}(\hat{Z},\hat{K})) \rangle   \\
\label{kl_blockdelta2_main_est_error}
 \Delta_2^{(k,l)}(\hat{Z},\hat{K}) & = &  2 \langle \Xi^{(k,l)}(\hat{Z},\hat{K}), \Pi_{\tilde{u}, \tilde{v}} 
\left(P_{*}^{(k,l)}(\hat{Z},\hat{K})\right)- P_{*} ^{(k,l)}(\hat{Z},\hat{K})\rangle\\
\label{kl_blockdelta3_main_est_error}
 \Delta_3^{(k,l)}(\hat{Z},\hat{K}) & = & 2\langle \Xi^{(k,l)}(\hat{Z},\hat{K}),\Pi_{\hat{u}, \hat{v}} 
(P_{*}^{(k,l)}(\hat{Z},\hat{K}) ) -
\Pi_{\tilde{u}, \tilde{v}} \left(P_{*}^{(k,l)}(\hat{Z},\hat{K})\right) \rangle 
\eeqn
With some abuse of notations,   for any matrix $B$, let $\Pi_{\tilde{u}, \tilde{v}} \left(B(\hat{Z},\hat{K})\right)$ be  the matrix 
with blocks $\Pi_{\tilde{u}, \tilde{v}} \left(B^{(k,l)}(\hat{Z},\hat{K})\right)$, and 
$\Pi_{\hat{u}, \hat{v}}  \left(B(\hat{Z},\hat{K})\right)$ be  the matrix with blocks 
$\Pi_{\hat{u}, \hat{v}} \left(B^{(k,l)}(\hat{Z},\hat{K})\right)$,  $k,l = 1,2,\cdots, \hat{K}$. 
Then, it follows from \eqref{eq:Delkl_sum}--\eqref{kl_blockdelta3_main_est_error} that 
\be \label{eq:Del_sum}
\Delta(\hat{Z},\hat{K}) = \Delta_1(\hat{Z},\hat{K}) + \Delta_2(\hat{Z},\hat{K}) + \Delta_3(\hat{Z},\hat{K})  
\ee
where 
\beqn 
\label{delta1_main_est_error}
\Delta_1(\hat{Z},\hat{K}) & = & 2 \langle (\Xi(\hat{Z},\hat{K}),\Pi_{\hat{u}, \hat{v}} ( \Xi(\hat{Z},\hat{K})) \rangle \\
\label{delta2_main_est_error}
 \Delta_2(\hat{Z},\hat{K}) & = &  2 \langle \Xi(\hat{Z},\hat{K}), \Pi_{\tilde{u}, \tilde{v}}\left( P_{*}(\hat{Z},\hat{K})\right)- 
P_{*} (\hat{Z},\hat{K})\rangle\\
\label{delta3_main_est_error}
 \Delta_3(\hat{Z},\hat{K}) & = &  2\langle \Xi(\hat{Z},\hat{K}),\Pi_{\hat{u}, \hat{v}} (P_{*}(\hat{Z},\hat{K}) ) -
\Pi_{\tilde{u}, \tilde{v}}\left( P_{*}(\hat{Z},\hat{K})\right) \rangle  
\eeqn
Now, we need to derive an upper bound for each component in \eqref{eq:Delkl_sum} and \eqref{eq:Del_sum}.

%%%%%%%%%%%%%%%%%%%%%%%%%%%%%%%  Delta 1  %%%%%%%%%%%%%%%%%%%%%%%%%%%%%%%%%%%%%%%%%

Observe that 
\begin{align*}
\Delta_1^{(k,l)}(\hat{Z},\hat{K})&   
= 2 \langle \Xi^{(k,l)}(\hat{Z},\hat{K}),\Pi_{\hat{u}, \hat{v}} ( \Xi^{(k,l)}(\hat{Z},\hat{K})) \rangle  
= 2 \norm{\Pi_{\hat{u}, \hat{v}} ( \Xi^{(k,l)}(\hat{Z},\hat{K}))}_F^2\\
& \leq 2\norm{ \Xi^{(k,l)}(\hat{Z},\hat{K}) }_{op}^2.  
\end{align*}  
Now, fix $t$ and let $\Omega_{1}$ be the set where $\di\sum_{k,l = 1}^{\hat{K}}\ \norm{\Xi (\hat{Z},\hat{K})}_{op}^2 \leq   F_1(n,\hat{K}) + C_3 t$.
According to Lemma~\ref{lem:prob_error_bound_main}, 
\be \label{eq:P_Omega1}
\PP(\Omega_{1}) \geq  1-   \exp(-t),
\ee
and,  for $\omega \in \Omega_{1}$, one has 
\be \label{eq:Del1_bound}
\vert \Delta_1(\hat{Z},\hat{K})\vert \leq  
  2 \di\sum_{k,l = 1}^{\hat{K}} \norm{\Xi^{(k,l)}(\hat{Z},\hat{K})}_{op}^2 \leq 2 F_1(n,\hat{K}) + 2 C_3 t
\ee

%%%%%%%%%%%%%%%%%%%%%%%%%%%%%%%  Delta 2  %%%%%%%%%%%%%%%%%%%%%%%%%%%%%%%%%%%%%%%%%

Now, consider  $\Delta_2(\hat{Z},\hat{K})$ given by \eqref{delta2_main_est_error}.
Note that 
\be \label{eqn_useful}
| \Delta_2(\hat{Z},\hat{K}) | = 2\|\Pi_{\tilde{u}, \tilde{v}}\left( P_{*}(\hat{Z},\hat{K})\right)- 
P_{*} (\hat{Z},\hat{K}) \|_F |  \langle \Xi(\hat{Z},\hat{K}),H_{\tilde{u}, \tilde{v}}(\hat{Z},\hat{K}) \rangle|
\ee  
where 
\bes
H_{\tilde{u},\tilde{v}}(\hat{Z},\hat{K})= \frac{\Pi_{\tilde{u}, \tilde{v}}\left( P_{*}(\hat{Z},\hat{K})\right)- 
P_{*} (\hat{Z},\hat{K})}{\|\Pi_{\tilde{u}\tilde{v}}\left( P_{*}(\hat{Z},\hat{K})\right)- P_{*} (\hat{Z},\hat{K})\|_F }
\ees
Since  for any $a,b$ and $\alpha_1 > 0$, one has $2ab \leq \alpha_1 a^2 +  b^2/\alpha_1$, obtain
\be \label{eq:Del2_sum} 
|\Delta_2(\hat{Z},\hat{K})| \leq  \alpha_1 \|\Pi_{\tilde{u}, \tilde{v}}\left( P_{*}(\hat{Z},\hat{K})\right)- 
P_{*} (\hat{Z},\hat{K}) \|_F^2 + 1/\alpha_1\, |  \langle \, \Xi(\hat{Z},\hat{K}),H_{\tilde{u},\tilde{v}}(\hat{Z},\hat{K})\, \rangle|^2
\ee 
Observe that if $K$ and $Z \in \calM_{n,K}$ are fixed, then $H_{\tilde{u},\tilde{v}}(Z,K)$ 
is  fixed and, for any $K$ and $Z$, one has  $\| H_{\tilde{u},\tilde{v}}(Z,K)\|_F =1$.
Note also that, for fixed $K$ and $Z$, permuted matrix  $\Xi(Z,K) \in [0,1]^{n\times n}$ 
contains independent Bernoulli errors. It is well known that if $\xi$ is a vector of independent Bernoulli errors
and $h$ is any fixed vector with $\di\sum_{i = 1}^{n} h_i^2 = 1$, then, for any $x>0$,   Hoeffding's inequality yields
$$ 
\PP(|\xi^T h|^2> x) \leq 2 \exp(- x/2) 
$$
Since
$\langle \Xi(Z,K),H_{\tilde{u},\tilde{v}}(Z,K)\rangle =  [\vect(\Xi(Z,K))]^T  \vect(H_{\tilde{u},\tilde{v}}(Z,K))$,
obtain  for any fixed $K$ and $Z$:
\bes 
\PP  \left( | \langle \Xi(Z,K),H_{\tilde{u},\tilde{v}}(Z,K)\rangle|^2 - x > 0 \right) \leq 2 \exp(-x/2) 
\ees
Now, applying the union bound, derive
\begin{align} 
& \PP  \left( | \langle \Xi(\hat{Z},\hat{K}),H_{\tilde{u},\tilde{v}}(\hat{Z},\hat{K}) \rangle|^2 - F_2(n,\hat{K}) > 2t \right) \nonumber\\ 
\leq  
& \PP  \left( \underset{1 \leq K \leq n}{\max} \di\underset  {Z \in \mathcal{M}_{n,k}} {\max} 
\left( |\langle \Xi(Z,K),H_{\tilde{u},\tilde{v}}(Z,K) \rangle|^2 - F_2(n,K) \right)  > 2t \right) \label{eq:Del2_union} \\
\leq  
& 2n K^{n} \exp \lfi -F_2(n,K)/2 - t \rfi =  2 \exp(-t), \nonumber
\end{align}
where $F_2(n,K)$ is defined in \eqref{eq:F2}.  
By Lemma~\ref{lem:Pi_orth}, one has 
\bes
\|\Pi_{\tilde{u}, \tilde{v}}\left( P_{*}(\hat{Z},\hat{K})\right)- P_{*} (\hat{Z},\hat{K}) \|_F^2 \leq 
\|\Pi_{\hat{u}, \hat{v}}\left( P_{*}(\hat{Z},\hat{K})\right)- P_{*} (\hat{Z},\hat{K}) \|_F^2 \leq 
\| \hat{P}  - P_{*}  \|_F^2.
\ees
Denote the set on which \eqref{eq:Del2_union} holds by $\Omega_{2}^C$, so that 
\be \label{eq:P_Omega2}
\PP(\Omega_{2}) \geq  1 - 2 \exp(-t).
\ee
Then inequalities \eqref{eq:Del2_sum}  and \eqref{eq:Del2_union} imply that, for any $\alpha_1>0$, $ t >0$ and 
any $\omega \in \Omega_{2}$,  one has
\be \label{eq:Del2_bound}
 |\Delta_2(\hat{Z},\hat{K})|  \leq  \alpha_1  \| \hat{P}  - P_{*}  \|_F^2 + 
1/\alpha_1 \, F_2(n,\hat{K}) + 2\,t/\alpha_1.
\ee

%%%%%%%%%%%%%%%%%%%%%%%%%%%%%%%  Delta 3  %%%%%%%%%%%%%%%%%%%%%%%%%%%%%%%%%%%%%%%%%

Now consider $\Delta_3(\hat{Z},\hat{K})$ defined in \eqref{delta3_main_est_error}
with components \eqref{kl_blockdelta3_main_est_error}.
Note that matrices $\Pi_{\hat{u}, \hat{v}} (P_{*}^{(k,l)}(\hat{Z},\hat{K}) ) -
\Pi_{\tilde{u}, \tilde{v}} \left(P_{*}^{(k,l)}(\hat{Z},\hat{K})\right)$ have rank at most two.
Use the fact that (see, e.g., Giraud  (2014), page 123)
\be \label{eq:Ky-Fan}
\langle A, B \rangle  \leq \|A\|_{(2,r)} \|B\|_{(2,r)} \leq 2 \|A\|_{op} \|B\|_F,\quad r = \min \{\text{rank}(A), \text{rank}(B)\}.
\ee
Here $\|A\|_{(2,q)}$ is the Ky-Fan $(2,q)$ norm 
$$
\|A\|^2_{(2,q)} = \sum_{j=1}^q \sig_j^2(A) \leq \|A\|^2_F,
$$
where $\sig_j(A)$ are the singular values of $A$.
Applying inequality \eqref{eq:Ky-Fan} with $r=2$ and taking into account that for any matrix $A$ one has 
$\|A\|^2_{(2,2)} \leq 2 \|A\|_{op}^2$, derive
\bes
|\Delta_3^{(k,l)}(\hat{Z},\hat{K})| \leq  4 \| \Xi^{(k,l)}(\hat{Z},\hat{K})\|_{op} \|\Pi_{\hat{u}, \hat{v}} (P_{*}^{(k,l)}(\hat{Z},\hat{K}) ) 
- \Pi_{\tilde{u}, \tilde{v}} \left(P_{*}^{(k,l)}(\hat{Z},\hat{K})\right) \|_F.
\ees
Then, for any $\alpha_2>0$, obtain 
\begin{align}
& |\Delta_3(\hat{Z},\hat{K}) | \leq  \di\sum_{k,l = 1}^{\hat{K}} |\Delta_3^{(k,l)}(\hat{Z},\hat{K})| \label{eq:Del3_sum} \\
& \leq \frac{2}{\alpha_2} \di\sum_{k,l = 1}^{\hat{K}} \| \Xi^{(k,l)}(\hat{Z},\hat{K})\|_{op}^2 + 
2\alpha_2 \di\sum_{k,l = 1}^{\hat{K}}  \|  \Pi_{\hat{u}, \hat{v}} (P_{*}^{(k,l)}(\hat{Z},\hat{K}) )  
- \Pi_{\tilde{u}, \tilde{v}} \left(P_{*}^{(k,l)}(\hat{Z},\hat{K})\right) \|_F^2. \nonumber
\end{align} 
Note that, by Lemma~\ref{lem:Pi_orth}, 
\begin{align*}
& \|\Pi_{\hat{u}, \hat{v}} (P_{*}^{(k,l)}(\hat{Z},\hat{K})) -\Pi_{\tilde{u}\tilde{v}} \left(P_{*}^{(k,l)}(\hat{Z},\hat{K})\right) \|_F^2\\
\leq 
& 2 \|\Pi_{\hat{u}, \hat{v}} (P_{*}^{(k,l)}(\hat{Z},\hat{K})) - P_{*}^{(k,l)}(\hat{Z},\hat{K})\|^2_F + 
2 \| \Pi_{\tilde{u}, \tilde{v}} (P_{*}^{(k,l)}(\hat{Z},\hat{K})) - P_{*}^{(k,l)}(\hat{Z},\hat{K})\|^2_F\\
\leq 
& 4 \|\Pi_{\hat{u}, \hat{v}} (P_{*}^{(k,l)}(\hat{Z},\hat{K})) - P_{*}^{(k,l)}(\hat{Z},\hat{K})\|^2_F \\
\leq 
& 4 \|\Pi_{\hat{u}, \hat{v}}(A^{(k,l)}(\hat{Z},\hat{K})) - P_{*}^{(k,l)}(\hat{Z},\hat{K}) \|_F^2  
= 4 \| \hat{\Theta}^{(k,l)}(\hat{Z},\hat{K}) - P_{*}^{(k,l)}(\hat{Z},\hat{K}) \|_F^2  
\end{align*}
Therefore,
\begin{align}  
& \sum_{k,l = 1}^{\hat{K}}  \|\Pi_{\hat{u}, \hat{v}} (P_{*}^{(k,l)}(\hat{Z},\hat{K})) -
\Pi_{\tilde{u}\tilde{v}} \left(P_{*}^{(k,l)}(\hat{Z},\hat{K})\right) \|_F^2 \leq  \nonumber\\
&  4  \norm{\hat{\Theta}(\hat{Z},\hat{K}) -P_{*}(\hat{Z},\hat{K})}_F^2 
= 4 \| \hat{P}  - P_{*}  \|_F^2  \label{eq:Del3_part}
\end{align}
Combine  inequalities  \eqref{eq:Del3_sum} and \eqref{eq:Del3_part}  and recall  that 
$\norm{\Xi (\hat{Z},\hat{K})}_{op}^2 \leq  F_1(n,\hat{K}) + C_3 \, t$ for 
$\omega \in \Omega_{1}$. Then,  for any $\alpha_2>0$  and  $\omega \in \Omega_{1}$, 
one has
\be \label{eq:Del3_bound}
  |\Delta_3(\hat{Z},\hat{K})| \leq 8 \alpha_2 \| \hat{P}  - P_{*}  \|_F^2 + 
2/\alpha_2 F_1(n,\hat{K}) + 2 C_3\, t/\alpha_2.
\ee

Now, let $\Omega =  \Omega_{1} \cap \Omega_{2}$. Then, \eqref{eq:P_Omega1} and \eqref{eq:P_Omega2}
imply that $\PP(\Omega) \geq 1 - 3 \exp(-t)$ and, for $\om \in \Omega$,
inequalities \eqref{eq:Del1_bound}, \eqref{eq:Del2_bound} and \eqref{eq:Del3_bound} simultaneously   hold.
Hence, by \eqref{eq:Del_sum}, derive that, for any $\om \in \Omega$,
\bes 
|\Delta(\hat{Z},\hat{K})| \leq (2 + 2/\alpha_2) F_1(n,\hat{K}))  + 1/\alpha_1\, F_2(n,\hat{K})  
+ (\alpha_1 + 8 \alpha_2) \| \hat{P}  - P_{*}  \|_F^2 + 2(C_3 + 1/\alpha_1 + C_3/\alpha_2)\, t.
\ees 
Combination of the last inequality and \eqref{eq:tot_err} yields that, for $\alpha_1 +  8\alpha_2 <1$ and any $\om \in \Omega$,
%\bes 
\begin{align*}
 (1-\alpha_1 - 8\alpha_2)\, \norm{\hat{P}  -P_{*}}_F^2 \leq  
\left( 2 + \frac{2}{\alpha_2}\right) F_1(n,\hat{K})  + \frac{1}{\alpha_1} F_2(n,\hat{K}) + \Pen(n,K_{*}) - \Pen(n,\hat{K})\\
 + 2(C_3 + 1/\alpha_1 + C_3/\alpha_2)\, t   
\end{align*}
%\ees
% provided $\alpha_1 +  8\alpha_2 <1$. 
Setting $\Pen(n,K) =  (2 + 2/\alpha_2) F_1(n,K)   + 1/\alpha_1  F_2(n,K)$ and dividing 
by $(1-\alpha_1 - 8\alpha_2)$, obtain that 
\be \label{eq:tot_err_new}
\PP \lfi \|{\hat{P}  -P_{*}}\|_F^2 \leq  (1-\alpha_1 - 8\alpha_2)^{-1} \, \Pen(n,K_{*}) + \tilde{C}\, t \rfi 
\geq 1 - 3 e^{-t}
\ee
where
\be \label{eq:tildeC}
\tilde{C} =  2\, (1-\alpha_1 - 8\alpha_2)^{-1} \,(C_3 + 1/\alpha_1 + C_3/\alpha_2) 
\ee
In order to derive  \eqref{eq:oracle}, set $\beta_1 = 8 \alpha_2$ and $\beta_2 =  \alpha_1$.
%
%
%%%%%%%%%%%%%%%%%%%%%%%%%%%%%%% Expectation %%%%%%%%%%%%%%%%%%%%%%%%%%%%%%%%%%%%%%%%%
%
%
In order to obtain  the upper bound  \eqref{eq:oracle Expectation}  note that 
for $\xi = \|{\hat{P} -P_{*}}\|_F^2 -  (1 - \beta_1 -\beta_2)^{-1} \,  \Pen (n,K_{*})$, 
one has 
$%$
\EE \|{\hat{P} - P_{*}}\|_F^2  = (1 - \beta_1 -\beta_2)^{-1} \,  \Pen (n,K_{*})  + \EE \xi, 
$%$
where
\bes
\EE \xi \leq \int_0^{\infty} \PP(\xi > z) dz = \tilde{C}\int_0^{\infty} \PP(\xi >  \tilde{C}t) dt  
\leq \tilde{C}\int_0^{\infty} 3 \, e^{-t}\, dt  = 3\tilde{C},
\ees
which yields \eqref{eq:oracle Expectation}.
\\

%%%%%%%%%%%%%%%%%%%%%%%%%%%%%%%%%%%%%%%%%%%%%%%%%%%%%%%%%%%%%%%%%%%%%%%%%%%%%%%%%%%%%%%%%%%%%%%%%%%%%%%%%%%%%%%%%%%%%%%%%%%%%%
%%%%%%%%%%%%%%%%%%%%%%%%%%%%%%%%%%%%%%%%%%%%%%%%%%%%%%%%%%%%%%%%%%%%%%%%%%%%%%%%%%%%%%%%%%%%%%%%%%%%%%%%%%%%%%%%%%%%%%%%%%%%%%

\subsection {Proof  of Theorem~\ref{th:oracle_sparse}}
The idea of the proof is to essentially repeat the steps in the proof of Theorem~\ref{th:oracle} 
without the penalty and the union bound over all possible values of $K$. The main difference here is that 
inequality \eqref{eq:lem3} of Lemma~\ref{lem:prob_bound_error} is replaced everywhere by the inequality 
\eqref{eq:sparse_prob_bou_er}  of Lemma~\ref{lem:sparse_prob_bound_error}.
\\

\noindent
Consider function  
\be  \label{eq:F3S}
 F_3(n) =  C_1 \tau_n K + C_2 K^2  +C_3(  n \ln K) 
% \label{eq:F2S}
% F_2(n) & = &  2 n \ln K,
\ee 
where $C_1, C_2$ and $C_3$  are  absolute constants.  Similar to the proof of Theorem~\ref{th:oracle}, denote
$\Xi(Z) = \mathscr{P}_{Z}^T\Xi \mathscr{P}_{Z}$ and  
$P_{*} (Z) = \mathscr{P}_{Z}^T P_{*} \mathscr{P}_{Z}$.
Obtain that
\be \label{eq:tot_errS}
  \|\hat{P}  -P_{*}\|_F^2 \leq  \Delta(\hat{Z}) =    2 \Tr\left[\Xi^T (\hat{P}  -P_{*})\right].
\ee
Similarly to the proof of  Theorem~\ref{th:oracle}, partition $\Delta(\hat{Z})$ as
\be   \label{eq:Del_sumS}
 \Delta(\hat{Z}) =  \displaystyle \sum_{k,l=1}^K \Delta^{(k,l)}(\hat{Z}) = 
 \Delta_1(\hat{Z}) + \Delta_2(\hat{Z}) + \Delta_3(\hat{Z})  
\ee
 where 
\beqn 
\label{delta1_main_est_errorS}
\Delta_1(\hat{Z}) & = & 2 \langle (\Xi(\hat{Z}),\Pi_{\hat{u}, \hat{v}} ( \Xi(\hat{Z})) \rangle \\
\label{delta2_main_est_errorS}
 \Delta_2(\hat{Z}) & = &  2 \langle \Xi(\hat{Z}), \Pi_{\tilde{u}, \tilde{v}}\left( P_{*}(\hat{Z})\right)- 
P_{*} (\hat{Z})\rangle\\
\label{delta3_main_est_errorS}
 \Delta_3(\hat{Z}) & = &  2\langle \Xi(\hat{Z}),\Pi_{\hat{u}, \hat{v}} (P_{*}(\hat{Z}) ) -
\Pi_{\tilde{u}, \tilde{v}}\left( P_{*}(\hat{Z})\right) \rangle  
\eeqn

%%%%%%%%%%%%%%%%%%%%%%%%%%%%%%%  Delta 1  %%%%%%%%%%%%%%%%%%%%%%%%%%%%%%%%%%%%%%%%%

\noindent
Observe that,  by Lemma~\ref{lem:prob_error_bound_mainS}, for $\omega \in \Omega_{3}$, 
where $\PP(\Omega_{3}) \geq  1-   \exp(-t)$, one has
\be \label{eq:Del1_boundS}
\vert\Delta_1(\hat{Z})\vert \leq  
  2 \di\sum_{k,l = 1}^K \norm{\Xi^{(k,l)}(\hat{Z})}_{op}^2 \leq 2 F_3(n) + 2 C_3 t
\ee

%%%%%%%%%%%%%%%%%%%%%%%%%%%%%%%  Delta 2  %%%%%%%%%%%%%%%%%%%%%%%%%%%%%%%%%%%%%%%%%

\noindent
For $\Delta_2(\hat{Z})$, given by \eqref{delta2_main_est_errorS}, for any 
 $\alpha_1 > 0$, one has  
\be \label{eq:Del2_sumS} 
|\Delta_2(\hat{Z})| \leq  \alpha_1 \|\Pi_{\tilde{u}, \tilde{v}}\left( P_{*}(\hat{Z})\right)- 
P_{*} (\hat{Z}) \|_F^2 + 1/\alpha_1\, |  \langle \, \Xi(\hat{Z}),H_{\tilde{u},\tilde{v}}(\hat{Z})\, \rangle|^2
\ee 
where  
\bes
H_{\tilde{u},\tilde{v}}(\hat{Z})= \frac{\Pi_{\tilde{u}, \tilde{v}}\left( P_{*}(\hat{Z})\right)- 
P_{*} (\hat{Z})}{\|\Pi_{\tilde{u}\tilde{v}}\left( P_{*}(\hat{Z})\right)- P_{*} (\hat{Z})\|_F }
\ees
Similarly to the proof of Theorem~\ref{th:oracle}, derive that there exists a set $\Om_4$ 
with $\PP(\Omega_{4}) \geq  1 - 2 \exp(-t)$ such that for $\om \in \Om_4$,  
\bes 
 |\langle \Xi(\hat{Z}),H_{\tilde{u},\tilde{v}}(\hat{Z}) \rangle|^2 \leq  2 n \ln K  + 2t.
\ees
Since 
\bes
\|\Pi_{\tilde{u}, \tilde{v}}\left( P_{*}(\hat{Z})\right)- P_{*} (\hat{Z}) \|_F^2 \leq 
\|\Pi_{\hat{u}, \hat{v}}\left( P_{*}(\hat{Z})\right)- P_{*} (\hat{Z}) \|_F^2 \leq 
\| \hat{P}  - P_{*}  \|_F^2, 
\ees
the above inequalities imply that, for any $\alpha_1>0$, $ t >0$ and 
any $\omega \in \Omega_{4}$,  one has
\be \label{eq:Del2_boundS}
 |\Delta_2(\hat{Z})|  \leq  \alpha_1  \| \hat{P}  - P_{*}  \|_F^2 + 
2/\alpha_1 \, n \ln K + 2\,t/\alpha_1.
\ee

%%%%%%%%%%%%%%%%%%%%%%%%%%%%%%%  Delta 3  %%%%%%%%%%%%%%%%%%%%%%%%%%%%%%%%%%%%%%%%%

\noindent
Now consider $\Delta_3(\hat{Z})$ defined in \eqref{delta3_main_est_errorS}.
Similarly to the proof of Theorem~\ref{th:oracle}, for any $\alpha_2>0$, obtain 
\be   \label{eq:Del3_SSS} 
 |\Delta_3(\hat{Z}) | \leq   \frac{2}{\alpha_2} \di\sum_{k,l = 1}^K \| \Xi^{(k,l)}(\hat{Z})\|_{op}^2 + 
2\alpha_2 \di\sum_{k,l = 1}^K  \|  \Pi_{\hat{u}, \hat{v}} (P_{*}^{(k,l)}(\hat{Z}) )  
- \Pi_{\tilde{u}, \tilde{v}} \left(P_{*}^{(k,l)}(\hat{Z})\right) \|_F^2,  
\ee  
where 
\bes 
  \sum_{k,l = 1}^K  \|\Pi_{\hat{u}, \hat{v}} (P_{*}^{(k,l)}(\hat{Z})) -
\Pi_{\tilde{u}\tilde{v}} \left(P_{*}^{(k,l)}(\hat{Z})\right) \|_F^2 \leq  
 4 \| \hat{P}  - P_{*}  \|_F^2  
\ees
Since  for $\omega \in \Omega_{3}$, the first term in \eqref{eq:Del3_SSS} 
is bounded above by $2/\alpha_2(F_3(n) +   C_3\, t)$, 
for any $\alpha_2>0$  and  $\omega \in \Omega_{3}$, we obtain
\be \label{eq:Del3_boundS}
  |\Delta_3(\hat{Z})| \leq 8 \alpha_2 \| \hat{P}  - P_{*}  \|_F^2 + 
2/\alpha_2 F_3(n) + 2 C_3\, t/\alpha_2.
\ee

%%%%%%%%%%%%%%%%%%%%%%%%%%%%%%%%%%% All %%%%%%%%%%%%%%%%%%%%%%%%%%%%%%%%%%%%%%%%%%%

Finally, combination of the terms in \eqref{eq:Del_sumS},   yields that, 
for $\alpha_1 +  8\alpha_2 <1$ and any $\om \in \Omega_3 \cap \Om_4$,
%\bes 
\begin{align*}
 (1-\alpha_1 - 8\alpha_2)\, \norm{\hat{P}  -P_{*}}_F^2 \leq  
\left( 2 + \frac{2}{\alpha_2}\right) F_1(n)  + \frac{2}{\alpha_1} n \ln K + 2(C_3 + 1/\alpha_1 + C_3/\alpha_2)\, t   
\end{align*}
To complete the proof of the \eqref{eq:oracle_sparse}, divide the last expression by $ (1-\alpha_1 - 8\alpha_2)$ and simplify,
arriving at 
\be \label{eq:tot_err_newS}
\PP \lfi   \|\hat{P}  -P_{*}\|_F^2    \leq    H_1 \tau_n n\, K + H_2  K^2   + H_3  n\, \ln K   
+   \tilde{H}\,t  \rfi \geq  1 - 3 e^{-t}, 
\ee 
where
\begin{align*}  
& H_1=  \frac{2C_1}{\alpha_2} (1 +\alpha_2)(1-\alpha_1 - 8\alpha_2)^{-1} \\
& H_2 = \frac{2C_2}{\alpha_2} (1 +\alpha_2)(1-\alpha_1 - 8\alpha_2)^{-1} \\
& H_3 =  2\left(\frac{1}{\alpha_1} + \frac{ C_3}{\alpha_2 }(1 +\alpha_2) \right) (1-\alpha_1 - 8\alpha_2)^{-1} \\
& \tilde{H} = 2( C_3+1 /\alpha_1 + C_3/\alpha_2)(1-\alpha_1 - 8\alpha_2)^{-1} 
\end{align*}
Finally, the proof of inequality \eqref{eq:oracle Expectation_sparse} is identical to the proof
of \eqref{eq:oracle Expectation}.

%%%%%%%%%%%%%%%%%%%%%%%%%%%%%%%%%%%%%%%%%%%%%%%%%%%%%%%%%%%%%%%%%%%%%%%%%%%%%%%%%%%%%%%%%%%%%%%%%%%%%%%%%%%%%%%%%%%%%%%%%%%%%%
%%%%%%%%%%%%%%%%%%%%%%%%%%%%%%%%%%%%%%%%%%%%%%%%%%%%%%%%%%%%%%%%%%%%%%%%%%%%%%%%%%%%%%%%%%%%%%%%%%%%%%%%%%%%%%%%%%%%%%%%%%%%%%

\subsection {Proof  of Lemma~\ref{lem:detect} }
Note that the left hand side of inequality \eqref{eq:detect} is equal to identical zero. 
Consider matrix  $Z \in \calM_{n, K_*}$ such that $Z$ cannot be obtained from $Z_*$ by
a permutation of columns. Let $i$ be a misclassified node, so that it belongs to communities $l_*$ and $l$
according to   $Z_*$ and $Z$, respectively. Let $t$ be such that vectors ${\Lambda}^{(t,1)}, \ldots,  {\Lambda}^{(t,K_*)}$ 
are linearly independent. Then, sub-matrix $P_*^{(t,k)}(Z)$ of matrix  $P_*(Z)$ will contain multiples of columns ${\Lambda}^{(t,l_*)}$ and 
${\Lambda}^{(t,l)}$. Under Assumption~{\bf A1}, those multiples cannot be identically equal to zero. Under Assumption~{\bf A1*}, one 
can choose $t$ such that they are not identically equal to zero, since matrix $P$ does not have zero columns. Then, the rank of  matrix 
$P_*^{(t,k)}(Z)$ is   at least two, so that, the right hand side of \eqref{eq:detect} is positive, which completes the proof.

%%%%%%%%%%%%%%%%%%%%%%%%%%%%%%%%%%%%%%%%%%%%%%%%%%%%%%%%%%%%%%%%%%%%%%%%%%%%%%%%%%%%%%%%%%%%%%%%%%%%%%%%%%%%%%%%%%%%%%%%%%%%%%
%%%%%%%%%%%%%%%%%%%%%%%%%%%%%%%%%%%%%%%%%%%%%%%%%%%%%%%%%%%%%%%%%%%%%%%%%%%%%%%%%%%%%%%%%%%%%%%%%%%%%%%%%%%%%%%%%%%%%%%%%%%%%%

\subsection {Proof  of Theorem~\ref{th:clust} }
Note that it follows from \eqref{eq:opt_ZK3} that  
\be \label{eq:clust1}  
\sum_{k,l = 1}^K \norm {A^{(k,l)}(\hat{Z}) -  \Pi_{\hat{u}, \hat{v}} \left(A^{(k,l)}(\hat{Z})\right)}_{F}^2  
\leq    \sum_{k,l = 1}^K \norm {A^{(k,l)}(Z_{*}) -  \Pi_{\hat{u}, \hat{v}} \left(A^{(k,l)}(Z_{*})\right)}_{F}^2  
\ee 
Observe  that for any $Z \in \mathcal{M}_{n,K}$, one has
\begin{align*}
\sum_{k,l = 1}^K \norm {A^{(k,l)}(Z) -  \Pi_{\hat{u}, \hat{v}} \left(A^{(k,l)}(Z)\right) }_{F}^2   
= \sum_{k,l = 1}^K   \left\{ \norm {A^{(k,l)}(Z)} _{F}^2  - \norm {\Pi_{\hat{u}, \hat{v}} \left(A^{(k,l)}(Z)\right)}_{F}^2 \right\}, 
\end{align*}
so that, due to $\displaystyle \sum_{k,l = 1}^K \norm {A^{(k,l)}(Z) }_{F}^2  = \norm{A}_{F}^2$, 
\eqref{eq:clust1} can be re-written as 
\be \label{eq:clust2}  
\sum_{k,l = 1}^K \norm {\Pi_{\hat{u}, \hat{v}} \left(A^{(k,l)}(\hat{Z})\right)}_{F}^2  
\geq    \sum_{k,l = 1}^K \norm {\Pi_{\hat{u}, \hat{v}} \left(A^{(k,l)}(Z_{*})\right)}_{F}^2  
\ee 
Applying  Proposition 6.2 of  Giraud \cite{Giraud:1999833}, obtain  
\bes
\|{ \Pi_{\hat{u}, \hat{v}}  (A^{(k,l)}(Z)) - P_{*}^{(k,l)} (Z) }\|_{F}^2 
\leq \frac{(2+\theta)^2}{\theta^2} \ \sum_{r =2}^{\min\{n_k,n_l\}}  \sigma_r^2 P_{*}^{(k,l)} (Z) + 
\frac{2(1+\theta)(2+\theta)}{\theta} \, \|{\Xi^{(k,l)}(Z)}\|_{op}^2,
\ees
where $\theta >0$ is an arbitrary constant, $P_{*}$ is the true matrix of probabilities, 
$\Xi^{(k,l)}(Z) = A^{(k,l)}(Z) - P_{*}^{(k,l)} (Z)$, and $\sigma_r(B)$ is the $r-$th largest singular value of $B$.
Since  matrix  $P_{*}^{(k,l)} (Z_{*})$ has rank one, the previous inequality yields for $\theta = \sqrt{ 2}$
\be\label{eq:opt_Z5}
 \|{ \Pi_{\hat{u}, \hat{v}} \left(A^{(k,l)}(Z_{*})\right) - P_{*}^{(k,l)} (Z_{*}) }\|_{F}^2 \leq 2(1+\sqrt{ 2})^2 \|{\Xi^{(k,l)}(Z_{*})}\|_{op}^2
 \ee

Using  Lemma~\ref{lem:prob_bound_error}, derive for any $t>0$ that
\be \label{eq:ineq1}
\PP \left\{\sum_{k,l =1}^K \norm{\Xi^{(k,l)}(Z_{*})}_{op}^2   \leq  C_1 nK + C_2 K^2 \ln (ne) + C_3\, t \right\}   \geq  1 - \exp (-t).
\ee 
Also, since   $\card(\mathcal{M}_{n,K}) = K^n$ , replacing $t$ by  $n \ln K + t $ and applying union bound, obtain
\be  \label{eq:ineq2}
\PP \left\{ \sum_{k,l =1}^K \|{\Xi^{(k,l)}(\hat{Z})}\|_{op}^2  
\leq C_1nK + C_2 K^2 \ln (ne) + C_3( n \ln K + t ) \right\}   \geq   1 - \exp (-t).
\ee

Note that for any $\alpha_1 \in (0,1)$,
\begin{align*}  
& \norm {\Pi_{\hat{u}, \hat{v}} \left(A^{(k,l)}(Z_{*})\right)}_{F}^2  =  
\norm { \Pi_{\hat{u}, \hat{v}} \left(A^{(k,l)}(Z_{*})\right) - P_{*}^{(k,l)} (Z_{*})  + P_{*}^{(k,l)} (Z_{*})}_{F}^2 \geq \\
% 
% &    \norm { \Pi_{\hat{u}, \hat{v}} \left(A^{(k,l)}(Z_{*})\right) - P_{*}^{(k,l)} (Z_{*})}_{F}^2  + \norm{P_{*}^{(k,l)} (Z_{*})}_{F}^2
%   + 2 \langle P_{*}^{(k,l)} (Z_{*}), \Pi_{\hat{u}, \hat{v}} \left(A^{(k,l)}(Z_{*})\right) - P_{*}^{(k,l)} (Z_{*})\rangle \\
%
&  (1-\alpha_1)  \norm{P_{*}^{(k,l)} (Z_{*})}_{F}^2 - \left( \alpha_1^{-1} - 1\right) 
\norm{\Pi_{\hat{u}, \hat{v}} \left(A^{(k,l)}(Z_{*})\right) - P_{*}^{(k,l)} (Z_{*})}_{F}^2 
\end{align*}
Combining the last inequality with \eqref{eq:opt_Z5} and taking a sum, obtain
 \be\label{eq:opt_Z6}
\sum_{k,l =1}^K\,  \norm { \Pi_{\hat{u}, \hat{v}} \left(A^{(k,l)}(Z_{*})\right)}_{F}^2  \geq (1-\alpha_1)  \norm{P_{*}}_{F}^2 - 
 2(1+\sqrt{ 2})^2 \left(\frac{1}{\alpha_1} - 1\right)  \sum_{k,l =1}^K\,  \norm{\Xi^{(k,l)}(Z_{*})}_{op}^2,
 \ee
where we used the fact that $\norm{P_{*}  (Z_{*})}_{F}  = \norm{P_{*}}_{F}$.
 On the other hand, for any $ Z \in \mathcal{M}_{n,K} $ and any $\alpha_2 > 0$,
 \bes
\norm{ \Pi_{\hat{u}, \hat{v}}\left( A^{(k,l)}(Z) \right)}_{F}^2 
\leq (1 + \alpha_2)\norm{ \Pi_{\hat{u}, \hat{v}}\left( P_{*}^{(k,l)} (Z) \right)}_{F}^2 + 
\left( 1 +  \alpha_2^{-1}\right) \norm{\Pi_{\hat{u}, \hat{v}}\left( \Xi^{(k,l)}(Z) \right)}_{F}^2, 
\ees
so that 
\be\label{eq:opt_Z7}
\norm{\Pi_{\hat{u}, \hat{v}}\left( A^{(k,l)}(Z) \right)}_{F}^2 \leq (1 + \alpha_2)\norm{P_{*}^{(k,l)} (Z) }_{op}^2 + 
\left( 1 + \alpha_2^{-1} \right) \norm{ \Xi^{(k,l)}(Z) }_{op}^2 
\ee
    
Now, we prove the theorem by contradiction. Assume that $\hat{Z} \in \Upsilon(Z_{*},\rho_n)$ is the solution of optimization problem 
\eqref{eq:opt_ZK3}. Then, inequality   \eqref{eq:clust2}  holds. Combining  \eqref{eq:clust2}, \eqref{eq:opt_Z6} and \eqref{eq:opt_Z7},
obtain that 
\begin{align*}
& (1-\alpha_1)  \norm{P_{*}}_{F}^2 - 
 2(1+\sqrt{ 2})^2 \left(\frac{1}{\alpha_1} - 1\right)  \sum_{k,l =1}^K\,  \norm{\Xi^{(k,l)}(Z_{*})}_{op}^2 \leq \\
& (1 + \alpha_2) \sum_{k,l =1}^K\, \norm{P_{*}^{(k,l)} (\hat{Z}) }_{op}^2 + 
\left( 1 + \alpha_2^{-1} \right) \sum_{k,l =1}^K\, \norm{ \Xi^{(k,l)}(\hat{Z}) }_{op}^2 
\end{align*}
Due to  \eqref{eq:ineq1} and \eqref{eq:ineq2}, with probability at least $1 - 2 \exp(-t)$, the last inequality yields 
\bes   
(1-\alpha_1) \,  \|{P_{*}}\|_{F}^2 - (1 + \alpha_2) \sum_{k,l =1}^K\, \|{P_{*}^{(k,l)} (\hat{Z}) }\|_{op}^2 \leq 
(1-\alpha_1) \,  C (\alpha_1,\alpha_2)\,  [ C_1  n K + C_2 K^2 \ln (ne)  + C_3(n \ln K +t)].
\ees
where  
$C (\alpha_1,\alpha_2) =   2(1 + \sqrt{2})^2\, \alpha_1^{-1}    + (1 + \alpha_2)\, (\alpha_2(1- \alpha_1))^{-1}.$
Now,  set $\alpha_1 = \alpha_2 = \frac{\alpha} {1+ \alpha}$   which yields
\bes
C (\alpha,\alpha) =  \alpha^{-1} \, [(1 + \sqrt{2})^2 (4 + 2\alpha) + (2 + \alpha)(1+\alpha)] \leq 
\alpha^{-1} \, [5\, (1 + \sqrt{2})^2 + 3.5].
\ees
Set $t=n$. Due to $n \ln K \leq  n K$ and $\ln (ne) \leq 2 \, \ln n$, one obtains
 \bes 
 \norm{P_{*}}_{F}^2 -  (1 + \alpha_n)\,   
\sum_{k,l = 1}^K \norm{  P_{*}^{(k,l)} (Z) }_{op}^2  
 \leq \frac{H}{\alpha_n}\,   (n K + K^2 \ln n). 
\ees
 where 
\bes
H = [5\, (1 + \sqrt{2})^2 + 3.5] \max\{ C_1 +2C_3, 2C_2\}.
\ees 
The latter  contradicts \eqref{eq:cond_check}, since  $\hat{Z} \in \Upsilon(Z_{*},\rho_n)$, which completes the proof.

%%%%%%%%%%%%%%%%%%%%%%%%%%%%%%%%%%%%%%%%%%%%%%%%%%%%%%%%%%%%%%%%%%%%%%%%%%%%%%%%%%%%%%%%%%%%%%%%%%%%%%%%%%%%%%%%%%%%%%%%%%
%%%%%%%%%%%%%%%%%%%%%%%%%%%%%%%%%%%%%%%%%%%%%%%%%%%%%%%%%%%%%%%%%%%%%%%%%%%%%%%%%%%%%%%%%%%%%%%%%%%%%%%%%%%%%%%%%%%%%%%%%%

\subsection{Proofs  of Lemma~\ref{lem:SBM_clust_er} and Corollary~\ref{cor:SBM_clust_er} }  

\noindent
{\it Proof of Lemma~\ref{lem:SBM_clust_er}. } Denote $N = n/2$ and let $s_{12} = N \rho_1$ and $s_{21}= N \rho_2$ 
be, respectively, the number of nodes in class 1 that are placed into class 2 and visa versa. 
Then, $\rho = (\rho_1 + \rho_2)/2$. Let $s_{11} = N(1 -\rho_1)$ and $s_{22} = N(1 - \rho_2)$.

Denote the version of  matrix $P$ arranged according to $Z$ by $\tilP$. Observe that matrix  $\tilP$
has four rank two blocks of sizes $(s_{11} + s_{21})\times (s_{11} + s_{21})$ (top  left), 
$(s_{22} + s_{12})\times (s_{22} + s_{12})$ (bottom right) and $(s_{22} + s_{12})\times (s_{11} + s_{21})$ (bottom left).
% with constant blocks equal to $b$ on the main diagonal and $r$ off the main diagonal. 

Consider diagonal blocks $\tilP_1$ and $\tilP_2$ of the   matrix $\tilP$.
Note that 
\be  \label{eq:norm_dif} 
 \Del = \|P_{*}\|_{F}^2 -  (1 + \alpha)\,  % \underset{Z \in \Upsilon(Z_{*},\rho)}   
\sum_{k,l = 1}^2 \|P_{*}^{(k,l)} (Z)\|_{op}^2  
 \geq \sum_{k=1}^2 \sig_{2,k}^2 - \alpha \sig_{1,k}^2,
  \ee
where $\sig_{1,k} \geq \sig_{2,k} >0$ are two nonzero singular values of matrix $\tilP_k$, $k=1,2$.
% Assume , without loss of generality, that $\rho_2 \neq 0$.
Then, the diagonal sub-matrices $\tilP_k$, $k=1,2$,   of   $\tilP$ have constant blocks   
$b 1_{s_{1k}}  1_{s_{1k}}^T$ and $b 1_{s_{2k}}  1_{s_{2k}}^T$
on the main diagonal and blocks $r 1_{s_{1k}}  1_{s_{2k}}^T$  and $r 1_{s_{2k}}  1_{s_{1k}}^T$  off the main diagonal.
It is easy to check that $\tilP_k \tilP_k^T$  is the symmetric matrix with the constant blocks $A 1_{s_{1k}} 1_{s_{1k}}^T$ 
and $B 1_{s_{2k}} 1_{s_{2k}}^T$ on the main diagonal and blocks $R 1_{s_{1k}}  1_{s_{2k}}^T$  and $R 1_{s_{2k}}  1_{s_{1k}}^T$  off the main diagonal,
$k=1,2$. Here,
\bes 
A = b^2  s_{1k}  + r^2 s_{2k},\quad B = r^2 s_{1k} + b^2 s_{2k}, \quad R = r b (s_{1k}+ s_{2k}).
\ees
Moreover, $\sig_{l,k}^2 = \lam_{l,k}$, $l,k=1,2$,  where $\lam_{l,k}$ is the $l$-th nonzero eigenvalue of matrix
$\tilP_k \tilP_k^T$. 

Apply   Lemma~\ref{lem:more_facts}  with $m_1 = s_{11}$ and $m_2 = s_{21}$; and then again with 
$m_1 = s_{12}$ and $m_2 = s_{22}$.    Direct calculations yield that
\begin{align}  
& A m_1 + B m_2 = b^2 N^2 \lkv 1 - 2 (\rho_1 + \rho_2\, r^2/b^2) + 
(\rho_1^2 + \rho_2^2 - 2 \rho_1 \rho_2\,  r^2/b^2)\rkv, \nonumber \\
& \label{eq:sum}\\
& m_1 m_2 (A B - R^2) = N^4 (1 - \rho_1)^2 \rho_2^2\, (b^2 - r^2)^2, \nonumber
\end{align} 
where the first expression in \eqref{eq:sum} is true for $\tilP_1 \tilP_1^T$ and the second one is true in both cases.
In order to obtain an expression for   $\tilP_2 \tilP_2^T$, one needs to interchange 
$\rho_1$ and $\rho_2$ in the first formula in \eqref{eq:sum}. 
It is easy to check that, if $\rho   <\min(1, r^2/b^2)$, then 
$A m_1 + B m_2 \leq b^2 N^2$ in both cases. Hence, it follows from \eqref{eq:lam_dif} that, 
for $\Del$   defined in \eqref{eq:norm_dif}, one has
\bes 
\Del  \geq C_{\lam}\, \lfi N^2\, (b^2 - r^2)^2 \, b^{-2}\, [(1 -\rho_1)^2 \rho_2^2 + (1 -\rho_2)^2 \rho_1^2] 
- 2 \alpha b^2 N^2 \rfi, 
\ees
where $\rho_1, \rho_2 \geq 0$ and  $\rho_1 + \rho_2 = \rho$. 
It is easy to check that minimum of the right hand side  
occurs when $\rho_1 = \rho_2 = \rho/2$, so that
\bes 
\Del  \geq C_{\lam}\, \lfi  N^2 \rho^2\, (b^2 - r^2)^2 \, b^{-2}/2  
- 2 \alpha b^2 N^2 \rfi. 
\ees
To complete  the proof, recall that $N = n/2$. 
 \\

\medskip

%%%%%%%%%%%%%%%%%%%%%%%%%%%%%%%%%%%%%%%%%%%%%%%%%%%%%%%%%%%%%%%%%%%%%%%%%%%%%%%%%%%%%%%%%%%%%%%%%%%%%%%%%%%%%%%%%%%%%%%%%%

\noindent
{\it Proof of Corollary~\ref{cor:SBM_clust_er}. } Since $K=2$, in order to apply Theorem~\ref{th:clust}, one needs 
$\Del \gtrsim n/\alpha$   in \eqref{eq:norm_dif}, which, due to Lemma~\ref{lem:SBM_clust_er}, is equivalent to
\be \label{eq:ineq_rho} 
\rho^2 (1-r^2/b^2)^2 - 4 \alpha \gtrsim (\alpha n b^2)^{-1}.
\ee 
Inequality \eqref{eq:ineq_rho} implies that $\alpha \asymp \rho^2 (1-r^2/b^2)^2$ which, in turn, leads to
\bes
\rho^2 (1-r^2/b^2)^2 \lesssim [n \rho^2 b^2  (1-r^2/b^2)^2]^{-1}. 
\ees
The latter is equivalent to \eqref{eq:rho_bound}.

%%%%%%%%%%%%%%%%%%%%%%%%%%%%%%%%%%%%%%%%%%%%%%%%%%%%%%%%%%%%%%%%%%%%%%%%%%%%%%%%%%%%%%%%%%%%%%%%%%%%%%%%%%%%%%%%%%%%%%%%%%
%%%%%%%%%%%%%%%%%%%%%%%%%%%%%%%%%%%%%%%%%%%%%%%%%%%%%%%%%%%%%%%%%%%%%%%%%%%%%%%%%%%%%%%%%%%%%%%%%%%%%%%%%%%%%%%%%%%%%%%%%%

\subsection {Proof  of Theorem~\ref{th:correctness}}
For simplicity, we re-order matrix $P_*$, so columns/rows corresponding to the same community 
are consecutive. Denote corresponding column sub-matrices of  $P_*$ by  $P^{(k)}_*$, $k=1, \ldots, K$.
Here, each  $P^{(k)}_*$ is the concatenation of sub-matrices $P^{(j,k)}_*$ with $j=1, \ldots, K$.

We start the proof with the definition. Following \cite{pmlr-v51-wang16b}, we say that 
 $P^{(k)}_*$   is {\it in general position} if for
all $l$, $1 \leq l \leq d_k$, where $d_k$ is the dimension of the column space of  $P^{(k)}_*$,   any subset of $l$  columns 
in $P^{(k)}_*$  are linearly independent. We say that  $P_*$ is in
general position if $P^{(k)}_*$ is in general position for all
$k=1, \ldots, K$.

Due to Theorem 1 of  \cite{Elhamifar:2013:SSC:2554063.2554078},
the similarity matrix $S$ obtained as a solution of 
optimization problem \eqref{mn:opt_prob1} satisfies the SEP. 
According to Theorem B.2 of \cite{pmlr-v51-wang16b}, Algorithm 1 partitions 
points into a minimal union-of-subspace structure where points in each subspace 
are in general position. Hence, by the arguments identical to those in the proof of 
Theorem 3.1 of \cite{pmlr-v51-wang16b}, Algorithm 1    recovers the correct partition 
up to a permutation $\pi$ on $\{1, \ldots, K\}$.
\\ 

%%%%%%%%%%%%%%%%%%%%%%%%%%%%%%%%%%%%%%%%%%%%%%%%%%%%%%%%%%%%%%%%%%%%%%%%%%%%%%%%%%%%%%%%%%%%%%%%%%%%%%%%%%%%%%%%%%%%%%%%%%%%%%
%%%%%%%%%%%%%%%%%%%%%%%%%%%%%%%%%%%%%%%%%%%%%%%%%%%%%%%%%%%%%%%%%%%%%%%%%%%%%%%%%%%%%%%%%%%%%%%%%%%%%%%%%%%%%%%%%%%%%%%%%%%%%%

%%%%%%%%%%%%%%%%%%%%%%%%%%%%%%%%%%%%%%%%%%%%%%%%%%%%%%%%%%%%%%%%%%%%%%%%%%%%%%%%%%%%%%%%%%%%%%%%%%%%%%%%%%%%%%%%%%%%%%%%%%
%%%%%%%%%%%%%%%%%%%%%%%%%%%%%%%%%%%%%%%%%%%%%%%%%%%%%%%%%%%%%%%%%%%%%%%%%%%%%%%%%%%%%%%%%%%%%%%%%%%%%%%%%%%%%%%%%%%%%%%%%%

\subsection{Supplementary statements and their proofs}

\begin{lem}  \label{lem:lowrank_approx}
For any matrices $A, B \in \RR^{m \times n}$ and any unit vectors $u \in \RR^m$ and $v \in \RR^n$, let 
\be \label{eq:Pi_uv}
\Pi_{u,v}(A) =  (uu^T)A(vv^T)   
\ee 
denote the projection of matrix $A$ on the vectors $(u,v)$. Then,
\be \label{eq:Inner_Prod_Pi_uv}
\langle \Pi_{u,v}(B), A - \Pi_{u,v}(A) \rangle =0.
\ee
Furthermore, if we let $\hat{u}$ and $\hat{v}$ be the singular vectors of matrix $A$ corresponding to its largest singular value $\sig$, the best rank one approximation of $A$ is given by 
\be \label{eq:Pi_uv_A}
\Pi_{\hat{u},\hat{v}}(A) = (\hat{u} \hat{u}^T)A(\hat{v}\hat{v}^T)=  \sigma \hat{u}  \hat{v} ^T.
\ee 
\end{lem}

\begin{lem}  \label{lem:Pi_orth}
Let $(\hat{u}, \hat{v})$ and $(u,v)$ denote the pairs of singular vectors of matrices $A$ and $P$, 
respectively, corresponding to their largest singular values. Then,
\be \label{eq:Pi_orth}
\left\|\Pi_{u,v} (P) - P\right\|_F \leq \left\|\Pi_{\hat{u},\hat{v}} (P) - P\right\|_F  \leq \left\|\Pi_{\hat{u},\hat{v}} (A) - P\right\|_F
\ee 
where $\Pi_{u,v} (\cdot)$ is defined in \eqref{eq:Pi_uv}.
\end{lem}

 %%%%%%%%%%%%%%%%%%%%%%%%%%%%

\noindent{\bf Proof.}
The first inequality in \eqref{eq:Pi_orth}  is true because $\Pi_{u,v} (P)$ is the best rank one approximation of $P$. Now let $A=P+\Xi$. Then 
 \bes
\begin{aligned}
 \left\|\Pi_{\hat{u},\hat{v}} (A) - P\right\|_F^2  
% = \|\Pi_{\hat{u},\hat{v}} (P + \Xi) - P\|_F^2\\
= \left\|\Pi_{\hat{u},\hat{v}} (P) - P + \Pi_{\hat{u},\hat{v}} (\Xi)\right\|_F^2 =  \left\|\Pi_{\hat{u},\hat{v}} (P) - P\right\|_F^2 +\left \|\Pi_{\hat{u},\hat{v}} (\Xi)\right\|_F^2
\end{aligned}
\ees 
which leads to the second inequality in  \eqref{eq:Pi_orth}.
 
%%%%%%%%%%%%%%%%%%%%%%%%%%%%%%%%%%%%%%%%%%%%%%%%%%%%%%%%%%%%%%%%%%%%%%%%%%%%%%%%%%%%%%%%%%%%%%%%%%%%%%%%%%%%%%%%%%%%%%%%%%
%%%%%%%%%%%%%%%%%%%%%%%%%%%%%%%%%%%%%%%%%%%%%%%%%%%%%%%%%%%%%%%%%%%%%%%%%%%%%%%%%%%%%%%%%%%%%%%%%%%%%%%%%%%%%%%%%%%%%%%%%%

\begin{lem} \label{lem:prob_bound_error}
Let    elements of matrix $\Xi \in (-1,1)^{n \times n}$ be independent Bernoulli errors. 
Let matrix   $ \Xi$ be partitioned into $K^2$ sub-matrices $\Xi^{(k,l)}$,  $k,l = 1, \cdots, K$.
Then, for any $x >0$  
\be \label{eq:lem3}
\PP \left\{ \sum_{k,l =1}^K \norm{\Xi^{(k,l)}}_{op}^2  \leq  C_1 nK + C_2 K^2 \ln (ne) + C_3x \right\}   \geq  1 - \exp (-x),
\ee 
where $C_1,C_2$ and $C_3 $ are absolute constants independent of $n$ and $K$. 
\end{lem}

%%%%%%%%%%%%%%%%%%%%%%%%%%%%%%%%%%%%%%%%%%%%%%%%%%%%%%%%%%%%%%%%%%%%%%%%%%%%%%

\noindent{\bf Proof.} 
Consider vectors $\xi$ and $\mu$ with elements $\xi_{k,l} = \|\Xi^{(k,l)}\|_{op}$ and 
$\mu_{k,l} =\EE \|\Xi^{(k,l)}\|_{op}$, $k,l = 1,\cdots, K$,   and let $\eta = \xi - \mu$.
Then, 
\be \label{norm_xi}
\Delta = \sum_{k,l =1}^K \norm{\Xi^{(k,l)}}_{op}^2  =  \| \xi \|^2 \leq  2 \| \eta \|^2 + 2 \| \mu \|^2
\ee 
Hence,  we need to construct the upper bounds for $\| \eta \|^2$ and $\| \mu \|^2$.

We start with constructing upper bounds for  $\| \mu \|^2$.
Let $\Xi_{i,j}^{(k,l)}$  be elements of the $(n_k \times n_l)$-dimensional matrix $\Xi^{(k,l)}$.
Then, $\EE(\Xi_{i,j}^{(k,l)}) = 0$ and, by Hoeffding's inequality,
$\EE \lfi \exp(\lambda \Xi_{i,j}^{(k,l)})\rfi \leq \exp\left( \lambda^2/8 \right)$.
Taking into account that Bernoulli errors are bounded by one in absolute value and applying 
Corollary 3.3 of  \cite{bandeira2016}  with $m = n_k$, $n = n_l $, $\sigma_{*}=1$, $\sigma_{1} = \sqrt{n_l}$ and 
$\sigma_2 = \sqrt{n_k}$, obtain
$$ 
\mu_{k,l}    \leq C_0(\sqrt{n_k} + \sqrt{n_l} + \sqrt{\ln(n_k \wedge n_l)} )
$$
where $C_0$ is an absolute constant independent of $n_k$ and $n_l$.
Therefore, 
\be\label{NORM_MU} 
\|\mu\|^2 \leq 3C_0^2 \sum_{k,l =1}^K (n_k + n_l + \ln(n_k \wedge n_l)) 
\leq 6C_0^2 nK + 3C_0^2 K^2 \ln n.
\ee

Next, we show that, for any fixed partition,   $\eta_{k,l} = \xi_{k,l} - \mu_{k,l}$ are independent sub-gaussian random variables
when $1\leq k \leq  l \leq K$. Independence follows from the conditions of Lemma ~\ref{lem:prob_bound_error}. To prove the 
sub-gaussian property, use Talagrand's concentration inequality (Theorem 6.10 of \cite{Boucheron2013}):
if $\Xi_1, \Xi_2,\Xi_3, \cdots , \Xi_n$ are independent random variables taking values in the interval $[0,1]$ and 
$f : [0,1]^n \rightarrow R$ is a separately convex function such that 
$|f(x) - f(y) | \leq \|x-y\|$  for  all  $x,y \in [0,1]^n$, 
then, for $Z = f(\Xi_1,\Xi_2,\Xi_3, \cdots , \Xi_n)$  and any $t > 0$, one has    
\be \label{eq:Talag}
\PP(Z > \EE Z +t) \leq \exp ( - t^2/2).  
\ee

Apply this theorem to vectors $\zeta_{k,l} = \vect (\Xi^{(k,l)})  \in [0,1]^{n_k \times n_l}$ and 
$f (\Xi^{(k,l)}) = f(\zeta_{k,l}) =  \norm{\Xi^{(k,l)}}_{op} $.
Note that, for any two matrices $\Xi$ and  $ \tilde{\Xi}$ of the same size, one has 
$\|\Xi - \tilde{\Xi} \|_{op}^2 \leq \|{\Xi - \tilde{\Xi}}\|_{F}^2  = \| \vect(\Xi) - \vect(\tilde{\Xi})\|^2$.
Then, applying Talagrand's inequality with  $Z = \|{\Xi^{(k,l)}}\|_{op}$ and $Z = -\|{\Xi^{(k,l)}}\|_{op}$,  obtain
\ignore{
\bes
 \PP\left(  \norm{\Xi^{(k,l)}}_{op} - \EE \norm{\Xi^{(k,l)}}_{op} > t \right)  \leq  \exp \left(\frac{ -t^2}{2} \right)
\ees
\bes
 \PP\left(  - \norm{\Xi^{(k,l)}}_{op} + \EE \norm{\Xi^{(k,l)}}_{op} > - t \right)  \leq  \exp \left(\frac{ -t^2}{2} \right)
\ees
is equivalent to 
\bes
 \PP\left(   \norm{\Xi^{(k,l)}}_{op} - \EE \norm{\Xi^{(k,l)}}_{op} <  - t \right)  \leq  \exp \left(\frac{ -t^2}{2} \right)
\ees
so that 
}
\bes 
 \PP\left( \left| \|{\Xi^{(k,l)}}\|_{op} - \EE \|{\Xi^{(k,l)}}\|_{op} \right| > t \right)  \leq  2 \exp (-t^2/2).
\ees
Now, use the Lemma 5.5 of  \cite{vershynin_2012}  which states that the latter implies that for any $t>0$ 
and some absolute constant $C_4>0$, 
\be \label{vershynin-app}
\EE\left[ \exp(t \eta_{k,l}) \right]  = \EE\left[ \exp(t(\xi_{k,l} - \mu_{k,l})) \right] \leq \exp (C_4 t^2/2),\ \  C_4 > 0.
\ee
Hence, $\eta_{k,l}$ are independent   sub-gaussian  random variables when $1\leq k \leq  l \leq K$.

Now, we obtain  an  upper bound for $\|{\eta}\|^2$.
Use Theorem 2.1 of \cite{Hsu2011} which states that for any matrix $A$, if for some $\sigma > 0$ and any vector  $h$
one has $\EE[\exp(h^T \tilde{\eta} )] \leq \exp  (\|{h}\|^2 \sigma^2/2)$,  then, for any $x > 0$,
\be \label{eq:hsu}
 \PP \lfi \| A\tilde{\eta} \|^2  \geq \sigma^2 (\Tr(A^TA) + 2\sqrt{\Tr((A^TA)^2)\, x} + 2 \|{A^TA}\|_{op}\, x)\rfi  \leq   \exp (-x).
\ee  
Applying  \eqref{eq:hsu} with  $A =  I_{K(K+1)/2}$ and $\sigma^2 = C_4$ to a sub-vector  
$\tilde{\eta}$ of $\eta$ which contains components $\eta_{k,l}$ with $1\leq k \leq  l \leq K$,
obtain 
\bes
 \PP\left\{\|\tilde{\eta}\|^2  \geq C_4 \lkr K(K+1)/2   +  \sqrt{2\, K(K+1)\, x} + 2 x \rkr \rfi \leq   \exp (-x).
\ees
Since $\|{\eta}\|^2  \leq  2 \|{\tilde{\eta}}\|^2 $, derive 
\be\label{bound_eta}
 \PP\left\{\norm{\eta}^2  \geq 2C_4 K(K+1) + 6C_4x   \right\} \leq   \exp \left(-x \right)
\ee
Combination of formulas \eqref{norm_xi}  and \eqref{bound_eta} yield 
\bes
 \PP\left\{\norm{\xi}^2  \leq  2 \norm{\mu}^2  + 4C_4 K(K+1) + 12C_4x   \right\} \geq   1-  \exp \left(-x \right)
\ees

Plugging in $\norm{\mu}^2$  from \eqref{NORM_MU} into the last inequality, derive for any $x>0$ that
\be\label{main_res}
 \PP\left\{  \norm{\xi}^2  \leq  12C_0^2 nK + 6C_0^2 K^2 \ln n + 4C_4 K(K+1) + 12C_4x   \right\} \geq   1-  \exp \left(-x \right).
\ee
Since $K(K+1) \leq  2K^2$ and $6C_0^2 K^2 \ln n +  8C_4 K^2 \leq \max(6C_0^2, 8C_4)  K^2 \ln (ne)$, inequality 
 \eqref{eq:lem3} holds with $C_1 = 12C_0^2$, $ C_2 = \max(6C_0^2, 8C_4)$  and  $C_3 = 12 C_4$.

%%%%%%%%%%%%%%%%%%%%%%%%%%%%%%%%%%%%%%%%%%%%%%%%%%%%%%%%%%%%%%%%%%%%%%%%%%%%%%%%%%%%%%%%%%%%%%%%%%%%%%%%%%%%%%%%%%%%%%%%%%
%%%%%%%%%%%%%%%%%%%%%%%%%%%%%%%%%%%%%%%%%%%%%%%%%%%%%%%%%%%%%%%%%%%%%%%%%%%%%%%%%%%%%%%%%%%%%%%%%%%%%%%%%%%%%%%%%%%%%%%%%%

\begin{lem}\label{lem:prob_error_bound_main}
For any $ t > 0$,
\be \label{rr:prob_bound_Xi_case1} 
\PP \left\{   \di\sum_{k,l = 1}^{\hat{K}} \norm{\Xi^{(k,l)}(\hat{Z},\hat{K})}_{op}^2  -F_1(n,\hat{K})  \leq   C_3 t \right\} \geq  1 - \exp{(-t)}.
\ee
where $F_1(n,K)$ is given by \eqref {eq:F1}.
%$F_1(n,\hat{K}) =   C_1n\hat{K} + C_2 \hat{K}^2 \ln (ne) +C_3( \ln n + n \ln \hat{K}) $
\end{lem}

%%%%%%%%%%%%%%%%%%%%%%%%%%%%

\noindent{\bf Proof. }
Using Lemma \ref{lem:prob_bound_error}, for any fixed $K$ and $Z \in \calM_{n,K}$, obtain 
\bes
\PP \left\{  \di\sum_{k,l = 1}^K \|{\Xi^{(k,l)}(Z,K)}\|_{op}^2 -  C_1nK - C_2 K^2 \ln (ne) - C_3x  \geq 0 \right\} \leq  \exp{(-x)}.
\ees
Application of the union bound  over  $Z \in \calM_{n,K}$ and  $K \in [1,n] $  and setting $x = t + \ln n + n \ln K $ yields
\begin{align*}
& \PP \left\{  \di\sum_{k,l = 1}^{\hat{K}} \|{\Xi^{(k,l)}(\hat{Z},\hat{K})}\|_{op}^2 -  
C_1 n\hat{K} - C_2 \hat{K}^2 \ln (ne) - C_3t - C_3 \ln n - C_3 n \ln \hat{K}  \geq 0 \right\}  \\
\leq &  
\PP \left\{ \underset{1 \leq K \leq n}{\max}\ \underset{Z \in \mathcal{M}_{n,K}}{\max}  \left( \di \sum_{k,l = 1}^K  \|{\Xi^{(k,l)}(Z,K)}\|_{op}^2
- F_1(n,K) \right) \geq C_3 t\right\} \\
\leq &
\sum_{k=1}^n\   \sum_{Z\in \mathcal{M}_{n,K }} \PP \left\{ \di\sum_{k,l = 1}^K  \norm{\Xi^{(k,l)}(Z,K)}_{op}^2 
- F_1(n,K) \geq C_3 t\right\}
% \\
\leq 
% & 
n K^n \exp\{-t - \ln n - n \ln K\} = \exp{(-t)},
\end{align*}
which completes the proof.
\\

%%%%%%%%%%%%%%%%%%%%%%%%%%%%%%%%%%%%%%%%%%%%%%%%%%%%%%%%%%%%%%%%%%%%%%%%%%%%%%%%%%%%%%%%%%%%%%%%%%%%%%%%%%%%%%%%%%%%%%%%%%
%
%%%%%%%%%%%%%%%%%%%%%%%%%%%%%%%%%%%%%%%%%%%%%%%%%%%%%%%%%%%%%%%%%%%%%%%%%%%%%%%%%%%%%%%%%%%%%%%%%%%%%%%%%%%%%%%%%%%%%%%%%%%
% THIS IS THE LEMMA EQUIVALENT TO ABOVE LEMMA FOR THEOREM 2

\begin{lem}\label{lem:prob_error_bound_mainS}
For any $ t > 0$ and fixed $K$,
\be \label{rr:prob_bound_Xi_case1S} 
\PP \left\{   \di\sum_{k,l = 1}^{K} \norm{\Xi^{(k,l)}(\hat{Z})}_{op}^2  - F_3(n)  \leq   C_3 t \right\} \geq  1 - \exp{(-t)}.
\ee
where $F_3(n)$ is given by \eqref {eq:F3S}.
\end{lem}

%%%%%%%%%%%%%%%%%%%%%%%%%%%%

\noindent{\bf Proof. }
Using Lemma \ref{lem:sparse_prob_bound_error}, for any fixed $Z \in \calM_{n,K}$, obtain 
\bes
\PP \left\{  \di\sum_{k,l = 1}^K \|{\Xi^{(k,l)}(Z)}\|_{op}^2 -  C_1 \tau_n K - C_2 K^2  - C_3x  \geq 0 \right\} \leq  \exp{(-x)}.
\ees
Application of the union bound  over  $Z \in \calM_{n,K}$  and setting $x = t + n \ln K $ yields
\begin{align*}
& \PP \left\{  \di\sum_{k,l = 1}^{K} \|{\Xi^{(k,l)}(\hat{Z})}\|_{op}^2 -  
C_1 \tau_n K - C_2 K^2  - C_3 n \ln K - C_3t   \geq 0 \right\}  \\
\leq &  
\PP \left\{  \underset{Z \in \mathcal{M}_{n,K}}{\max}  \left( \di \sum_{k,l = 1}^K  \|{\Xi^{(k,l)}(Z)}\|_{op}^2
- F_3(n) \right) \geq C_3 t\right\} \\
\leq &
  \sum_{Z\in \mathcal{M}_{n,K }} \PP \left\{ \di\sum_{k,l = 1}^K  \norm{\Xi^{(k,l)}(Z)}_{op}^2 
- F_3(n) \geq C_3 t\right\}
% \\
\leq 
% & 
 K^n \exp\{-t - n \ln K\} = \exp{(-t)},
\end{align*}
which completes the proof.

%%%%%%%%%%%%%%%%%%%%%%%%%%%%%%%%%%%%%%%%%%%%%%%%%%%%%%%%%%%%%%%%%%%%%%%%%%%%%%%%%%%%%%%%%%%%%%%%%%%%%%%%%%%%%%%%%%%%%%%%%%%
%
%%%%%%%%%%%%%%%%%%%%%%%%%%%%%%%%%%%%%%%%%%%%%%%%%%%%%%%%%%%%%%%%%%%%%%%%%%%%%%%%%%%%%%%%%%%%%%%%%%%%%%%%%%%%%%%%%%%%%%%%%%

% This lemma is an equivalent of Bandeira and Handel

\begin{lem} \label{lem:Bern_small_p}
Let $P, A \in [0,1]^{N \times M}$ where $A_{i,j} \sim  Bernoulli(P_{i,j})$ are independent
and $\|P\|_{\infty} \leq \tau \equiv \tau_{M,N}$. Let $\Xi = A -P$.  
If $N+M$ be large enough, so that 
\be \label{eq:large_NM}
\log(N+M) \leq (N+M)^{2/13},
\ee
then, for some absolute constant $C>0$ one has 
\be \label{eq:Expect_First}
\EE \|\Xi\|_{op} \leq C \lkr \sqrt{\tau(M+N)}  +  \sqrt{\ln(M+N)} \rkr.
\ee
Moreover, if $\tau \geq C_{\tau} \, (M+N)^{-1}\, \ln(M+N)$  for some absolute constant $C_{\tau}  >0$,
then 
\be \label{eq:Expect_Second}
\EE \|\Xi\|_{op} \leq C \,  \sqrt{\tau(M+N)}.
\ee
\end{lem}

%%%%%%%%%%%%%%%%%%%%%%%%%%%%

\noindent{\bf Proof. }
First, symmetrize matrix $\Xi$, similarly to \cite{bandeira2016}, by replacing it by the 
$(N+M)\times (N+M)$ matrix $\tilXi$ with zero sub-matrices on the main diagonal 
and matrices $\Xi$  and $\Xi^T$ as the bottom left and top right sub-matrices.
Then, $\|\Xi\|_{op} = \|\tilXi\|_{op}$ and 
\bes 
\max_i \sum_{j} \EE(\tilXi_{i,j})^2 = \max_i \lkr \sum_{j=1}^N P_{i,j}(1 - P_{i,j}),
 \sum_{j=1}^M P_{j,i}(1 - P_{j,i})\rkr \leq \tau(M+N).
\ees
Set 
$H = (\tau(N+M))^{-1/2}\, \tilXi$.   Obtain
\bes
\max_i \sum_{j} \EE(H_{i,j})^2 \leq 1, \quad
\max_{i,j} \EE (H_{i,j}^2) \leq  (M+N)^{-1}.
\ees 
Now, apply Theorem~2.6 of \cite{benaychgeorges2017spectral} to matrix $H$
with $\kappa =1$, $n = N+M$ and $q = \min(\sqrt{\tau(M+N)}, \sqrt{\ln(M+N)})$.
Note that if $N+M$ is large enough, so \eqref{eq:large_NM} holds, 
Theorem~2.6 of \cite{benaychgeorges2017spectral}  yields 
$\EE \|H\|_{op} \leq 2 + C  \sqrt{\ln(M+N)}\,/q$. To complete the 
proof, recall that 
$\EE \|\Xi\|_{op} =  \sqrt{\tau(M+N)}\, \EE \|H\|_{op}$.
\\

%%%%%%%%%%%%%%%%%%%%%%%%%%%%%%%%%%%%%%%%%%%%%%%%%%%%%%%%%%%%%%%%%%%%%%%%%%%%%%%%%%%%%%%%%%%%%%%%%%%%%%%%%%%%%%%%%%%%%%%%%%
%%%%%%%%%%%%%%%%%%%%%%%%%%%%%%%%%%%%%%%%%%%%%%%%%%%%%%%%%%%%%%%%%%%%%%%%%%%%%%%%%%%%%%%%%%%%%%%%%%%%%%%%%%%%%%%%%%%%%%%%%%

\begin{lem} \label{lem:sparse_prob_bound_error}
Let    elements of matrix $\Xi \in (-1,1)^{n \times n}$ be independent Bernoulli errors. 
Let matrix   $ \Xi$ be partitioned into $K^2$ sub-matrices $\Xi^{(k,l)}$,   of sizes $n_k \times n_l$, $k,l = 1, \cdots, K$,
where $\sum n_k =n$.
Assume that $n_{\min} = \min_k (n_k)$ is large enough, so that  $\log(2 n_{\min}) \leq (2 n_{\min})^{2/13}$.
If $\tau_n \geq C_{\tau}\, \log(2 n_{\min})/(2 n_{\min})$ for some absolute constant $C_{\tau}  >0$,  then, for any $x >0$ 
\be \label{eq:sparse_prob_bou_er}
\PP \left\{ \sum_{k,l =1}^K \norm{\Xi^{(k,l)}}_{op}^2  \leq  C_1 \tau_n n K + C_2 K^2   + C_3 x \right\}   \geq  1 - \exp (-x),
\ee 
where $C_1,C_2$ and $C_3 $ are absolute constants independent of $n$ and $K$ and $\tau_n = \di \max_{i,j} P_{i,j}$.
\end{lem}

%%%%%%%%%%%%%%%%%%%%%%%%%%%%

\noindent{\bf Proof. } 
The proof is identical to the proof of Lemma~\ref{lem:prob_bound_error} with the only exception that 
$\|\mu\|^2 \leq C \tau_n n K$ due to Lemma~\ref{lem:Bern_small_p}. 
\\

%%%%%%%%%%%%%%%%%%%%%%%%%%%%%%%%%%%%%%%%%%%%%%%%%%%%%%%%%%%%%%%%%%%%%%%%%%%%%%%%%%%%%%%%%%%%%%%%%%%%%%%%%%%%%%%%%%%%%%%%%%
%%%%%%%%%%%%%%%%%%%%%%%%%%%%%%%%%%%%%%%%%%%%%%%%%%%%%%%%%%%%%%%%%%%%%%%%%%%%%%%%%%%%%%%%%%%%%%%%%%%%%%%%%%%%%%%%%%%%%%%%%%

\begin{lem} \label{lem:simple_facts}
Let $I_m$  be   the identity matrix of size $m$ and $E_m = 1_m 1_m^T$.
Let $a$, $b$  and  $\lam$ be nonzero scalars.  
Then,  
\be \label{eq:simple}
(a E_m + b I_m)^{-1} = \frac{1}{b} I_m - \frac{a}{b\, (am + b)}\, E_m,
\quad 
\det(a E_m - \lam I_m) = (-1)^m \lam^{m-1} (\lam - m a)
\ee
\end{lem}

%%%%%%%%%%%%%%%%%%%%%%%%%%%%

\noindent{\bf Proof. } 
The first statement in \eqref{eq:simple} can be verified directly. 
The second statement can be proved by induction. Indeed, the statement is correct for $m=1$.
Assume that \eqref{eq:simple} is true for $m=M$. By Theorem~1.2.6  of  Gupta and Nagar \cite{gupta_nagar},
one has
\be \label{eq:sim1} 
\det(a E_{M+1} - \lam I_{M+1}) = S \det(a E_{M} - \lam I_{M})
\ee 
where 
\beqns  
S  &=& a - \lam - a^2 \, 1_M^T (a E_{M} -  \lam I_{M})^{-1}  1_M  \\
   &= & a - \lam - a^2 M/(aM - \lam) = (\lam^2 - a\lam(M+1))/(aM - \lam).
\eeqns 
Plugging $S$ into \eqref{eq:sim1} obtain that
$\det(a E_{M+1} - \lam I_{M+1}) = (-1)^{M+1} \lam^{M} (\lam  - a\,(M+1))$,
which completes the proof. 
\\

%%%%%%%%%%%%%%%%%%%%%%%%%%%%%%%%%%%%%%%%%%%%%%%%%%%%%%%%%%%%%%%%%%%%%%%%%%%%%%%%%%%%%%%%%%%%%%%%%%%%%%%%%%%%%%%%%%%%%%%%%%
%%%%%%%%%%%%%%%%%%%%%%%%%%%%%%%%%%%%%%%%%%%%%%%%%%%%%%%%%%%%%%%%%%%%%%%%%%%%%%%%%%%%%%%%%%%%%%%%%%%%%%%%%%%%%%%%%%%%%%%%%%

\begin{lem} \label{lem:more_facts}
Let $A$,$B$ and $R$ be arbitrary scalars and $E_{k,l} = 1_k 1_l^T$. 
Consider a symmetric matrix $X$ with constant  blocks $A E_{m_1}$ and 
$B E_{m_2}$  on the main diagonal and $R E_{m_2, m_1}$ off the main diagonal. 
Then the only two non-zero eigenvalues  of $X$ are $\lam_1 > \lam_2 >0$ where
\be \label{eq:eigenvalues}
\lam_{1,2} = \frac{A m_1 + B m_2}{2} \pm \frac{1}{2} \sqrt{(A m_1 + B m_2)^2 - 4 m_1 m_2 (A B - R^2)}.
\ee
Moreover,   for $0 < \alpha < 1,$ one has
\be \label{eq:lam_dif}
\lam_2 - \alpha \lam_1 \geq %= C_{\lam}\, 
\frac{m_1 m_2 (A B - R^2) - \alpha (A m_1 + B m_2)^2}{A m_1 + B m_2}
% \quad \mbox{with} \quad C_{\lam} \in (1,9).
\ee
\end{lem}

%%%%%%%%%%%%%%%%%%%%%%%%%%%%

\noindent{\bf Proof. } 
In order to find the eigenvalues, we find $\det(X - \lam I)$ and equate it to zero. 
Using Theorem~1.2.6  of  Gupta and Nagar \cite{gupta_nagar} and Lemma~\ref{lem:simple_facts}, we write 
$\det(X - \lam I) = S_1 S_2$, where
\bes
S_1   =   \det(A E_{m_1}  - \lam I_{m_1}), \quad 
S_2   =   \det(B  E_{m_2} - \lam I_{m_2} - R^2 E_{m_2, m_1}  (A E_{m_1}  - \lam I_{m_1})^{-1} E_{m_1, m_2}).
\ees
According to Lemma~\ref{lem:simple_facts}, 
\bes
S_1   = (-1)^{m_1} \lam^{m_1 -1} (\lam - m_1 A), \quad 
S_2 = (-1)^{m_2} \lam^{m_2 -1} (\lam - m_2 \tilA)
\ees
where $\tilA = [B \lam - (AB - R^2) m_1)]/(\lam - m_1 A)$.
Finally, combining all terms, obtain that
\bes
\det(X - \lam I)  = (-1)^{m_1 + m_2} \, \lam^{m_1+m_2 -2} \,
\lkv (\lam - m_1 A) \lam - B \lam m_2 + m_1 m_2  (AB - R^2)\rkv,
\ees
so that the two nonzero eigenvalues of $X$ are the solutions of the quadratic equation
\bes
\lam^2 - (A m_1 + B  m_2)\lam +  m_1 m_2  (AB - R^2)=0,
\ees
and, hence, are of the form \eqref{eq:eigenvalues}.

In order to prove the second statement, observe that
\bes
\lam_2 - \alpha \lam_1 = \frac{2 \lkv (1 + \alpha)^2 \, m_1 m_2 (A B - R^2) - \alpha (A m_1 + B m_2)^2 \rkv}
{(1 - \alpha)\,(A m_1 + B m_2) + (1 + \alpha)\, \sqrt{(A m_1 + B m_2)^2 - 4  m_1 m_2 (A B - R^2)}\,   }
\ees
which yields \eqref{eq:lam_dif}. % if $0 < \alpha < 1/2$.

%%%%%%%%%%%%%%%%%%%%%%%%%%%%%%%%%%%%%%%%%%%%%%%%%%%%%%%%%%%%%%%%%%%%%%%%%%%%%%%%

\ignore{
 \section*{acknowledgements}
All three authors of the paper were  partially supported by National Science Foundation
(NSF)  grant DMS-1712977.
%\printendnotes

% Submissions are not required to reflect the precise reference formatting of the journal (use of italics, bold etc.), however it is important that all key elements of each reference are included.
%\bibliography{PABM}
}

\end{document}

%% file: Def1.tex
%%%%%%%%%%%%%%%%%%%%%%%%%%%%%%%%%%%%%%%%%%%%%

\newcommand{\be}{\begin{equation}}
\newcommand{\ee}{\end{equation}}
\newcommand{\bes}{\begin{equation*}}
\newcommand{\ees}{\end{equation*}}
\newcommand{\beqn}{\begin{eqnarray}}
\newcommand{\eeqn}{\end{eqnarray}}
\newcommand{\beqns}{\begin{eqnarray*}}
\newcommand{\eeqns}{\end{eqnarray*}}

\newcommand{\lkr}{\left(}
\newcommand{\lkv}{\left[}
\newcommand{\rkv}{\right]}
\newcommand{\rkr}{\right)}
\newcommand{\lfi}{\left\{}
\newcommand{\rfi}{\right\}}
\newcommand{\di}{\displaystyle}

\newcommand{\bm}[1]{{\mbox{\mathversion{bold}$#1$}}}

\newcommand{\R}{\mathbb{R}}

\newcommand{\vart}{\vartheta}
\newcommand{\ro}{\varrho}
\newcommand{\ph}{\varphi}

\newcommand{\Del}{\Delta}
\newcommand{\dn}{\delta_n}
\newcommand{\af}{\alpha}
\newcommand{\eps}{\epsilon}
\newcommand{\Ga}{\Gamma}
\newcommand{\ga}{\gamma}
\newcommand{\te}{\theta}
\newcommand{\om}{\omega}
\newcommand{\lam}{\lambda}
\newcommand{\Up}{\Upsilon}
\newcommand{\up}{\upsilon}
\newcommand{\dz}{\zeta}
\newcommand{\sig}{\sigma}
\newcommand{\sgmd}{\sigma^2}
\newcommand{\Lam}{\Lambda}
\newcommand{\Om}{\Omega}

\newcommand{\EE}{\ensuremath{{\mathbb E}}}
\newcommand{\JJ}{\ensuremath{{\mathbb J}}}
\newcommand{\II}{\ensuremath{{\mathbb I}}}
\newcommand{\ZZ}{\ensuremath{{\mathbb Z}}}
\newcommand{\PP}{\ensuremath{{\mathbb P}}}
\newcommand{\QQ}{\ensuremath{{\mathbb Q}}}
\newcommand{\KK}{\ensuremath{{\mathbb K}}}
\newcommand{\RR}{\ensuremath{{\mathbb R}}}

\newcommand{\MISE}{\mbox{MISE}}
\newcommand{\Span}{\mbox{Span}}
\newcommand{\intg}{\mbox{int}}
\newcommand{\card}{\mbox{card}}
\newcommand{\Range}{\mbox{Range}}
\newcommand{\Var}{\mbox{Var}}
\newcommand{\Cov}{\mbox{Cov}}
\newcommand{\diag}{\mbox{diag}}
\newcommand{\supp}{\mbox{supp}}
\newcommand{\xil}{ {\it et. al }}
\newcommand{\std}{\mbox{std}}
\newcommand{\SNR}{\mbox{SNR}}
\newcommand{\Tr}{\mbox{Tr}}
\newcommand{\proj}{\mbox{proj}}
\newcommand{\Ber}{\mbox{Ber}}
\newcommand{\ERR}{\mbox{ERR}}
\newcommand{\rank}{\mbox{rank}}
\newcommand{\dist}{\mbox{dist}}
\newcommand{\vect}{\mbox{vec}}
\newcommand{\Pen}{\mbox{Pen}}
\newcommand{\Err}{\mbox{Err}}

\newtheorem{thm}{Theorem}
\newtheorem{lem}{Lemma}
\newtheorem{cor}{Corollary}
\newtheorem{rem}{Remark}

\newcommand{\bo}{\mathbf{}}
\newcommand{\ba}{\mathbf{a}}
\newcommand{\bb}{\mathbf{b}}
\newcommand{\bc}{\mathbf{c}}
\newcommand{\bd}{\mathbf{d}}
\newcommand{\boe}{\mathbf{e}}
\newcommand{\bof}{\mathbf{f}}
\newcommand{\bh}{\mathbf{h}}
\newcommand{\bi}{\mathbf{i}}
\newcommand{\bq}{\mathbf{q}}
\newcommand{\bt}{\mathbf{t}}
\newcommand{\bu}{\mathbf{u}}
\newcommand{\bv}{\mathbf{v}}
\newcommand{\bw}{\mathbf{w}}
\newcommand{\bx}{\mathbf{x}}
\newcommand{\by}{\mathbf{y}}
\newcommand{\bz}{\mathbf{z}}

\newcommand{\bA}{\mathbf{A}}
\newcommand{\bB}{\mathbf{B}}
\newcommand{\bD}{\mathbf{D}}
\newcommand{\bI}{\mathbf{I}}
\newcommand{\bM}{\mathbf{M}}
\newcommand{\bQ}{\mathbf{Q}}
\newcommand{\bR}{\mathbf{R}}

\newcommand{\bH}{\mathbf{H}}
\newcommand{\bU}{\mathbf{U}}
\newcommand{\bV}{\mathbf{V}}
\newcommand{\bW}{\mathbf{W}}
\newcommand{\bX}{\mathbf{X}}
\newcommand{\bY}{\mathbf{Y}}
\newcommand{\bZ}{\mathbf{Z}}

\newcommand{\bzero}{\mathbf{0}}

\newcommand{\bte}{\mbox{\mathversion{bold}$\te$}}
\newcommand{\bupsilon}{\mbox{\mathversion{bold}$\upsilon$}}
 \newcommand{\btheta}{\mbox{\mathversion{bold}$\theta$}}
\newcommand{\bxi}{\mbox{\mathversion{bold}$\xi$}}
\newcommand{\boeta}{\mbox{\mathversion{bold}$\xi$}}
\newcommand{\bbe}{\mbox{\mathversion{bold}$\beta$}}
\newcommand{\bzeta}{\mbox{\mathversion{bold}$\zeta$}}
\newcommand{\bphi}{\mbox{\mathversion{bold}$\ph$}}
\newcommand{\bpsi}{\mbox{\mathversion{bold}$\psi$}}
\newcommand{\beps}{\mbox{\mathversion{bold}$\eps$}}
\newcommand{\bobeta}{\mbox{\mathversion{bold}$\beta$}}
\newcommand{\bgamma}{\mbox{\mathversion{bold}$\gamma$}}
\newcommand{\bmu}{\mbox{\mathversion{bold}$\mu$}}

\newcommand{\hbtheta}{\widehat{\btheta}}

\newcommand{\bbJ}{(\bb_J)}

\newcommand{\bPhi}{\mbox{\mathversion{bold}$\Phi$}}
\newcommand{\bUp}{\mbox{\mathversion{bold}$\Up$}}
\newcommand{\bPsi}{\mbox{\mathversion{bold}$\Psi$}}
\newcommand{\bLam}{\mbox{\mathversion{bold}$\Lambda$}}
\newcommand{\bOm}{\mbox{\mathversion{bold}$\Om$}}

\newcommand{\Jc}{{J_{*}^c}}

\newcommand{\calD}{{\mathcal{D}}}
\newcommand{\calG}{{\mathcal G}}
\newcommand{\calJ}{{\mathcal{J}}}
\newcommand{\calmu}{{\mathcal M}_{\mu}}
\newcommand{\calL}{{\mathcal{L}}}
\newcommand{\calM}{{\mathcal{M}}}
\newcommand{\calP}{{\mathcal{P}}}
\newcommand{\calS}{{\mathcal{S}}}
\newcommand{\calW}{{\mathcal{W}}}
\newcommand{\calX}{{\cal{X}}}
\newcommand{\calY}{{\cal{Y}}}
\newcommand{\calZ}{{\cal{Z}}}

\newcommand{\calH}{{\cal H}}
\newcommand{\calT}{{\cal T}}
\newcommand{\calN}{{\cal N}}
\newcommand{\calF}{{\mathcal F}}

\newcommand{\sumjni}{\displaystyle\sum_{j\neq i}}
\newcommand{\sumjn}{\displaystyle\sum_{j=1}^n}
\newcommand{\sumin}{\displaystyle\sum_{i=1}^n}
\newcommand{\sumJ}{\sum_{j \in J_{*}}}
\newcommand{\sumJc}{\sum_{j \in \Jc}}
\newcommand{\sumijn}{\displaystyle\sum_{i,j =1}^n}

\newcommand{\scrP}{\mathscr{P}}
\newcommand{\scrPZK}{\mathscr{P}_{Z,K}}

\newcommand{\tilXi}{\tilde{\Xi}}
\newcommand{\tilP}{\tilde{P}}
\newcommand{\tilA}{\tilde{A}}
\newcommand{\tilC}{\tilde{C}}

\long\def\ignore#1{}

\newcommand {\colred}[1] {\textcolor{red}{#1}}
\newcommand {\colblue}[1] {\textcolor{blue}{#1}}
\newcommand {\colyellow}[1] {\textcolor{yellow}{#1}}
\newcommand {\colblack}[1] {\textcolor{black}{#1}}
\newcommand {\colgreen}[1] {\textcolor{OliveGreen}{#1}}
\newcommand {\colsred}[1] {\textcolor{OrangeRed}{#1}}
\newcommand {\colcn}[1] {\textcolor{cyan}{#1}}

%% file: PABM_ArXiv.bbl
\begin{thebibliography}{10}

\bibitem{JMLR:v18:16-480}
E.~Abbe.
\newblock Community detection and stochastic block models: Recent developments.
\newblock {\em J. Mach. Learn. Res.}, 18(177):1--86, 2018.

\bibitem{P.AgarwalandN.Mustafa2004}
P.~K. Agarwal and N.~H. Mustafa.
\newblock K-means projective clustering.
\newblock In {\em Proceedings of the twenty-third ACM SIGMOD-SIGACT-SIGART
  symposium on Principles of database systems}, pages 155--165. ACM, 2004.

\bibitem{Airoldi:2008:MMS:1390681.1442798}
E.~M. Airoldi, D.~M. Blei, S.~E. Fienberg, and E.~P. Xing.
\newblock Mixed membership stochastic blockmodels.
\newblock {\em J. Mach. Learn. Res.}, 9:1981--2014, June 2008.

\bibitem{DBLP:journals/corr/AminiL14}
A.~A. Amini and E.~Levina.
\newblock On semidefinite relaxations for the block model.
\newblock {\em Ann. Statist.}, 46(1):149--179, 02 2018.

\bibitem{bandeira2016}
A.~S. Bandeira and R.~van Handel.
\newblock Sharp nonasymptotic bounds on the norm of random matrices with
  independent entries.
\newblock {\em Ann. Probab.}, 44(4):2479--2506, 07 2016.

\bibitem{benaychgeorges2017spectral}
F.~Benaych-Georges, C.~Bordenave, and A.~Knowles.
\newblock Spectral radii of sparse random matrices, 2017.

\bibitem{Bickel21068}
P.~J. Bickel and A.~Chen.
\newblock A nonparametric view of network models and newman{\textendash}girvan
  and other modularities.
\newblock {\em Proceedings of the National Academy of Sciences},
  106(50):21068--21073, 2009.

\bibitem{Boucheron2013}
S.~Boucheron, G.~Lugosi, and P.~Massart.
\newblock {\em Concentration inequalities: A nonasymptotic theory of
  independence}.
\newblock Oxford university press, 2013.

\bibitem{inproceedings}
T.~Boult and L.~Gottesfeld~Brown.
\newblock Factorization-based segmentation of motions.
\newblock pages 179 -- 186, 11 1991.

\bibitem{Bradley:2000:KC:596077.596262}
P.~S. Bradley and O.~L. Mangasarian.
\newblock k-plane clustering.
\newblock {\em J. of Global Optimization}, 16(1):23--32, Jan. 2000.

\bibitem{celisse2012}
A.~Celisse, J.-J. Daudin, and L.~Pierre.
\newblock Consistency of maximum-likelihood and variational estimators in the
  stochastic block model.
\newblock {\em Electron. J. Statist.}, 6:1847--1899, 2012.

\bibitem{chen2018}
Y.~Chen, X.~Li, and J.~Xu.
\newblock Convexified modularity maximization for degree-corrected stochastic
  block models.
\newblock {\em Ann. Statist.}, 46(4):1573--1602, 08 2018.

\bibitem{doi:10.1080/10618600.2016.1237362}
J.~Cheng, T.~Li, E.~Levina, and J.~Zhu.
\newblock High-dimensional mixed graphical models.
\newblock {\em Journal of Computational and Graphical Statistics},
  26(2):367--378, 2017.

\bibitem{crossley2013cognitive}
N.~A. Crossley, A.~Mechelli, P.~E. V{\'e}rtes, T.~T. Winton-Brown, A.~X. Patel,
  C.~E. Ginestet, P.~McGuire, and E.~T. Bullmore.
\newblock Cognitive relevance of the community structure of the human brain
  functional coactivation network.
\newblock volume 110, pages 11583--11588. National Acad Sciences, 2013.

\bibitem{Vidal:2009aa}
E.~{Elhamifar} and R.~{Vidal}.
\newblock Sparse subspace clustering.
\newblock In {\em 2009 IEEE Conference on Computer Vision and Pattern
  Recognition}, pages 2790--2797, June 2009.

\bibitem{Elhamifar:2013:SSC:2554063.2554078}
E.~Elhamifar and R.~Vidal.
\newblock Sparse subspace clustering: Algorithm, theory, and applications.
\newblock {\em IEEE Trans. Pattern Anal. Mach. Intell.}, 35(11):2765--2781,
  Nov. 2013.

\bibitem{Favaro:2011:CFS:2191740.2191857}
P.~Favaro, R.~Vidal, and A.~Ravichandran.
\newblock A closed form solution to robust subspace estimation and clustering.
\newblock CVPR '11, pages 1801--1807, Washington, DC, USA, 2011. IEEE Computer
  Society.

\bibitem{Gao:2017:AOM:3122009.3153016}
C.~Gao, Z.~Ma, A.~Y. Zhang, and H.~H. Zhou.
\newblock Achieving optimal misclassification proportion in stochastic block
  models.
\newblock {\em J. Mach. Learn. Res.}, 18(1):1980--2024, Jan. 2017.

\bibitem{gao2018community}
C.~Gao, Z.~Ma, A.~Y. Zhang, H.~H. Zhou, et~al.
\newblock Community detection in degree-corrected block models.
\newblock {\em The Annals of Statistics}, 46(5):2153--2185, 2018.

\bibitem{Giraud:1999833}
C.~Giraud.
\newblock {\em {Introduction to high-dimensional statistics}}.
\newblock Chapman \& Hall/CRC Monographs on Statistics \& Applied Probability.
  CRC Press, Hoboken, NJ, 2015.

\bibitem{MAL-005}
A.~Goldenberg, A.~X. Zheng, S.~E. Fienberg, and E.~M. Airoldi.
\newblock A survey of statistical network models.
\newblock {\em Foundations and Trends{\textregistered} in Machine Learning},
  2(2):129--233, 2010.

\bibitem{gupta_nagar}
A.~K. Gupta and D.~K. Nagar.
\newblock {\em Matrix Variate Distributions}.
\newblock Chapman \& Hall/CRC, 1999.

\bibitem{Hsu2011}
D.~Hsu, S.~Kakade, and T.~Zhang.
\newblock A tail inequality for quadratic forms of subgaussian random vectors.
\newblock {\em Electron. Commun. Probab.}, 17:6 pp., 2012.

\bibitem{2017arXiv170807852J}
J.~{Jin}, Z.~T. {Ke}, and S.~{Luo}.
\newblock {Estimating network memberships by simplex vertex hunting}.
\newblock {\em arXiv e-prints}, page arXiv:1708.07852, Aug. 2017.

\bibitem{joseph2016}
A.~Joseph and B.~Yu.
\newblock Impact of regularization on spectral clustering.
\newblock {\em Ann. Statist.}, 44(4):1765--1791, 08 2016.

\bibitem{Karrer2011StochasticBA}
B.~Karrer and M.~E.~J. Newman.
\newblock Stochastic blockmodels and community structure in networks.
\newblock {\em Physical review. E, Statistical, nonlinear, and soft matter
  physics}, 83 1 Pt 2:016107, 2011.

\bibitem{KlTsVe2017}
O.~Klopp, A.~B. Tsybakov, and N.~Verzelen.
\newblock Oracle inequalities for network models and sparse graphon estimation.
\newblock {\em Ann. Statist.}, 45(1):316--354, 2017.

\bibitem{Kolaczyk:2009:SAN:1593430}
E.~D. Kolaczyk.
\newblock {\em Statistical Analysis of Network Data: Methods and Models}.
\newblock Springer Publishing Company, 1st edition, 2009.

\bibitem{le2016}
C.~M. Le, E.~Levina, and R.~Vershynin.
\newblock Optimization via low-rank approximation for community detection in
  networks.
\newblock {\em Ann. Statist.}, 44(1):373--400, 02 2016.

\bibitem{lei2015}
J.~Lei and A.~Rinaldo.
\newblock Consistency of spectral clustering in stochastic block models.
\newblock {\em Ann. Statist.}, 43(1):215--237, 02 2015.

\bibitem{Liu:2013:RRS:2412386.2412936}
G.~Liu, Z.~Lin, S.~Yan, J.~Sun, Y.~Yu, and Y.~Ma.
\newblock Robust recovery of subspace structures by low-rank representation.
\newblock {\em IEEE Trans. Pattern Anal. Mach. Intell.}, 35(1):171--184, Jan.
  2013.

\bibitem{Liu2010RobustSS}
G.~Liu, Z.~Lin, and Y.~Yu.
\newblock Robust subspace segmentation by low-rank representation.
\newblock In {\em Proceedings of the 27th International Conference on
  International Conference on Machine Learning}, ICML'10, pages 663--670, USA,
  2010. Omnipress.

\bibitem{Ma:2008:ESA:1405158.1405160}
Y.~Ma, A.~Y. Yang, H.~Derksen, and R.~Fossum.
\newblock Estimation of subspace arrangements with applications in modeling and
  segmenting mixed data.
\newblock {\em SIAM Rev.}, 50(3):413--458, Aug. 2008.

\bibitem{mairal2014spams}
J.~Mairal, F.~Bach, J.~Ponce, G.~Sapiro, R.~Jenatton, and G.~Obozinski.
\newblock Spams: A sparse modeling software, v2.3.
\newblock {\em URL http://spams-devel. gforge. inria. fr/downloads. html},
  2014.

\bibitem{Mallat:1993:MPT:2198030.2203996}
S.~Mallat and Z.~Zhang.
\newblock Matching pursuits with time-frequency dictionaries.
\newblock {\em Trans. Sig. Proc.}, 41(12):3397--3415, Dec. 1993.

\bibitem{ICPR_2011_Nasihatkon}
B.~Nasihatkon and R.~Hartley.
\newblock Graph connectivity in sparse subspace clustering.
\newblock CVPR '11, pages 2137--2144. IEEE Computer Society, 06 2011.

\bibitem{noroozi2019sparse}
M.~Noroozi, R.~Rimal, and M.~Pensky.
\newblock Sparse popularity adjusted stochastic block model, 2019.

\bibitem{rao_rao_1998}
C.~Rao and M.~Rao.
\newblock {\em Matrix Algebra and Its Applications to Statistics and
  Econometrics}, volume 528.
\newblock World Scientific, 1998.

\bibitem{rohe2011spectral}
K.~Rohe, S.~Chatterjee, B.~Yu, et~al.
\newblock Spectral clustering and the high-dimensional stochastic blockmodel.
\newblock {\em Ann. Statist.}, 39(4):1878--1915, 2011.

\bibitem{RePEc:bla:jorssb:v:80:y:2018:i:2:p:365-386}
S.~Sengupta and Y.~Chen.
\newblock A block model for node popularity in networks with community
  structure.
\newblock {\em Journal of the Royal Statistical Society Series B},
  80(2):365--386, 2018.

\bibitem{Shi_2019}
B.~Shi and S.~Iyengar.
\newblock {\em Mathematical Theories of Machine Learning - Theory and
  Applications}.
\newblock Springer, 2019.

\bibitem{soltanolkotabi2012}
M.~Soltanolkotabi and E.~J. Candes.
\newblock A geometric analysis of subspace clustering with outliers.
\newblock {\em Ann. Statist.}, 40(4):2195--2238, 08 2012.

\bibitem{soltanolkotabi2014}
M.~Soltanolkotabi, E.~Elhamifar, and E.~J. Candes.
\newblock Robust subspace clustering.
\newblock {\em Ann. Statist.}, 42(2):669--699, 04 2014.

\bibitem{tseng2000nearest}
P.~Tseng.
\newblock Nearest q-flat to m points.
\newblock {\em Journal of Optimization Theory and Applications},
  105(1):249--252, 2000.

\bibitem{vershynin_2012}
R.~Vershynin.
\newblock {\em Introduction to the non-asymptotic analysis of random matrices},
  pages 210--268.
\newblock Cambridge University Press, 2012.

\bibitem{vidal2011subspace}
R.~Vidal.
\newblock Subspace clustering.
\newblock {\em IEEE Signal Processing Magazine}, 28(2):52--68, 2011.

\bibitem{vidal2005generalized}
R.~Vidal, Y.~Ma, and S.~Sastry.
\newblock Generalized principal component analysis (gpca).
\newblock {\em IEEE Trans. Pattern Anal. Mach. Intell.}, 27(12):1945--1959,
  2005.

\bibitem{wang2018network}
B.~Wang, A.~Pourshafeie, M.~Zitnik, J.~Zhu, C.~D. Bustamante, S.~Batzoglou, and
  J.~Leskovec.
\newblock Network enhancement as a general method to denoise weighted
  biological networks.
\newblock {\em Nature Communications}, 9(1):3108, 2018.

\bibitem{pmlr-v51-wang16b}
Y.~Wang, Y.-X. Wang, and A.~Singh.
\newblock Graph connectivity in noisy sparse subspace clustering.
\newblock In A.~Gretton and C.~C. Robert, editors, {\em Proceedings of the 19th
  International Conference on Artificial Intelligence and Statistics},
  volume~51 of {\em Proceedings of Machine Learning Research}, pages 538--546,
  Cadiz, Spain, 09--11 May 2016. PMLR.

\bibitem{JMLR:v17:13-354}
Y.-X. Wang and H.~Xu.
\newblock Noisy sparse subspace clustering.
\newblock {\em Journal of Machine Learning Research}, 17(12):1--41, 2016.

\bibitem{weisberg2005applied}
S.~Weisberg.
\newblock {\em Applied linear regression}, volume 528.
\newblock John Wiley \& Sons, 2005.

\bibitem{zhao2012consistency}
Y.~Zhao, E.~Levina, J.~Zhu, et~al.
\newblock Consistency of community detection in networks under degree-corrected
  stochastic block models.
\newblock {\em Ann. Statist.}, 40(4):2266--2292, 2012.

\end{thebibliography}
